\documentclass[11pt]{article}
\usepackage{latexsym}
\usepackage{amssymb, amsmath, xspace, lscape,  latexsym, color}
\usepackage{graphicx} 
\graphicspath{{figures/}}

\renewcommand{\ge}{\geqslant}

\newcommand{\bea}{\begin{eqnarray}}
\newcommand{\eea}{\end{eqnarray}}
\def\beq#1#2\eeq{
        \begin{equation}
        \label{#1}
            #2
        \end{equation}}

\newcommand{\R}{\mathbb R}

\newcommand{\C}{\mathbb C}

\newcommand{\Z}{\mathbb Z}

\renewcommand{\hat}{\widehat}
\renewcommand{\tilde}{\widetilde}

\def\btheor#1\etheor{
        \begin{theor}
            #1
        \end{theor}
    }

    \def\bsled#1\esled{
        \begin{sled}
            #1
        \end{sled}   }

\def\btheor#1\elemma{
        \begin{lemma}
            #1
        \end{lemma}
    }

    \def\bsled#1\esled{
        \begin{sled}
            #1
        \end{sled}   }

\newtheorem{theorem}{Theorem}
\newtheorem{definit}{Definition}
\newtheorem{lemma}{Lemma}
\newtheorem{prop}{Proposition}

\newtheorem{remark}{Remark}

\def\hm#1{#1\nobreak\discretionary{}{\hbox{\m@th$#1$}}{}}
\def\mi#1{\discretionary{\hbox{\m@th$#1$}}{\hbox{\m@th$#1$}}{}}

\vspace{4ex}
\begin{document}
%\small{
\title{%On the discrete $Z^{{a}}$ }
The asymptotic behaviour of the discrete holomorphic map $Z^a$ via the Riemann-Hilbert method}
%}
\author{Alexander I. Bobenko\\
        Institut f\"r Mathematik\\
        Technische Universit\"at Berlin,\\
        Strasse des 17. June 136,\\
        10623 Berlin, Germany\\
        {}\\
        Alexander Its\\
        Department of Mathematical Sciences\\
        Indiana University-Purdue University Indianapolis\\
        402 N. Blackford Street\\
        Indianapolis, IN 46202-3216
        USA}

%\date{\textsf{09-05-2014}}
\maketitle
%%\begin{Summary of research}
\centerline{\bf{Abstract}} We study the asymptotic behavior of the discrete analogue of the 
holomorphic map $z^{{a}}$. The analysis is based on the use of the Riemann-Hilbert approach.
Specifically, using the Deift-Zhou nonlinear steepest descent method we prove
the asymptotic formulae which was conjectured in 2000 by the first co-author
and S.I.~Agafonov.

%%\end{Summary of research}
\vfill\eject

\noindent

\setcounter{equation}{0}
\section{Introduction.}

The nonlinear theory of discrete complex analysis goes back to 1985 Thurston's talk \cite{Thu}
at Purdue University  and declares circle patterns to be natural discrete analogs of analytic functions \cite{Sch, Ste}. The word ``nonlinear'' refers to the basic feature of equations describing circle patterns. Often, the so-called cross-ratio system is used for this. In \cite{BP} a {\em discrete conformal} map was defined as a complex valued function on the square grid $f: \Z^{2} \to \R^2 = \C$ with the property that the cross ratio on each elementary quadrilateral is -1:
\begin{equation}\label{def1}
\frac{(f_{n,m} - f_{n+1,m})(f_{n+1,m+1} - f_{n,m+1})}
{(f_{n+1,m} - f_{n+1,m+1})(f_{n,m+1} - f_{n,m})} = -1.
\end{equation}
Here and below we abbreviate $f_{n,m}=f(n,m)$. The boundary data $f(n,0), f(0,m)$ and the evolution equation (\ref{def1}) determine the whole map uniquely. A discrete conformal map is called {\em embedded} if the interiors of different elementary quadrilaterals are disjoint.
 
Note that the definition of a discrete conformal map is M\"obius invariant and is motivated by the following characterization for smooth mappings: A smooth map $f:D\to{\mathbb C}$ is conformal (holomorphic or antiholomorphic) if and only if $\forall z\in D\subset{\mathbb C}$
$$
\lim_{\epsilon\to 0}\frac{(f(z) - f(z+\epsilon)(f(z+\epsilon+i\epsilon) - f(z+i\epsilon))}
{(f(z+\epsilon) - f(z+\epsilon+i\epsilon))(f(z+i\epsilon) - f(z))} = -1.
$$

\begin{figure}[t]
  \begin{center}
\includegraphics[width=0.49\linewidth]{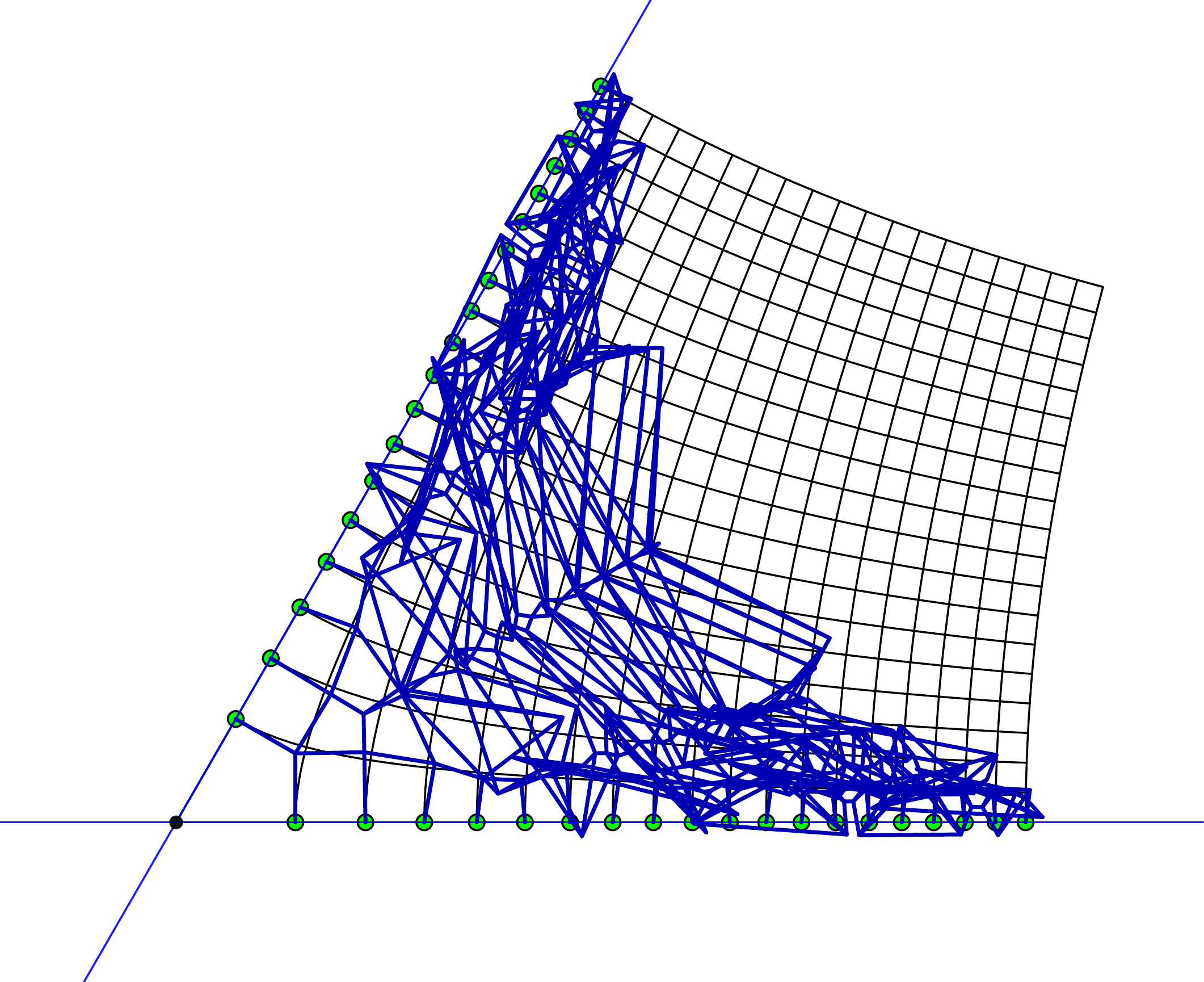} 
\includegraphics[width=0.49\linewidth]{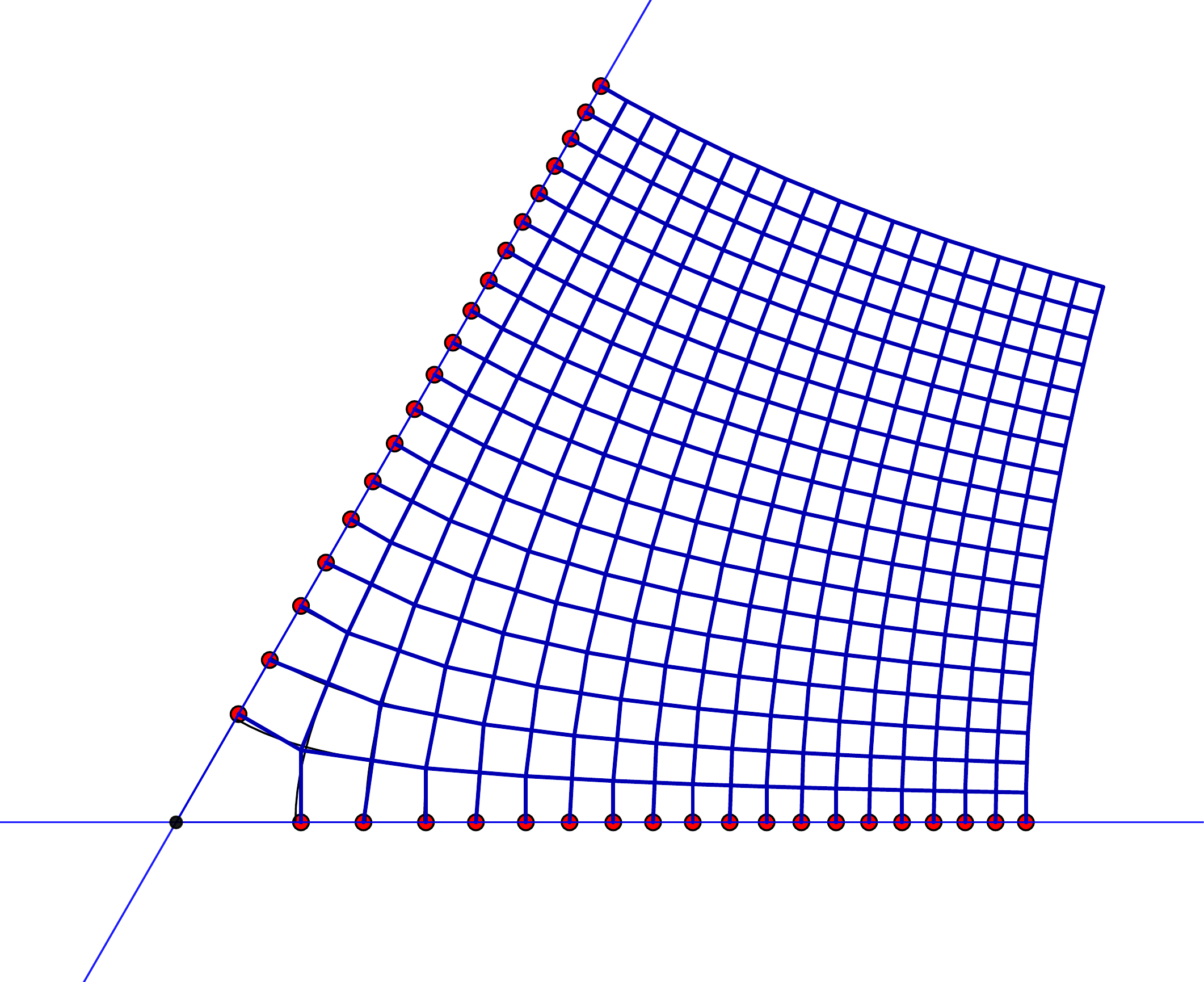}
  \end{center}
  \caption{Two discrete conformal maps with close initial data $n=0, m=0$: 
  (Left) Continuous holomorphic mapping $z^{2/3}$ and the discrete holomorphic mapping with the ``naive'' boundary data $f(n,0)=n^{2/3}$ and $f(0,m)=(im)^{2/3}$. 
  (Right) The discrete $Z^{2/3}$: the boundary data is slightly different from the ``naive'' ones. {\small [Images by J. Richter-Gebert]}}\label{f.smooth-discrete}
\end{figure}

It is a very appealing problem to find discrete conformal maps corresponding to classical holomorphic functions. In the following we discuss a discretization of the holomorphic map $z^a$. A naive way to construct a discrete analogue of the holomorphic map $z^a$  would be to take (\ref{def1}) with the boundary data $f(n,0)=n^{2/3}$ and $f(0,m)=(im)^{2/3}$.  However, as demonstrated in Figure~\ref{f.smooth-discrete}(left), the resulting lattice is not embedded and is far from its continuous counterpart. Hence this map cannot be treated as a discrete $z^a$.

The discrete embedded analog, $Z^a$, of the function $z^a$ exists and is shown in Figure~\ref{f.smooth-discrete}(right).
To construct discrete $Z^a$ more involved methods coming from the theory of integrable systems are required. Indeed, a crucial property of equation (\ref{def1}) is its integrability \cite{NC,BP}. In the following we summarize some known facts about the discrete conformal map $Z^a$, see \cite{BobSurDDG, AB, Bob1} for more details. 

The discrete map $Z^a$ was introduced in \cite{Bob1}. In order to construct an embedded discrete analog of $z^a$ the following approach is used. Equation (\ref{def1}) can be supplemented with the nonautonomous constraint
\begin{equation}\label{def2}
a f_{n,m} =2n\frac{(f_{n+1,m} - f_{n,m})(f_{n,m} - f_{n-1,m})}
{(f_{n+1,m} - f_{n-1,m})}
+ 
2m\frac{(f_{n,m+1} - f_{n,m})(f_{n,m} - f_{n,m-1})}
{(f_{n,m+1} - f_{n,m-1})}.
\end{equation}
This constraint is derived within the theory of integrable systems. Solutions of (\ref{def1}) satisfying (\ref{def2}) are singled out by an auxiliary special Fuchsian system, which yields formula (\ref{def2}) (see Section \ref{s.RH} and \cite{Bob1, AB} for more details). This constraint is compatible with (\ref{def1}). A proof of the compatibility based on the analysis of the corresponding Lax representation and the above mentioned Fuchsian system is given in \cite{BobSurDDG}. 

We assume that $0<a<2$ and denote $\Z^{2}_{+}=\{ (n,m)\in\Z : n,m\ge 0 \}$.
To demonstrate that the constraint (\ref{def2}) indeed corresponds to a discrete $Z^a$
we investigate its continuous limit. The right hand side of (\ref{def2}) in the limit $\epsilon\to 0$ for $z=x+iy=\epsilon (n+im)$ gives
$$
\lim_{\epsilon\to 0}2\frac{x}{\epsilon}\frac{(f(z+\epsilon)-f(z))(f(z)-f(z-\epsilon))}{f(z+\epsilon)-f(z-\epsilon)}+
2\frac{y}{\epsilon}\frac{(f(z+i\epsilon)-f(z))(f(z)-f(z-i\epsilon))}{f(z+i\epsilon)-f(z-i\epsilon)}=xf_x+yf_y=zf_z,
$$
where we have used the holomorphicity of the limiting mapping. The corresponding limit of (\ref{def2}) becomes $af=zf_z$, and its general solution is $f(z)=z^a$ up to scaling.  

This consideration and the properties $z^a(\R_+)=\R_+$ and $z^a(i\R_+)=e^{a\pi i/2}\R_+$ of the holomorphic mapping $z^a$ motivate the following definition \cite{Bob1} of its discrete analog.
\begin{definit}
For $0<a<2$ the dicsrete conformal map $Z^a:\Z_+^2\to\C$ is the solution of equations (\ref{def1}) and (\ref{def2}) with the initial conditions 
\begin{equation}\label{incond}
Z^a(0,0) = 0,\quad Z^a(1,0) = 1, \quad Z^a(0,1) = e^{a \pi i/2}.
\end{equation}
\end{definit}
The properties $Z^a(n,0)\in\R_+$ and $Z^a(0,m)\in e^{a \pi i/2}\R_+$ are obvious. 
The existence of this map was proven using the methods of the theory of integrable systems.

\begin{figure}
  \begin{center}
\includegraphics[width=0.5\linewidth]{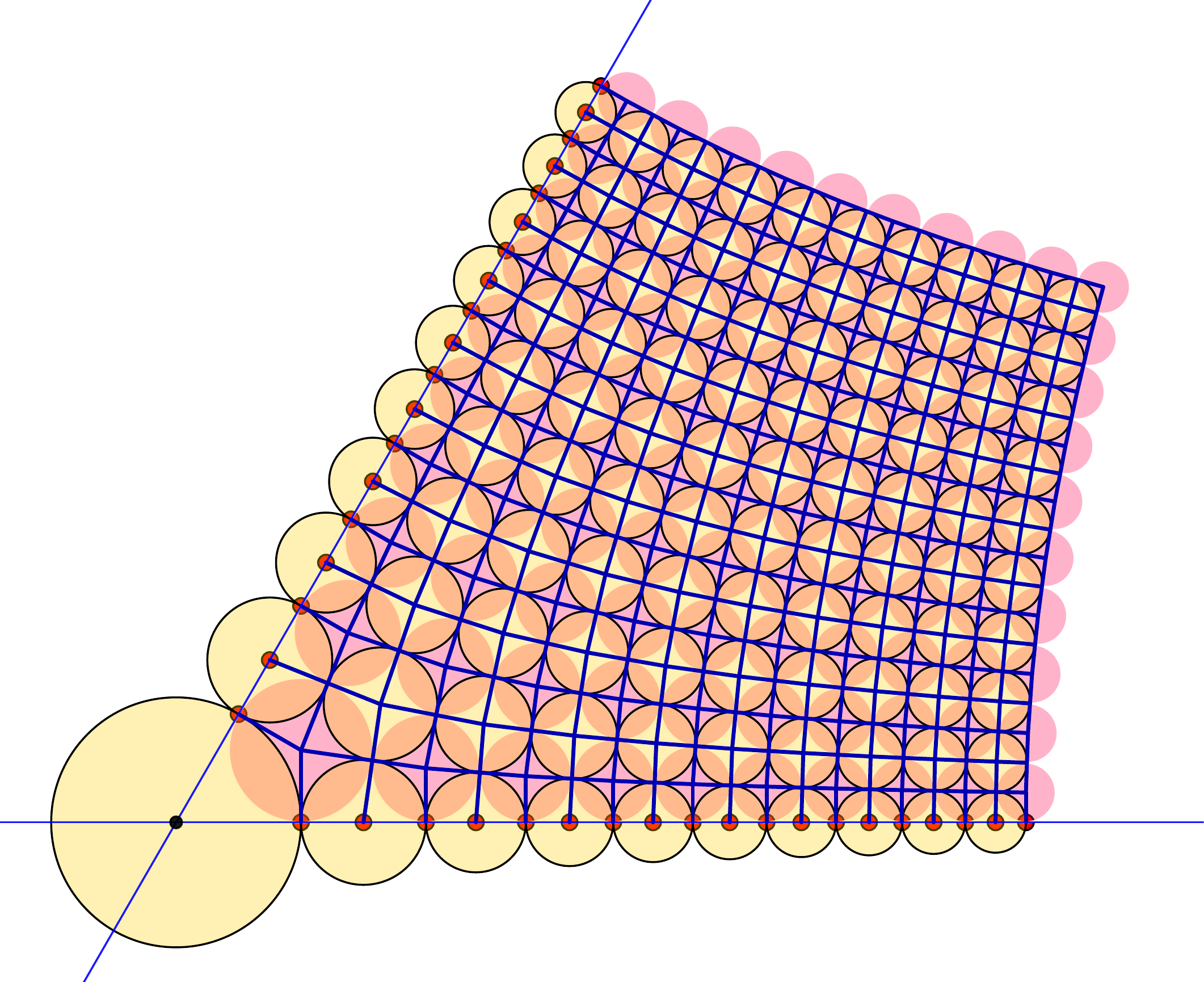} 
\includegraphics[width=0.45\linewidth]{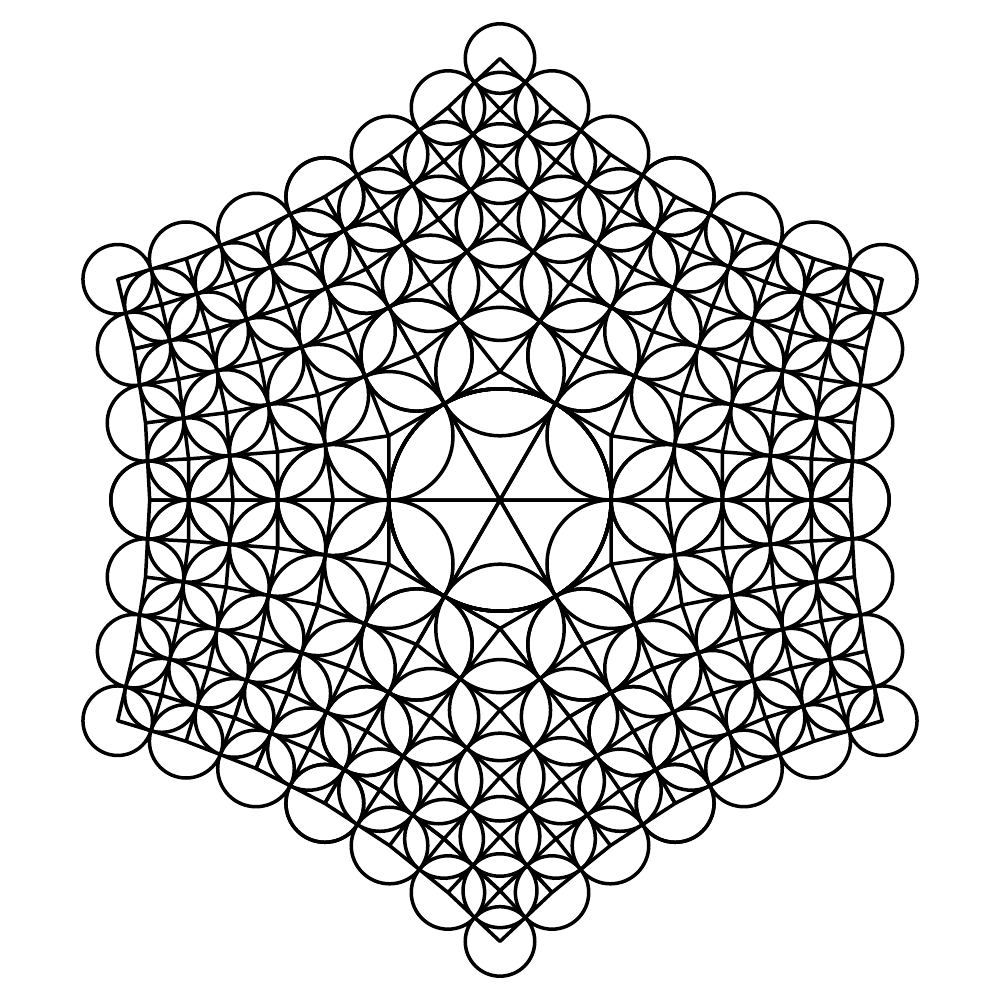}
  \end{center}
  \caption{The discrete $z^{2/3}$ as an infinite orthogonal circle pattern: one sector and the covering of the plane. {\small [Images by J. Richter-Gebert and T. Hoffmann]}}\label{f.circle_pattern}
\end{figure}

As it was shown in \cite{AB}, the discrete conformal map $Z^a$ determines a circle pattern of Schramm type, i.e. an orthogonal circle pattern with the combinatorics of the square grid. The points $Z^a(n,m)$ with even and odd $n+m$ are the centers of the circles and their intersection points respectively (see Figure~\ref{f.circle_pattern}). Moreover, this discrete conformal mapping was also proven to be immersed, i.e. the neighboring elementary quadrilaterals do not overlap. Finally the embeddedness of this mapping was proven in \cite{Aga}.

It turns out that the orthogonal $Z^a$-circle pattern can be defined in a pure geometric way without referring to integrable equations.  The corresponding rigidity result was obtained in \cite{Bue} by analysis methods. It reads as follows. For $a\in (0,2)$ the infinite orthogonal circle pattern corresponding to the discrete conformal mapping $Z^a$ is the unique embedded orthogonal circle pattern (up to global scaling) with the following two properties (see Figure~\ref{f.circle_pattern}(left)):
\begin{itemize}
\item[(i)] The union of the corresponding kites (elementary quadrilaterals) of the $Z^a$-circle pattern covers the infinite sector $\{z=re^{i\phi}\in\C :r\ge 0, \phi\in [0,a\pi/2]\}$ with angle $a\pi/2$.
\item[(ii)] The centers of the boundary circles lie on the boundary half lines $\R_+$ and $e^{a\pi i/2}\R_+$.
\end{itemize}

For rational $a=\frac{4}{N},\ N\in\{2,3\ldots\}$ the rigidity of $Z^a$ follows from the rigidity results obtained in \cite{He}. For example for the infinite circle pattern in Figure~\ref{f.circle_pattern}(right) it reads as follows. Consider an orthogonal circle pattern with the combinatorics shown in this figure, i.e. there is one circle intersected by six neighboring circles and all other circles have exactly four intersecting neighbors. Then an orthogonal embedded circle pattern that covers the whole plane and possesses the described combinatorics is unique.  

Our goal is to prove the following asymptotic behavior of $f_{n,m} \equiv Z^a_{n,m}$
as $n, m \to \infty$, which was conjectured in \cite{AB}. 
\begin{theorem}\label{theorem1}
Let $Z^a_{n,m}$ be the above defined  discrete analog of the power function $z^a$.
Assume that $0 < a < 2$. Then, 
\begin{equation}\label{ABconj}
Z^a_{n,m} = c(a)\left(\frac{n+im}{2}\right)^a
 \left( 1 + O\left(\frac{1}{n^2 + m^2}\right)\right), \quad n^2 + m^2 \to \infty,
\end{equation}
with 
$$
c(a) = \frac{\Gamma\left(1 - \frac{a}{2}\right)}{\Gamma\left(1 + \frac{a}{2}\right)}.
$$
\end{theorem}
%In what follows we shall prove this statement 
%using the Riemann-Hilbert approach. We shall start with putting the problem into the Riemann-Hilbert
%formalism. 

This asymptotics was proven for $n = 0, 1$ in  \cite{AB}. 
Also by elementary methods  the corresponding asymptotics without a formula for $c(a)$ was proven  for $n-m={\rm const}$ in \cite{Aga_rigidity}.

Theorem \ref{theorem1} is the  statement about the asymptotics  of the solution of the Cauchy problem
for equations (\ref{def1}) and (\ref{def2}) determined by the initial data (\ref{incond}).
Equations (\ref{def1}) and (\ref{def2}) are nonlinear  difference equations. It is difficult, if not impossible, 
unless solution is explicit or given in terms of contour integrals,
to perform global asymptotic analysis of nonlinear equations, both difference and differential. 
The  reason we  are able to do this in the case of the  Cauchy
problem for equations (\ref{def1}) and (\ref{def2}) is their integrability. The latter allows us
to use the Riemann-Hilbert approach - a noncomutative analog of contour integral
representation and apply  the nonlinear steepest descent method  of Deift and Zhou \cite{DZ}
in our investigation.  

%Although the integrability of $Z^a_{n,m}$ has been known for many years, the proof of 
%(\ref{ABconj}) came only now. The reason is in the Fuchshian nature of the associated
%Lax pair  (Section \ref{s.RH} for more detail).  At  the time the conjecture (\ref{ABconj})
%was formulated the nonlinear steepest descent method had not been developed enough,
%while the original WKB-based approach to the integrable asymptotics did not
%work well since nobody knows how to develop the WKB scheme in the presence
%of the regular singularities only. It took a decade for the nonlinear steepest descent 
%method to mature enough to become applicable to the asymptotic analysis of 
%$Z^a_{n,m}$

As it will be shown in the main text,  the function $Z^a_{n,m}$ is intimately related
to a certain collection of orthogonal polynomials. Hence the necessity to use 
the orthogonal polynomial version of the Deift-Zhou method \cite{DKMVZ}. 
The Riemann-Hilbert problem corresponding to $Z^a_{n,m}$ is the problem
of a Fuchsian type - the associated 
system of linear differential equations has the regular singular points only. 
Simultaneously, the problem is posed on a half-line. This is a rather 
rare situation which leads to certain peculiarities  in the implementation 
of the nonlinear steepest descent method. 
In particular, support of  the
relevant equilibrium measure coincides with the whole half-line, and  
the so-called ``lenses opening'' is not a local operation. Also,  what is 
usually appear as  a ``global parametrix'',
here becomes a ``local parametrix''  near infinity. One more deviation  from
the standard situation is the need to use at some point (the proof of Theorem \ref{logas}) 
a rather sophisticated error term estimates in the Hankel asymptotic
series. More details 
on the Riemann-Hilbert problem we are working with  are  in the main text.

Since we address the paper to a broad geometric audience we decided to make it self-contained. In our presentations, we included all the details of the nonlinear steepest descent scheme, although some of them are standard to the experts.

%With all its technical specifics, some of the considerations of this paper will be familiar
%to the experts in the Riemann-Hilbert method. Since we address the paper to a broader geometric 
%audience we decided to make it self-contained and included  the details of the nonlinear steepest descent scheme.
%  so that the ``discrete
%geometry''  part of our readers could  find in the text an introduction to a
%powerful Riemann-Hilbert techniques which we sure will become one 
%of the key analytical tools in discrete geometry.

The proof with the use of the Riemann-Hilbert method needs a lot
of preparatory steps which in itself are of considerable interest.
First, we need the Lax-pair formulation, then the setting of
the relevant monodromy data which is followed by its conversion
into the Riemann-Hilbert setting. In the course of these steps we will
reveal the above mentioned connection to the theory of orthogonal polynomials and
the theory of discrete Painlev\'e equations. These connections do not help to prove 
formula (\ref{ABconj}), while the results of
our paper might be of interest in both these theories. With this in mind,  
we make a  detour from our main goal and discuss in  Sections \ref{opol} 
the orthogonal polynomials related to $Z^{a}$.

Finally in Section \ref{s.log} we define two discrete analogs of the logarithm function: 
the function $L(n,m)$ defining an orthogonal circle pattern  (nonlinear theory) and  Green's function $\ell(n,m)$ (linear theory of discrete holomorphicity).
The latter was introduced by Kenyon in \cite{Ken}. We derive their asymptotics at $ r^2 \equiv n^2+m^2\to\infty$ from (\ref{ABconj}):
\begin{eqnarray*}
L(n,m)&=&\log (n+im)+\gamma -\log 2 + O\left(\frac{\log r}{r^2}\right),\\
\ell(n,m)&=&\log \sqrt{n^2+m^2}+\gamma +\log 2 + O\left(\frac{\log r}{r}\right), \quad n+m \ {\rm even}. 
\end{eqnarray*}
 Here $\gamma$ is Euler's constant. The last formula has already been obtained by a different method in \cite{Ken}.

As already been said, we start with putting the problem of investigation
of the discrete conformal map $Z^a$ into the Riemann-Hilbert formalism.

\section{The Riemann-Hilbert representation for $Z^{{a}}$.} \label{s.RH}
\subsection{Isomonodromity of $Z^{{a}}$ }
The possibility to apply the Riemann-Hilbert technique to the  asympotical analysis of $Z^{{a}}$ is based  
on the integrability of the system  (\ref{def1}) - (\ref{def2}). The latter exactly means the
following two facts.
\begin{prop} \label{prop1} (\cite{NC, BP}) 
The nonlinear difference equation (\ref{def1}) is the compatibility
condition of the following system of linear difference equations - the  {\it Lax pair },
\begin{equation}\label {lax}
\Psi_{n+1,m} = U_{n,m}\Psi_{n.m},\quad \Psi_{n,m+1} = V_{n,m}\Psi_{n.m},
\end{equation}
where
\begin{equation}\label{UV}
U_{n,m} \equiv U_{n,m}(\lambda) 
= \begin{pmatrix}1&-u_{n,m}\cr
\frac{\lambda}{u_{n,m}}&1\end{pmatrix},\quad
V_{n,m} \equiv V_{n,m}(\lambda) 
= \begin{pmatrix}1&-v_{n,m}\cr
-\frac{\lambda}{v_{n,m}}&1\end{pmatrix},
\end{equation}
and
\begin{equation}\label{uv}
u_{n,m} = f_{n+1,m} - f_{n,m},\quad v_{n,m} = f_{n,m+1} - f_{n,m}.
\end{equation}
\end{prop}
In particular, this statement means that equation (\ref{def1}) implies the
matrix relation,
\begin{equation}\label{UVVU}
U_{n,m+1}(\lambda)V_{n,m}(\lambda) = V_{n+1,m}(\lambda)U_{n,m}(\lambda),\quad \forall \lambda.
\end{equation}

The following proposition was proven in (\cite{AB}), note that the isomonodromic constraint (\ref{def2}) was obtained for $a=1$ in \cite {Nij}.
\begin{prop}\label{prop2}
The addition  constraint (\ref{def2}) is equivalent to the existence 
of a solution $\Psi_{n,m}$ to (\ref{lax}) satisfying also the 
following ``$\lambda$ - equation'',
\begin{equation}\label{lambda}
 \frac{d}{d\lambda}\Psi_{n,m}  = A_{n,m}\Psi_{n,m},\quad
 A_{n,m} = -\frac{B_{n,m}}{1+\lambda} + \frac{C_{n,m}}{1-\lambda}
 +\frac{D_{n,m}}{\lambda},
 \end{equation}
where the independent of $\lambda$ matrices $B_{n,m}$, $C_{n,m}$,
and $D_{n,m}$ are of the following structure,
\begin{equation}\label{B}
B_{n,m} = -\frac{n}{u_{n,m} + u_{n-1,m}}\begin{pmatrix}u_{n,m}&u_{n,m}u_{n-1,m}\cr
1&u_{n-1,m}\end{pmatrix}
\end{equation}
\begin{equation}\label{C}
C_{n,m} = -\frac{m}{v_{n,m} + v_{n,m-1}}\begin{pmatrix}v_{n,m}&v_{n,m}v_{n,m-1}\cr
1&v_{n,m-1}\end{pmatrix}
\end{equation}
\begin{equation}\label{D}
D_{n,m} = \begin{pmatrix}-\frac{{a}}{4}&-\frac{{a}}{2}f_{n,m}\cr
0&\frac{{a}}{4}\end{pmatrix}
\end{equation}
\end{prop}

In particular, this statement means that the system (\ref{def1}) - (\ref{def2}) implies, in
addition to (\ref{UVVU}), two more matrix equations, 
\begin{equation}\label{UAAU}
\frac{dU_{n,m}(\lambda)}{d\lambda} =A_{n+1,m}(\lambda)U_{n,m}(\lambda)
- U_{n,m}(\lambda)A_{n,m}(\lambda),\quad \forall \lambda,
\end{equation}
and
\begin{equation}\label{VAAV}
\frac{dV_{n,m}(\lambda)}{d\lambda} =A_{n,m+1}(\lambda)V_{n,m}(\lambda)
- V_{n,m}(\lambda)A_{n,m}(\lambda),\quad \forall \lambda.
\end{equation}

Equation (\ref{lambda}) is a Fuchsian liner system with four regular
points ($\pm1$, $0$ and $\infty$).   The above statements imply that  equations
(\ref{def1}) - (\ref{def2}) describe discrete isomonodromy deformations
of system (\ref{lambda}), and that the monodromy data of
this system are the first integrals of (\ref{def1}) - (\ref{def2}) (cf. \cite{JMU}). 
Our first step will be  the  evaluation of these integrals for the particular choice 
of the initial data (\ref{incond}) corresponding to  $Z^{{a}}$.
We shall start with the definition of the matrix valued  function $\Psi_{n,m}(\lambda)$ - a carrier 
of the monodromy data in question, by the equations,
$$
\Psi_{0,0}(\lambda)=\lambda^{-\frac{{a}}{4}\sigma_{3}},\quad \Psi_{0,1}(\lambda) = V_{0,0}(\lambda)\Psi_{0,0}(\lambda),
\quad \Psi_{1,1}(\lambda) = U_{0,1}(\lambda)V_{0,0}(\lambda)\Psi_{0,0}(\lambda),
$$
$$
\Psi_{n,m}(\lambda)=U_{n-1,m}(\lambda)U_{n-2,m}(\lambda)...U_{0,m}(\lambda)
$$
\begin{equation}\label{Psinmdef}
\times V_{0,m-1}(\lambda)V_{0,m-2}(\lambda)...V_{0,0}(\lambda)\Psi_{0,0}(\lambda),
\quad  n, m \geq 1,
\end{equation}
$$
\sigma_3 = \begin{pmatrix}1&0\cr0&-1\end{pmatrix}.
$$
In these equations, $u_{n,m}$ and $v_{n,m}$ are defined via (\ref{uv}) with $f_{n,m} \equiv Z^{{a}}_{n,m}$.
The  function $\lambda^{-\frac{{a}}{4}\sigma_3}$ as defined  on the 
$\lambda$ -plane cut along the negative imaginary axis and fixed by the condition,
$$
-\frac{\pi}{2} < \arg\lambda < \frac{3\pi}{2}.
$$
 It is also worth noticing that
\begin{equation}\label{detPsi}
\det\Psi_{n,m}(\lambda) = (\lambda+1)^{n}(1-\lambda)^{m}.
\end{equation} 

\begin{prop} The function $\Psi_{n,m}(\lambda)$ is the common solution of linear equations (\ref{lax}),
and (\ref{lambda}).
\end{prop}
{\it Proof.} The first equation in (\ref{lax}) is satisfied by construction. In order to see that 
the second equation in (\ref{lax}) is satisfied it is enough to observe that 
matrix equation (\ref{UVVU}) allows to switch the matrices $U$ and $V$ in
the definition of the function $\Psi_{n,m}(\lambda)$ and  re-write it in the form,
$$
\Psi_{n,m}(\lambda)=V_{n,m-1}(\lambda)V_{n,m-2}(\lambda)...V_{n,0}(\lambda)
$$
\begin{equation}\label{Psinmdef2}
\times U_{n-1,0}(\lambda)U_{n-2,0}(\lambda)...U_{0,0}(\lambda)\Psi_{0,0}(\lambda),
\quad  n, m \geq 1.
\end{equation}
Verification of equation (\ref{lambda}) needs a little bit more work. Put
$$
F_{n,m}(\lambda) :=  \frac{d}{d\lambda}\Psi_{n,m}(\lambda)  - A_{n,m}(\lambda)\Psi_{n,m}(\lambda).
$$
In view of (\ref{UAAU}), we have that
$$
F_{n+1,m} =  \frac{d}{d\lambda}\Psi_{n+1,m}  - A_{n+1,m}\Psi_{n+1,m}
= \frac{d}{d\lambda}\left(U_{n,m}\Psi_{n,m}\right)  - A_{n+1,m}U_{n,m}\Psi_{n,m}
$$
$$
=\left(A_{n+1,m}U_{n,m} - U_{n,m}A_{n,m}\right)\Psi_{n,m} + U_{n,m}\frac{d}{d\lambda}\Psi_{n,m} 
- A_{n+1,m}U_{n,m}\Psi_{n,m} = U_{n,m}F_{n,m},
$$
which means that
$$
F_{n,m}(\lambda) \equiv \Psi_{n,m}(\lambda)C_{m}(\lambda),
$$
where the matrix $C_{m}(\lambda)$ does not depend on $n$, but  might depend on $m$
and $\lambda$. Similar arguments based on the relation (\ref{VAAV}) yields the $m$ - independence 
of the matrix $C_{m}(\lambda)$,
$$
C_{m} (\lambda) \equiv C_{0}(\lambda) \equiv C(\lambda).
$$
It remains to notice that $F_{0,0} \equiv 0$ and hence $C(\lambda)\equiv 0$.
This completes the proof of the proposition.

We shall now proceed with the establishing of the monodromy properties of the
function $\Psi_{n,m}(\lambda)$. 
\begin{itemize}
\item{\it The neighborhood of the point $\lambda = 0$.} This is the easiest. Indeed,
from the definition (\ref{Psinmdef}), we immediately conclude that
\begin{equation}\label {00}
\Psi_{n,m}(\lambda) = \hat{\Psi}^{(0)}_{n,m}(\lambda)\lambda^{-\frac{{a}}{4}\sigma_3},
\end{equation}
where $\hat{\Psi}^{(0)}_{n,m}(\lambda)$ and   $[\hat{\Psi}^{(0)}_{n,m}(\lambda)]^{-1}$
are  holomorphic at $\lambda = 0$.
Moreover,
$$ 
\hat{\Psi}^{(0)}_{n,m}(0) = \begin{pmatrix}1&-f_{n,m}\cr
 0&1\end{pmatrix}.
$$
\item{\it The neighborhood of the point $\lambda = \infty$.} With the help
of a straightforward induction, one can easy check that in the 
neighborhood of infinity, the function $\Psi_{n,m}(\lambda)$ admits
the following representation. 
\begin{equation}\label{Psinminfty0}
\Psi_{n,m}(\lambda) = \hat{\Psi}^{(\infty)}_{n,m}(\lambda)\lambda^{-T_{\infty}},
\end{equation}
where
$$
T_{\infty} = \frac{{a}}{4}\sigma_3 -
\begin{cases}
\begin{pmatrix}[\frac{m+n}{2}]+1&0\cr
0&[\frac{m+n}{2}]\end{pmatrix} & \mbox{if}\,\, m+n\,\, \mbox{is odd}\cr\cr
\frac{m+n}{2}I& \mbox{if}\,\, m+n\,\, \mbox{is even}
\end{cases}.
$$
The functions  $\hat{\Psi}^{(\infty)}_{n,m}(\lambda)$ and   $[\hat{\Psi}^{(\infty)}_{n,m}(\lambda)]^{-1}$
are  holomorphic at $\lambda = \infty$.
Moreover,
\begin{equation}\label{Psihatinf}
 \hat{\Psi}^{(\infty)}_{n,m}(\infty) =
 \begin{cases}
  \begin{pmatrix}0&\bullet\cr
 \bullet&\bullet\end{pmatrix}& \mbox{if}\,\, m+n\,\, \mbox{is odd}\cr\cr
 \begin{pmatrix}\bullet&0\cr
 \bullet&\bullet\end{pmatrix}& \mbox{if}\,\, m+n\,\, \mbox{is even}
 \end{cases}
\end{equation}
%and (cf. (\ref{detPsi}))
%\begin{equation}\label{detPsihatinf}
%\det \hat{\Psi}^{(\infty)}_{n,m}(\infty) = (-1)^{m}.
%\end{equation} 
The symbol ``$\bullet$''  indicates that no specific conditions are imposed 
on the corresponding entry. In other words, description (\ref{Psihatinf}) of the 
matrix  $\hat{\Psi}^{(\infty)}_{n,m}(\infty)$  is equivalent to the statement that this
matrix satisfies the following property,
$$
\left[\hat{\Psi}^{(\infty)}_{n,m}(\infty)\right]_{11} = 0,\quad\mbox{if}\quad n+m\quad\mbox{is odd}
\quad \mbox{and}\quad 
\left[\hat{\Psi}^{(\infty)}_{n,m}(\infty)\right]_{12} = 0,\quad\mbox{if}\quad n+m\quad\mbox{is even}.
$$
\item{\it The neighborhood of the point $\lambda = -1$, $n \geq 1$.} 
The eigenvalues of the residue matrix $-B_{n,m}$ 
at the point $\lambda =-1$ are $n$ and $0$. Hence, by the general theory of 
differential equations with rational coefficients (see e.g. \cite{JMU}),
\begin{equation}\label {-10}
\Psi_{n,m}(\lambda) = \hat{\Psi}^{(-1)}_{n,m}(\lambda)\begin{pmatrix}1&0\cr
                               0&(\lambda+1)^{n}\end{pmatrix}E.
\end{equation}
where $\hat{\Psi}^{(-1)}_{n,m}(\lambda)$ and  $[\hat{\Psi}^{(-1)}_{n,m}(\lambda)]^{-1}$ 
are  holomorphic at $\lambda = -1$, and $E$ does not depend on $\lambda$.
(We note that the absence of the logarithmic terms at  $\lambda = -1$
follows  from the  very definition of the function $\Psi_{n,m}(\lambda)$.) Matrix 
$E$ in formula (\ref{-10}) is defined  up to the left multiplication by a lower triangular
matrix factor, and it can be brought either to the form, 
\begin{equation}\label{choice1}
E = \begin{pmatrix}1&c_{n,m}\cr
0&1\end{pmatrix}
\end{equation}
($E_{11} \neq 0$), or to the form,
\begin{equation}\label{choice2}
E = \begin{pmatrix}0&1\cr
1&0\end{pmatrix}
\end{equation}
($E_{11} = 0$). We argue that the structure of the matrix $E$ must be the same  for all $n$, $m$.
Indeed, let us suppose that
$$
E_{n,m} = \begin{pmatrix}1&c_{n,m}\cr
0&1\end{pmatrix}, \quad \mbox{while}\quad E_{n+1, m} = \begin{pmatrix}0&1\cr
1&0\end{pmatrix}.
$$
Then, from (\ref{-10}) it follows that 
\begin{equation}\label{check}
\Psi_{n+1,m}(\lambda)\Psi_{n,m}^{-1}(\lambda) = H(\lambda)
\begin{pmatrix}1& (\lambda+1)^{-n}\cr
 (\lambda+1)^{n+1}&c_{n,m}(\lambda + 1)\end{pmatrix}\tilde{H}(\lambda),
\end{equation}
where $H(\lambda)$ and $\tilde{H}(\lambda)$ are holomorphic 
at $\lambda = -1$ functions. On the other hand, the left hand side of
the last equation is nothing else but $U_{n,m}(\lambda)$, which is holomorphic.
Moreover, the matrices $H(-1)$ and $\tilde{H}(-1)$ are invertible.
Therefore, the cancellation of singularity at $\lambda = -1$ in the right
hand side of (\ref{check}) is  not possible, and we ran into a contradiction. The reader
can easily check that the similar contradiction arrises if we assume that
$$
E_{n,m} = \begin{pmatrix}0&1\cr
1&0\end{pmatrix}, \quad E_{n+1, m} = \begin{pmatrix}1&c_{n+1,m}\cr
0&1\end{pmatrix},
$$
as well as if we assume that
$$
E_{n,m} = \begin{pmatrix}1&c_{n,m}\cr
0&1\end{pmatrix}, \quad E_{n, m+1} = \begin{pmatrix}0&1\cr
1&0\end{pmatrix}.
$$
or
$$
E_{n,m} = \begin{pmatrix}0&1\cr
1&0\end{pmatrix}, \quad E_{n, m+1} = \begin{pmatrix}1&c_{n,m+1}\cr
0&1\end{pmatrix},
$$
Observe now that in the  case of option (\ref{choice2}), the matrix $E$ does not depend on $m,n$.  In fact,
the same is true  even  if it is the option (\ref{choice1}) that is realized for all $m,n$.
Indeed, using again  (\ref{-10}) we  see that 
$$
\Psi_{n+1,m}(\lambda)\Psi_{n,m}^{-1}(\lambda) = H(\lambda)
\begin{pmatrix}1& (c_{n+1,m} -c_{n,m})(\lambda+1)^{-n}\cr
0&\lambda + 1\end{pmatrix}\tilde{H}(\lambda),
$$ 
where, as before, $H(\lambda)$ and $\tilde{H}(\lambda)$ are holomorphic 
at $\lambda = -1$ functions. Once again, the left hand side of
the last equation is $U_{n,m}(\lambda)$, which is holomorphic, while 
the matrices $H(-1)$ and $\tilde{H}(-1)$ are invertible. 
The cancellation of singularity at $\lambda = -1$
is now possible, and it is possible only  if,
$$
c_{n+1,m} = c_{n,m}.
$$
Similarly,
$$
\Psi_{n,m+1}(\lambda)\Psi_{n,m}^{-1}(\lambda) = J(\lambda)
\begin{pmatrix}1& (c_{n,m+1} -c_{n,m})(\lambda+1)^{-n}\cr
0&1\end{pmatrix}\tilde{J}(\lambda) \equiv V_{n,m}(\lambda)
$$ 
where $J(\lambda)$ and $\tilde{J}(\lambda)$ are again holomorphic 
at $\lambda = -1$ functions. Holomorphicity of $V_{n,m}(\lambda)$
would then imply that
$$
c_{n,m+1} = c_{n,m}.
$$

Just established independence of the parameter $E$, in the both its possible 
forms on $n$ and $m$
allows us to evaluate it by analyzing  the  function $\Psi_{1,0}(\lambda)$.
We have,
$$
\Psi_{1,0}(\lambda) = U_{0,0}(\lambda)\Psi_{0,0}(\lambda) =
 \begin{pmatrix}1&-1\cr
\lambda&1\end{pmatrix}\lambda^{-\frac{{a}}{4}\sigma_3}
= \begin{pmatrix}\lambda^{-\frac{{a}}{4}}&- \lambda^{\frac{{a}}{4}}\cr
\lambda^{-\frac{{a}}{4} +1}&\lambda^{\frac{{a}}{4}}
\end{pmatrix}.
$$
Consider the product
$$
\begin{pmatrix}\lambda^{-\frac{{a}}{4}}&- \lambda^{\frac{{a}}{4}}\cr
\lambda^{-\frac{{a}}{4} +1}&\lambda^{\frac{{a}}{4}}
\end{pmatrix}
\begin{pmatrix}1&-c\cr
0&1\end{pmatrix}\begin{pmatrix}1&0\cr
0&(\lambda+1)^{-1}\end{pmatrix} \equiv \hat{\Psi}^{(-1)}_{1,0}(\lambda)
$$
It is straightforward that this product is holomorphic at
$\lambda =-1$ iff
$$
c=-e^{\frac{\pi i{a}}{2}}.
$$
Hence the matrix $E$ in representation (\ref{-10}) is
given by the equation,
\begin{equation}\label{C-10}
E = \begin{pmatrix}1&-e^{\frac{i\pi{a}}{2}}\cr
                            0&1\end{pmatrix}.
\end{equation}
\item {\it The neighborhood of the point $\lambda = 1$, $m \geq 1$.} 
Repeating exactly the same arguments as in the previous case,
we arrive at the following representation of the function
$\Psi_{n,m}(\lambda)$ in the neighborhood of the point $\lambda = 1$.
\begin{equation}\label {10}
\Psi_{n,m}(\lambda) = \hat{\Psi}^{(1)}_{n,m}(\lambda)\begin{pmatrix}1&0\cr
                               0&(\lambda-1)^{m}\end{pmatrix}E.
\end{equation}
where $\hat{\Psi}^{(1)}_{n,m}(\lambda)$ and  $[\hat{\Psi}^{(1)}_{n,m}(\lambda)]^{-1}$ 
are  holomorphic at $\lambda = 1$, and the constant matrix $E$ 
is given by the {\it same}  equation (\ref{C-10}). This completes the
evaluation of the monodromy data of the linear system (\ref{lambda})
corresponding to the discrete map $Z^{{a}}$.  
\end{itemize}
The branch of the function $\lambda^{-{a}/4}$ appearing in (\ref{00}), (\ref{Psinminfty0})  
is defined on the $\lambda$ - plane
cut along the ray $[0, -i\infty)$ and it is fixed by the condition $-\frac{\pi}{2} < \arg \lambda< \frac{3\pi}{2}$.
It also should be noticed that  the product,
\begin{equation}\label{prod}
\Psi_{n,m}(\lambda)\lambda^{\frac{{a}}{4}\sigma_3},
\end{equation}
is analytic and single valued on the whole finite $\lambda$ - plane. It is, in fact, a matrix
polynomial. 

\begin{prop}\label{unique} Representations  (\ref{00}), (\ref{Psinminfty0}),
(\ref{-10}) and (\ref{10}) (with the matrix $E$ defined in  (\ref{C-10})) together
with the single-validness of the product (\ref{prod}) determine 
the function $\Psi_{n,m}(\lambda)$
uniquely. 
\end{prop}
{\it Proof.} Suppose that $\tilde{\Psi}_{n,m}(\lambda)$ is another matrix valued function 
which admits representations (\ref{00}), (\ref{Psinminfty0}),
(\ref{-10}) and (\ref{10}) (with the matrix $E$ defined in  (\ref{C-10})) and such
that the product 
$$
\tilde{\Psi}_{n,m}(\lambda)\lambda^{\frac{{a}}{4}\sigma_3},
$$
is analytic and single valued. Put
$$
\Theta(\lambda):= \tilde{\Psi}_{n,m}(\lambda)[\Psi_{n,m}(\lambda)]^{-1}.
$$
We first notice that
$$
\tilde{\Psi}_{n,m}(\lambda)[\Psi_{n,m}(\lambda)]^{-1} = 
\left(\tilde{\Psi}_{n,m}(\lambda)\lambda^{\frac{{a}}{4}\sigma_3}\right)
[\Psi_{n,m}(\lambda)\lambda^{\frac{{a}}{4}\sigma_3}]^{-1}
$$
and hence the function $\Theta(\lambda)$ is single valued. Secondly,
because of (\ref{detPsi}), the inversion of the matrix $\Psi_{n,m}(\lambda)$ does not
produce new singularities. Therefore, one can conclude that {\it a priori},
the function $\Theta(\lambda)$ is analytic on $\C\setminus \{0, \infty, 1, -1\}$.
At the same time, at the points $0$, $\infty$ and $\pm 1$,
the both functions which form the product $\Theta(\lambda)$  have exactly the same right 
singular factors which cancel out in the product. Therefore, we conclude that the function $\Theta(\lambda)$
is, in fact, a constant function,
$$
\Theta(\lambda) \equiv \mbox{constant}
$$
Now, evaluating this constant matrix at $\lambda =0$ we see that
\begin{equation}\label{Theta1}
\Theta(\lambda) \equiv \begin{pmatrix}1&\bullet\cr
                               0& 1\end{pmatrix},
\end{equation}
while the evaluation at $\lambda = \infty$ yields
\begin{equation}\label{Theta2}
\Theta(\lambda) \equiv
 \begin{pmatrix}\bullet&0\cr
 \bullet&\bullet
 \end{pmatrix}
 \end{equation}
(regardless the parity of $m+n$). Comparing (\ref{Theta1}) and (\ref{Theta2}) we conclude that
$$
\Theta(\lambda) \equiv I.
$$
The proposition is proven.

It is important to emphasize that we {\it do not need} to prescribe {\it a priori}  
the 12  matrix entry of the matrix $\hat{\Psi}^{(0)}_{n,m}(0).$ In fact,
we shall use the equation,
\begin{equation}\label{fdef}
f_{n,m} = - [\hat{\Psi}^{(0)}_{n,m}(0)]_{12},
\end{equation}
as an {\it independent } definition of the map $Z^{{a}}$.

We conclude this section by noticing that the independence on $n$ and $m$
of the matrix  $E$ means again  that the discrete map $Z^{{a}}$ describes 
a special one parameter family of discrete isomonodromy deformations
of system (\ref{lambda}).  In fact, this one parameter family can be also identified
with a special solution of a discrete Painlev\'e equation, namely of the d-PII equation,
see \cite{AB}. Moreover, there is also a connection to the continuous Painlev\'e
equations. The map $Z^{{a}}$ can be also obtained via the Backulnd transformation
of a special solution of the sixth Painlev\'e equation. This '' Painlev\'e connections'', however,
does not help in our main problem, which is the evaluation of the large $n$, $m$ asymptotics
of the map $Z^{{a}}$. Rather, the results of our paper might be used in building up a comprehensive 
asymptotic theory of Painlev\'e functions. For the modern state of the art in this area we refer
the reader to the monograph \cite{FIKN} and to the more recent source \cite{CA}.

\subsection{The Riemann-Hilbert setting.}
>From now on we shall assume that $m+n$ is even. That is, we will first prove Theorem \ref{theorem1}
for this case. An extension of the statement of the theorem to the case of  arbitrary parity of $m+n$  will
be done in the last section of the paper, in Section \ref{oddcase}.

We start with summarizing  the previous section's considerations as the
following theorem.
\begin{theorem}\label{theorem2}
Let $\Psi_{n,m}(\lambda)$ be the matrix valued function defined by the discrete conformal map $Z^{{a}}_{n,m}$ according to the equations (\ref{Psinmdef}).
Then, the function $\Psi_{n,m}(\lambda)\equiv\Psi(\lambda)$ is the unique solution of the following analytical problem. 
%consisting in finding 
%of a matrix valued function $\Psi(\lambda)$ satisfying the following properties.
%Then, according to the
%previous subsection, the {\it isomonodromy representation }
%of the map $Z^{{a}}(n,m)$ is obtained as follows.
%Let the $2\times 2$ matrix valued function  $\Psi(\lambda) \equiv \Psi_{n,m}(\lambda)$ be the
%(unique) solution of the following analytic problem.
\begin{itemize}
\item In the vicinity of $\lambda = 1$, the function  $\Psi(\lambda)$ admits
the representation,
\begin{equation}\label {1}
\Psi(\lambda) = \hat{\Psi}^{(1)}(\lambda)\begin{pmatrix}1&0\cr
                               0&(\lambda-1)^{m}\end{pmatrix}
                                \begin{pmatrix}1&-e^{\frac{i\pi{a}}{2}}\cr
                            0&1\end{pmatrix},
\end{equation}
where $\hat{\Psi}^{(1)}(\lambda)$ is holomorphic and invertible at $\lambda = 1$.
\item In the vicinity of $\lambda = -1$, the function  $\Psi(\lambda)$ admits
the representation,
\begin{equation}\label {-1}
\Psi(\lambda) = \hat{\Psi}^{(-1)}(\lambda)\begin{pmatrix}1&0\cr
                               0&(\lambda+1)^{n}\end{pmatrix}
                                \begin{pmatrix}1&-e^{\frac{i\pi{a}}{2}}\cr
                            0&1\end{pmatrix},
\end{equation}
where $\hat{\Psi}^{(-1)}(\lambda)$ is holomorphic and invertible at $\lambda = -1$.
\item In the vicinity of $\lambda = \infty$, the function  $\Psi(\lambda)$ admits
the representation,
\begin{equation}\label {infty}
\Psi(\lambda) = \hat{\Psi}^{(\infty)}(\lambda)\lambda^{-\frac{{a}}{4}\sigma_3}\lambda^{\frac{m+n}{2}},
\end{equation}
where $\hat{\Psi}^{(\infty)}(\lambda)$ is holomorphic and invertible at $\lambda = \infty$.
Moreover,
$$ 
\hat{\Psi}^{(\infty)}(\infty) = \begin{pmatrix}\bullet&0\cr
 \bullet&\bullet\end{pmatrix}.
$$
\item In the vicinity of $\lambda = 0$, the function  $\Psi(\lambda)$ admits
the representation,
\begin{equation}\label {0}
\Psi(\lambda) = \hat{\Psi}^{(0)}(\lambda)\lambda^{-\frac{{a}}{4}\sigma_3},
\end{equation}
where $\hat{\Psi}^{(0)}(\lambda)$ is holomorphic and invertible at $\lambda = 0$.
Moreover,
$$ 
\hat{\Psi}^{(0)}(0) = \begin{pmatrix}1&\bullet\cr
 0&1\end{pmatrix}.
$$
The branch of the function $\lambda^{-{a}/4}$ in (\ref{infty}) and (\ref{0}) is defined on the $\lambda$ - plane
cut along the ray $[0, -i\infty)$ and it is fixed by the condition $-\frac{\pi}{2} < \arg \lambda< \frac{3\pi}{2}$.
\item   The product,
$$
\Psi(\lambda)\lambda^{\frac{{a}}{4}\sigma_3},
$$
is analytic and single valued on the whole finite $\lambda$ - plane (it is, in fact, a matrix
polynomial). 
\end{itemize}
The map $Z^{{a}}$ itself can be recovered from the known function $\Psi$ by the relation,
%Having a solution of the above monodromy  $\Psi$ - problem, the map $Z^{{a}}$
%is given by the relation (cf. (\ref{fdef})),
\begin{equation}\label{Zgamma}
Z^{{a}}_{n,m} = - [\hat{\Psi}^{(0)}_{n,m}(0)]_{12}.
\end{equation}
\end{theorem}

We shall call the problem (\ref{1}) -- (\ref{0}) - the   {\it monodromy problem}, and we will  be saying  that
formula (\ref{Zgamma}) gives the {\it monodromy representation} of the discrete
power function $Z^{{a}}$.
We shall now perform  a series of  equivalent reformulations of the monodromy problem 
which will eventually  transform it to a Riemann-Hilbert factorization problem posed on the ray 
$[0, -i\infty)$.

\vskip .2 in
\noindent
{\it Step 1.} Put
\begin{equation}\label{Phi}
\Phi(\lambda) = \Psi(\lambda)\lambda^{\frac{{a}}{4}\sigma_3}.
\end{equation}
This simple transformation makes the new object - the function $\Phi(\lambda)$, a single valued 
function on the whole $\lambda$ - plane.  In terms of the  function $\Phi(\lambda)$ the monodromy problem 
reads as follows.

\begin{figure}[h]
  \begin{center}
     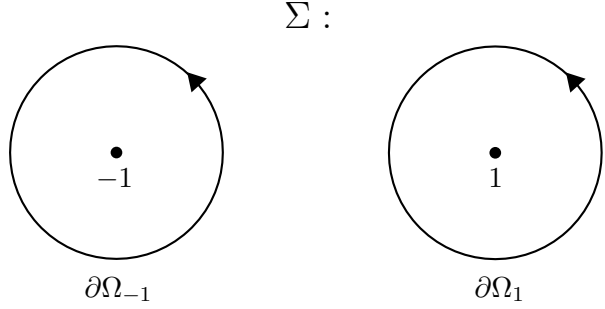
  \end{center}
  \caption{Contour for the $X$ - RH problem}\label{f.contour_X}
\end{figure}

\begin{itemize}
\item The function $\Phi(\lambda)$ is analytic on the finite $\lambda$ - plane.
\item In the vicinity of $\lambda = 1$, the function  $\Phi(\lambda)$ admits
the representation,
\begin{equation}\label {11}
\Phi(\lambda) = \hat{\Phi}^{(1)}(\lambda)\begin{pmatrix}1&0\cr
                               0&(\lambda-1)^{m}\end{pmatrix}
                                \begin{pmatrix}1&-\lambda^{-\frac{{a}}{2}}e^{\frac{i\pi{a}}{2}}\cr
                            0&1\end{pmatrix},
\end{equation}
where $\hat{\Phi}^{(1)}(\lambda)$ is holomorphic and invertible at $\lambda = 1$.
\item In the vicinity of $\lambda = -1$, the function  $\Phi(\lambda)$ admits
the representation,
\begin{equation}\label {-11}
\Phi(\lambda) = \hat{\Phi}^{(-1)}(\lambda)\begin{pmatrix}1&0\cr
                               0&(\lambda+1)^{n}\end{pmatrix}
                                \begin{pmatrix}1&-\lambda^{-\frac{{a}}{2}}e^{\frac{i\pi{a}}{2}}\cr
                            0&1\end{pmatrix},
\end{equation}
where $\hat{\Phi}^{(-1)}(\lambda)$ is holomorphic and invertible at $\lambda = -1$.
\item In the vicinity of $\lambda = \infty$, the function  $\Phi(\lambda)$ admits
the representation,
\begin{equation}\label {infty1}
\Phi(\lambda) = \hat{\Phi}^{(\infty)}(\lambda)\lambda^{\frac{m+n}{2}},
\end{equation}
where $\hat{\Phi}^{(\infty)}(\lambda)$ is holomorphic and invertible at $\lambda = \infty$.
Moreover,
$$ 
\hat{\Phi}^{(\infty)}(\infty) = \begin{pmatrix}\bullet&0\cr
 \bullet&\bullet\end{pmatrix}.
$$
\item The function $\Phi(\lambda)$ is normalized by the condition,
$$
\Phi(0) = 
 \begin{pmatrix}1&\bullet\cr
 0&1\end{pmatrix}.
 $$
\end{itemize}
In the vicinity of $\lambda = 0$ we have that (see (\ref{0}))
$$
\Phi(\lambda) = \hat{\Psi}^{(0)}(\lambda).
$$
Therefore, equation (\ref{Zgamma}) becomes the equation,
\begin{equation}\label{Zgamma2}
Z^{({a})} = -\Phi_{12}(0).
\end{equation}
(From now on and until  Section \ref{oddcase}, we will usually suppress the indication of the
$n,m$ - dependence.)
\vskip .2 in
\noindent
{\it Step 2.} Observe that  equations (\ref{11}) and (\ref{-11}) can
be rewritten in more uniform way, i.e. 
\begin{equation}\label {111}
\Phi(\lambda) = \tilde{\Phi}^{(1)}(\lambda)\begin{pmatrix}1&0\cr
                               0&(\lambda-1)^{m}(\lambda+1)^{n}\end{pmatrix}
                                \begin{pmatrix}1&-\lambda^{-\frac{{a}}{2}}e^{\frac{i\pi{a}}{2}}\cr
                            0&1\end{pmatrix},
\end{equation}
and 
\begin{equation}\label {-111}
\Phi(\lambda) = \tilde{\Phi}^{(-1)}(\lambda)\begin{pmatrix}1&0\cr
                               0&(\lambda-1)^{m}(\lambda+1)^{n}\end{pmatrix}
                                \begin{pmatrix}1&-\lambda^{-\frac{{a}}{2}}e^{\frac{i\pi{a}}{2}}\cr
                            0&1\end{pmatrix},
\end{equation}
with the functions  $ \tilde{\Phi}^{(1)}(\lambda)$ and  $\tilde{\Phi}^{(-1)}(\lambda)$
possessing the same properties as the functions 
$ \hat{\Phi}^{(1)}(\lambda)$ and  $\hat{\Phi}^{(-1)}(\lambda)$, respectively.
 Let now $\Omega_{1}$ and $\Omega_{-1}$ denote the discs of radius
$1/2$ and centered at $\lambda = 1$ and $\lambda = -1$, respectively. Define,
\begin{equation}\label{X}
X(\lambda ) = \begin{cases}
\tilde{\Phi}^{(1, -1)} (\lambda) & \lambda \in \Omega_{1,-1}\cr\cr
\Phi(\lambda)\begin{pmatrix}1&0\cr
                               0&(\lambda-1)^{-m}(\lambda+1)^{-n}\end{pmatrix}
                               & \lambda \in \C\setminus \Omega_{1}\cup\Omega_{-1}
                               \end{cases}
\end{equation}
The function $X(\lambda)$ satisfies a certain  factorization {\it Riemann-Hilbert problem}
posed on the contour,
$$
\Sigma = \partial \Omega_{1}\cup\partial \Omega_{-1},
$$
which is depicted in Figure \ref{f.contour_X}. The orientation of the circles which form the contour is counterclockwise. 
As usual, the orientation defines a $+$ and a $-$ 
side on each part of the contour, where the $+$ side is on the left when traversing the contour according to its orientation.
The Riemann-Hilbert problem  which the function  $X(\lambda)$  solves reads as follows.

\vskip .2in
\noindent
{\bf Riemann-Hilbert problem for $X(\lambda)$} 
\begin{itemize}
\item $X(\lambda)$ is analytic on $\C\setminus \Sigma$.
\item The boundary values, 
%$X_{\pm}(\lambda)$, of $X(\lambda)$ on $\Sigma$, exist 
%point-wise, i.e.,
\begin{equation}\label{bound}
X_{\pm}(\lambda) := \lim_{\lambda' \to \lambda, \,\, \lambda' \in \pm\mbox{side of}\,\Sigma}X(\lambda'),
\end{equation}
of $X(\lambda)$  on $\Sigma$ exist point-wise, the limits in (\ref{bound}) are uniform,  and  the
functions $X_{\pm}(\lambda)$ are continuous. Moreover, the functions $X_{\pm}(\lambda)$ satisfy the 
jump condition,
\begin{equation}\label{Xjump}
X_+(\lambda) = X_-(\lambda)G(\lambda), \quad \lambda \in \Sigma,
\end{equation}
where
$$
G(\lambda) =  \begin{pmatrix}1&e^{\frac{i\pi{a}}{2}}\lambda^{-\frac{{a}}{2}}
(\lambda - 1)^{-m}(\lambda+1)^{-n}\cr
                            0&1\end{pmatrix}.
$$
%The boundary values are defined point-wise,
%$$
%X_{\pm}(\lambda) = \lim_{\lambda' \to \lambda, \,\, \lambda' \in \pm\mbox{side of}\,\Sigma}X(\lambda'),
%$$
%and the limits are uniform.
\item The behavior of the function $X(\lambda)$ at the point $\lambda = \infty$
is described by the equation, 
\begin{equation}\label{Xinfty}
X(\lambda) = \begin{pmatrix}\bullet&0\cr
 \bullet&\bullet\end{pmatrix}\left(I + O\left(\frac{1}{\lambda}\right)\right)\lambda^{\frac{m+n}{2}\sigma_3},
 \quad \lambda \to \infty.
\end{equation}
\item The function $X(\lambda)$ is normalized by the condition,
\begin{equation}\label{X0}
X(0) = \begin{pmatrix}1&\bullet\cr
 0&1\end{pmatrix}\begin{pmatrix}1&0\cr
 0&(-1)^{m}\end{pmatrix}.
 \end{equation}
 \end{itemize}
We shall call the problem (\ref{Xjump}) -- (\ref{X0}) - the {\it X - RH problem}. In terms of the $X$-RH problem,
relation (\ref{Zgamma2}) becomes the equation,
\begin{equation}\label{Zgamma3}
Z^{({a})} = (-1)^{m+1}X_{12}(0).
\end{equation}

%The contour $\Sigma$ for the  X - RH problem is depicted in
%Figure 1.  The orientation of the circles which form the contour is counterclockwise. As usual, the orientation defines a $+$ and a $-$ 
%side on each part of the contour, where the $+$ side is on the left when traversing the contour according to its orientation. For a 
%function $F \in \C\setminus \Sigma$, the symbols  $F_{\pm}$ denote the limiting values of the function $F$ on $\Sigma$ from the
%$\pm$ -side. We note, that for the function $X(\lambda)$ its limiting values exist in the point-wise sense, i.e.
%$$
%X_{\pm}(\lambda) = \lim_{\lambda' \to \lambda, \,\, \lambda' \in \pm\mbox{side of}\,\Sigma}X(\lambda'),
%$$
%and the limits are uniform.
 
\vskip .2 in
\noindent
 {\it Step 3.} Put
 \begin{equation}\label{Y}
 Y(\lambda) = \begin{cases}
 X(\lambda)G(\lambda) & \lambda \in \C \setminus \Omega_{1}\cup \Omega_{-1}\cup [0, -i\infty) \cr\cr
 X(\lambda) & \lambda \in \Omega_{1}\cup \Omega_{-1}
 \end{cases}
 \end{equation}
This transformation moves the jumps from the circles $\Sigma$ to the ray $\Sigma_0 = [0, -i\infty)$.
Assuming that the ray $\Sigma_0$ is oriented towards infinity, we arrive at the {\it $Y$ - RH problem}.

\vskip .2 in
\noindent
{\bf Riemann-Hilbert problem for $Y(\lambda)$}
\begin{itemize}
\item $Y(\lambda)$ is analytic on $\C\setminus \Sigma_0$.
\item The boundary values of $Y(\lambda)$ on $\Sigma_0\setminus\{0\}$ satisfy the 
jump condition,
\begin{equation}\label{Yjump}
Y_+(\lambda) = Y_-(\lambda)
\begin{pmatrix}1&\omega(\lambda)e^{-\varphi(\lambda)}\cr\cr
                            0&1\end{pmatrix}, \quad \lambda \in \Sigma_0\setminus\{0\},
\end{equation}
where
$$
\omega(\lambda) = 2i\sin{\frac{{a}\pi}{2}}\lambda_{+}^{-\frac{{a}}{2}}
(-1)^{m+n} 
= 2i\sin{\frac{{a}\pi}{2}}\lambda_{+}^{-\frac{{a}}{2}},
$$
and
$$
\varphi(\lambda) = m\log(\lambda -1) +n\log(\lambda+1).
$$
\item The behavior of the function $Y(\lambda)$ at the point $\lambda = \infty$
is described by the equation, 
\begin{equation}\label{Yinfty}
Y(\lambda) = \hat{Y}^{(\infty)}(\lambda)
 \begin{pmatrix}1&e^{\frac{i\pi{a}}{2}}\lambda^{-\frac{{a}}{2}}
\cr\cr
 0&1\end{pmatrix}(\lambda - 1)^{\frac{m}{2}\sigma_3}(\lambda+1)^{\frac{n}{2}\sigma_3}
\end{equation} 
where $\hat{Y}^{(\infty)}(\lambda)$ is holomorphic at $\lambda = \infty$, and
$$
\hat{Y}^{(\infty)}(\infty) = 
 \begin{pmatrix}\bullet&0\cr
 \bullet&\bullet\end{pmatrix}.
 $$
  \item The behavior of the function $Y(\lambda)$ at the point $\lambda = 0$
is described by the equation, 
\begin{equation}\label{Yzero}
Y(\lambda) = \hat{Y}^{(0)}(\lambda)
 \begin{pmatrix}1&e^{\frac{i\pi{a}}{2}}\lambda^{-\frac{{a}}{2}}
 (\lambda - 1)^{-m}(\lambda+1)^{-n}
\cr\cr
 0&1\end{pmatrix}
\end{equation} 
where $\hat{Y}^{(0)}(\lambda)$ is holomorphic at $\lambda = 0$, and
$$
\hat{Y}^{(0)}(0) = 
 \begin{pmatrix}1&\bullet\cr
 0&1\end{pmatrix}\begin{pmatrix}1&0\cr
 0&(-1)^{m}\end{pmatrix}.
 $$
 \end{itemize}
One can easily see that the function $\hat{Y}^{(0)}(\lambda)$ in representation (\ref{Yzero}) is
just the function $X(\lambda)$. Hence, from (\ref{Zgamma3}), we conclude that in terms of the $Y$ - RH problem, 
the discrete $Z^{{a}}$ is given by the equation,
\begin{equation}\label{ZY}
Z^{{a}} = (-1)^{m+1}\hat{Y}^{(0)}_{12}(0).
\end{equation} 

\begin{figure}[h]
  \begin{center}
     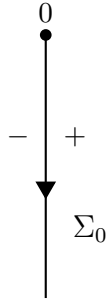
  \end{center}
  \caption{Contour for the $Y$ - RH problem}\label{f.contour_Y}
\end{figure}

The Y-RH problem is depicted in Figure \ref{f.contour_Y}. This problem is the final step in the series  of transformations of the original monodromy problem
(\ref{1}) -- (\ref{0}).
The Y-RH problem provides the {\it Riemann-Hilbert representation } for the conformal map $Z^{{a}}$
which we formulate as the following theorem.

\begin{theorem}\label{theorem3} Let $Y(\lambda)$ be the matrix valued function defined by the  discrete 
conformal map $Z^{{a}}$ according to the equations
(\ref{Y}), (\ref{X}), (\ref{Phi}), and (\ref{Psinmdef}). Then, the function $Y(\lambda)$ is the unique solution 
of the Riemann-Hilbert factorization problem (\ref{Yjump})  -- (\ref{Yzero}). The map $Z^{{a}}$ itself can be
recovered from the known function $Y$ by relation (\ref{ZY}).
\end{theorem}
 
\begin{remark}\label{rem1} In the setting of the $Y$ - RH problem, equation  (\ref{Yinfty}) can be replaced by the relation,
 \begin{equation}\label{Yinfty2}
 Y(\lambda) =  \begin{pmatrix}\bullet&0\cr
 \bullet&\bullet\end{pmatrix}\Bigl(I + o(1)\Bigr)\lambda^{\frac{n+m}{2}\sigma_3},\quad \lambda \to \infty,
\end{equation} 
while equation (\ref{Yzero}) --  by the relation:
\begin{equation}\label{Yzero2}
Y(\lambda) = \begin{pmatrix}1&\bullet\cr
 0&1\end{pmatrix}\begin{pmatrix}1&0\cr
 0&(-1)^{m}\end{pmatrix}\Bigl(I + o(1)\Bigr)
 \begin{pmatrix}1&(-1)^{m}e^{\frac{i\pi{a}}{2}}\lambda^{-\frac{{a}}{2}}
\cr\cr
 0&1\end{pmatrix}, \quad \lambda \to 0.
\end{equation} 
More precisely, this means that the more detailed formula (\ref{Yinfty}) follows from 
condition (\ref{Yinfty2}) and jump relation (\ref{Yjump}), and similar  is true for formula (\ref{Yzero}).
\end{remark}
%In terms of the $Y$ - RH problem, the discrete $Z^{{a}}$ is given by the equation,
%\begin{equation}\label{ZY}
%Z^{{a}} = (-1)^{m+1}\hat{Y}^{(0)}_{12}(0).
%\end{equation} 

The asymptotic solution of the $Y$ - RH problem will be our goal now. However, before we proceed 
with the relevant analysis, we want to make a brief algebraic detour and to
discuss a little bit  the  connection of the $Y$ - RH problem (\ref{Yjump}) - (\ref{Yzero}) to
a certain system of orhtogonal polynomials. 

\subsection{Connection to orthogonal polynomials.}\label{opol}

%\vskip .2 in
%\noindent
%{\bf Connection to orthogonal polynomials}
%\vskip .2in
The behavior of the function $Y(\lambda)$ at infinity, which is indicated in Remark \ref{rem1},
together with the upper triangularity of the jump matrix $G(\lambda)$,
shows that the $Y$ - RH problem belongs to the type of the Riemann-Hilbert problems which
appear  in the theory of orthogonal polynomials and random matrices \cite{FIK} (see also
monograph \cite{Deift} and survey \cite{I}). Indeed, the solution $Y(\lambda)$ of the $Y$ - RH problem
admits the following orthogonal polynomial representation.

\begin{equation}\label{Yort1}
Y(\lambda) = M\begin{pmatrix}P_k(\lambda)&\frac{1}{2\pi i}\int_{\Sigma_0}\frac{P_k(\mu)\omega(\mu)e^{-\varphi(\mu)}}{\mu-\lambda}
d\mu\cr\cr
-\frac{2\pi i}{h_{k-1}}P_{k-1}(\lambda)&-\frac{1}{h_{k-1}}\int_{\Sigma_0}\frac{P_{k-1}(\mu)\omega(\mu)e^{-\varphi(\mu)}}{\mu-\lambda}
d\mu\end{pmatrix},
\end{equation}
where 
$$
k =\frac{n+m}{2},
$$
$P_k(\lambda)$ and $P_{k-1}(\lambda)$ are the last two members of the collection $\{P_{l}(\lambda)\}_{l =0}^{k}$
of the  monic polynomials determined by the orthogonality conditions,
\begin{equation}\label{Yort2}
\int_{\Sigma_0}P_{l}(\lambda)\lambda^j\omega(\lambda)e^{-\varphi(\lambda)}d\lambda =0,
\quad j =0, ..., l-1, \quad l=0, ..., k,
\end{equation}
and $h_{k-1}$ is the squre of the  norm of the polynomial  $P_{k-1}(\lambda)$, i.e.
\begin{equation}\label{hl}
h_{k-1} = \int_{\Sigma_0}P_{k-1}(\lambda)\lambda^{k-1}\omega(\lambda)e^{-\varphi(\lambda)}d\lambda
\equiv \int_{\Sigma_0}P^{2}_{k-1}(\lambda)\omega(\lambda)e^{-\varphi(\lambda)}d\lambda.
\end{equation}
The  pre-factor $M$in (\ref{Yort1}) is the constant matrix uniquely determined by the structure
of the  matrices $\hat{Y}^{(\infty)}(\infty)$ and $\hat{Y}^{(0)}(0)$, see (\ref{Yinfty}) and (\ref{Yzero}).
The map $Z^{{a}}$ can be  expressed in terms of the polynomials $P_k(\lambda)$ as well.
The corresponding formulae are,
\begin{equation}\label{ZP}
Z^{{a}} = -\frac{1}{P_k(0)}\sum\mbox{res}_{\lambda =\pm 1}
 \left(P_k(\lambda)e^{\frac{i\pi{a}}{2}}\lambda^{-1 -\frac{{a}}{2}}(\lambda -1)^{-m}(\lambda +1)^{-n}\right).
\end{equation}
$$
= -\frac{1}{P_k(0)}\frac{e^{\frac{i\pi{a}}{2}}}{(m-1)!}\frac{d^{m-1}}{d\lambda^{m-1}}\left(
P_k(\lambda)\lambda^{-1 -\frac{{a}}{2}}(\lambda +1)^{-n}\right)\Bigl|_{\lambda=1}\Bigr.
$$
$$
-\frac{1}{P_k(0)}\frac{e^{\frac{i\pi{a}}{2}}}{(n-1)!}\frac{d^{n-1}}{d\lambda^{n-1}}\left(
P_k(\lambda)\lambda^{-1 -\frac{{a}}{2}}(\lambda -1)^{-m}\right)\Bigl|_{\lambda=-1}\Bigr.
$$
It is worth noticing that the orthogonality condition (\ref{Yort2}) can be also  re-written in
a simple residue form
$$
\sum \mbox{res}_{\lambda = \pm 1}\left(P_{l}(\lambda)e^{\frac{i\pi {a}}{2}}\lambda^{j-\frac{{a}}{2}}e^{-\varphi(\lambda)}\right)
\equiv \sum\mbox{res}_{\lambda = \pm 1}\left(P_{l}(\lambda)e^{\frac{i\pi {a}}{2}}\lambda^{j-\frac{{a}}{2}}(\lambda-1)^{-m}
(\lambda+1)^{-n}\right)
$$
$$
= \frac{e^{\frac{i\pi{a}}{2}}}{(m-1)!}\frac{d^{m-1}}{d\lambda^{m-1}}\left(
P_l(\lambda)\lambda^{j -\frac{{a}}{2}}(\lambda +1)^{-n}\right)\Bigl|_{\lambda=1}\Bigr.
$$
\begin{equation}\label{Yort3}
+
\frac{e^{\frac{i\pi{a}}{2}}}{(n-1)!}\frac{d^{n-1}}{d\lambda^{n-1}}\left(
P_l(\lambda)\lambda^{j -\frac{{a}}{2}}(\lambda -1)^{-m}\right)\Bigl|_{\lambda=-1}\Bigr.
 =0\quad j =0, ..., l-1, 
\end{equation}
which, in fact, allows one to extend the original  finite orthogonal polynomials system 
$\{P_{l}(\lambda)\}_{l =0}^{k}$ to the infinite system $\{P_{l}(\lambda)\}_{l =0}^{\infty}$.

The orthogonal polynomial connection  described in this 
subsection will make a little appearance in the rest of the paper. Mostly, we will use  it as  a motivation for certain steps in our asymptotic
analysis. Therefore,
we skip the formal discussion of the existence of the orthogonal polynomials $P_l(\lambda)$
which can be also considered as a direct consequence of  a prior existence 
of the function $Y(\lambda)$. We also skip the derivation of the  formula (\ref{Yort1}) itself.
It is standard (cf. \cite{FIK}, \cite{Deift}).  One only have to be a little bit more
careful, comparing with the usual cases, when deriving the asymptotic condition (\ref{Yinfty2}) from formula
(\ref{Yort1}) . Usually, in the Riemann-Hilbert approach to orthogonal polynomials the
weghts participajting in the orthogonality conditions are appeared to decay very fast at infinity, i.e., faster 
then any power. This is not the case with our weight, which itself has a power-like decay at infinity.
This, in particular,  means that the error $o(1)$ in (\ref{Yinfty2}) can not be replaced by $O(\lambda^{-1})$, as
it possible to do in the usual situation. Indeed,  in our case, $o(1) = 0(\lambda^{-{a}/2})$.

%\begin{remark} 
As it is always the case with the orthogonal polynomial Riemann-Hilbert 
problems (see e.g. \cite{Deift}),  one can extract from fromula
(\ref{Yort1}) and  orthogonality conditions (\ref{Yort2}) or (\ref{Yort3})   a Hankel type  
determinant representations for both, the solution of 
the Riemann-Hilbert problem  (\ref{Yjump}) - (\ref{Yzero}) and for our main object - the map $Z^{{a}}_{n,m}$.
Indeed, let 
\begin{equation}\label{moments}
H_{s} := \int_{\Sigma_0}\lambda^s\omega(\lambda)e^{-\varphi(\lambda)}d\lambda,\quad
s= 0, 1, ..., 2k-1 \equiv n+m -1,
\end{equation}
be the moments of the weight $\omega(\lambda)e^{-\varphi(\lambda)}$ and let
$$
{\cal{H}}_{l} = \{H_{k+j}\}_{j,k =0, ..., l-1}
\equiv
\begin{pmatrix} H_0& H_1&...&H_{l-1}\cr\cr
H_1& H_2&...&H_{l}\cr\cr
.....& ...&...&...\cr\cr
H_{l-1}& H_l&...&H_{2l-2}\end{pmatrix},\quad l \leq k
$$
be the corresponding $l\times l$ Hankel matrix. Define also the augmented Hankel
matrix,
$$
{\cal{H}}_{l+1}(\lambda) = \det\begin{pmatrix} H_0& H_1&...&H_{l-1}&H_l\cr\cr
H_1& H_2&...&H_{l}&H_{l+1}\cr\cr
.....& ...&...&...&...\cr\cr
H_{l-1}& H_l&...&H_{2l-2}&H_{2l-1}\cr\cr
1& \lambda &...&\lambda^{l-1}&\lambda^l\end{pmatrix},\quad l \leq k
$$
i.e. the Hankel matrix ${\cal{H}}_{l+1}$ with the last row replaced by the row of the successive  powers
of $\lambda$. Orthogonality condition (\ref{Yort2}) is a linear system for the coefficients
of polynomial $P_l(\lambda)$. Applying to this system  Cramer's rule  and after some simple
manipulations (see again, e.g. \cite{Deift}), we will arrive at the equations,
\begin{equation}\label{Pldet}
P_l(\lambda) = \frac{\det {\cal{H}}_{l+1}(\lambda)}{\det{\cal{H}}_{l}}, 
\quad h_l = \frac{\det {\cal{H}}_{l+1}}{\det{\cal{H}}_{l}},
\end{equation}
which, in conjunction with the formulae (\ref{Yort1}) and (\ref{ZP}), provide the determinant
representations for the solution of the $Y$ - RH problem and for the discrete power function $Z^{{a}}$.

It is worth noticing that the integral in (\ref{moments}) can be evaluated by residues, so that
the moments $H_s$ can be expressed in the following form,
$$
H_s
= 2\pi i\frac{e^{\frac{i\pi{a}}{2}}}{(m-1)!}\frac{d^{m-1}}{d\lambda^{m-1}}\left(
\lambda^{s-\frac{{a}}{2}}(\lambda +1)^{-n}\right)\Bigl|_{\lambda=1}\Bigr.
$$
\begin{equation}\label{moments2}
+
2\pi i\frac{e^{\frac{i\pi{a}}{2}}}{(n-1)!}\frac{d^{n-1}}{d\lambda^{n-1}}\left(
\lambda^{s -\frac{{a}}{2}}(\lambda -1)^{-m}\right)\Bigl|_{\lambda=-1}\Bigr..
\end{equation}
%$$
%\equiv 2\pi i\frac{e^{\frac{i\pi{a}}{2}}}{(n-1)!}\sum_{j=0}^{m-1}
%(-1)^{m-1-j}2^{-n-m+j+1}\frac{(n+m-j-2)!}{j!(m-1-j)!}\frac{\Gamma\left(s -\frac{{a}}{2} -j +2\right)}
%{\Gamma\left(s -\frac{{a}}{2} \right)}
%$$
%$$
%+2\pi i\frac{e^{\frac{i\pi s}{2}}}{(m-1)!}\sum_{j=0}^{n-1}
%(-1)^{m-j}2^{-n-m+j+1}\frac{(n+m-j-2)!}{j!(n-1-j)!}\frac{\Gamma\left(s -\frac{{a}}{2} -j +2\right)}
%{\Gamma\left(s -\frac{{a}}{2} \right)}
%$$
Alternatively, the moments $H_s$ can be expreseed in terms of the hypergeometric functions,
\begin{equation}\label{moments3} 
H_s = (-1)^s2\pi i\frac{\Gamma\left(m+n -1 +\frac{{a}}{2} -s\right)}
{\Gamma\left(\frac{{a}}{2} - s \right)(n+m)!}
F\left(m, 1 -\frac{{a}}{2} +s;m+n; 2\right).
\end{equation}

The determinant formulae for $Z^{{a}}_{n,m}$ similar to the ones presented above 
have already been obtained (without any use of the Riemann-Hilbert
analysis) in \cite{ando}.
However, the size of the determinants $\det {\cal H}_{k}$, $\det {\cal H}_{k\pm 1}$,
$\det {\cal H}_{k}(\lambda)$, and $\det {\cal H}_{k+1}(\lambda)$ which appear  in the representation 
of $Y(\lambda)$ and $Z^{{a}}$, grows unboundedly  as $n, m \to \infty$ which makes these determinant
formulae useless for the asymptotic analysis. We want to stress that this is a general feature of the orthogonal polynomial
theory. That is,  the Riemann-Hilbert problem is used to evaluate the asymptotics of the determinants  appearing in the
representations of orthogonal polynomials; not the other way round. This is also the reason why we decided to point out
at the relation of the Riemann-Hilbert problem  (\ref{Yjump}) - (\ref{Yzero}) to the system of orthogonal
polinomials (\ref{Yort2}) - (\ref{Yort3}). The latter might be of interest of their own, and the results
of our paper might be used for the describtion of the large $k$ behavior of the polynomials $P_k(\lambda)$,
as well as of the Hankel determinants $\det{\cal H}_{k}$ whose generating moments $H_s$ are given by the
formulae (\ref{moments2}) or (\ref{moments3}). 
%\end{remark}

\begin{remark}\label{gamma1remark}
In the special case ${a} =1$ everything of course trivializes. The unique solution of  (\ref{def1})-(\ref{incond})
is, as expected, $f_{n.m} = n+im \equiv Z^{{a}}|_{{a}=1}$ and the $Y$-RH problem admits an explicit (i.e.,
no growing with $n$ and $m$ nontrivial matrix products) solution. We discuss this issue in detail in Appendix C.
\end{remark}

\section{Asymptotic analysis.}
In the asymptotic analysis  of the $Y$ - RH problem we will follow the Deift-Zhou 
nonlinear steepest descent method for oscillatory Riemann-Hilbert problems \cite{DZ}.
More precisely, we shall use  the adaptation of the method to the Riemann-Hilbert problems arising  in the theory
of orthogonal polynomials and random matrices which was developed in \cite{DKMVZ}
(for a pedagogical exposition of the method see again \cite{Deift} and \cite{I}).

In our presentation we will use the specific terminology accepted in the nonlinear steepest descent 
method, such as ``$g$ - function'', local and global ``parametricies'', etc. (see e.g. \cite{Deift}).
% We won't always explain the origin
%of the terms and the motivation of the particular intermediate steps we use. Otherwise, the paper would become
%too big. We refer the interested reader for all the historical and methodological details to \cite{Deift}. 
%From the formal  mathematical point of view, our presentation will be, of course, self-contained.

Following the methodology of the  nonlinear steepest descent method, we will perform a series of additional
transformations of the $Y$ - RH problem. The aim is to arrive  at the RH problem whose
jump matrix is approaching the identity as $n, m \to \infty$. In the process of these transformations, 
we will solve in closed form certain local  Riemann-Hilbert problems and assemble 
these local solutions into a piece-wise analytic matrix valued function which will  
approximate solution of the whole $Y$ - RH problem. This, in turn, will produce our main results -  the asymptotic formula (\ref{ABconj}).

%In the process of these transformations
%we will construct an explicit  {\it parametrix}  which will approximate the solution of the $Y$ - RH problem
%and hence produce our main results -  the asymptotic formula (\ref{ABconj}). The parametrix, or the approximate
%solution, will itself be shown to consist of three parts: the two {\it local} parametrices which solve the local RH problems
%in the vicinities of the points $\lambda =0$ and $\lambda = \infty$ and the {\it global} parametrix which solves
%the RH-problem which appear  in between these two points.  

%A somewhat novel feature of the problem we are solving in this paper is that the global parametrix will
%not, as it is usually the case in the analysis of the orthogonal polynomial RH problems, coincide with the 
%parametrix at $\lambda = \infty$. This is because
%of the fact that the RH problem we are analyzing in this paper  came from the Fuchsian linear system.
%Equivalently, one can say that the origin of the new technical features we are dealing with here is the weak, indeed logarithmic,
%growth  of the potential defining the weight of the orthogonal polynomials associated with the $Y$ - RH problem.

%General scheme:
%$$
%Y(\lambda) -(g-\mbox{function})\to T(\lambda)-(\mbox{analysis of the $\Re h$ and lenses opening}) \to S(\lambda)
%$$
%$$
%-(\mbox{parametrix construction}) \to R(\lambda),\quad G_{R} \sim I
%$$
\subsection{First transformation $Y \to T$ }
The first step in the method of \cite{DKMVZ} is the  introduction of  the so-called  $g$-function.
Let us briefly describe this notion. For  more detailed exposition we refer the reader
to  monograph \cite{Deift}).

Orthogonal polynomial representation (\ref{Yort1}) of the function $Y(\lambda)$ implies
that
$$
Y_{11}(\lambda) = M_{11}P_k(\lambda), \quad k =\frac{n+m}{2}.
$$
On the other hand, taking a hint from  the general theory of  orthogonal
polynomials on the line with positive weights (see e.g. \cite{Deift}), one can
suggest that, as $n^2+m^2 \to \infty$, 
\begin{equation}\label{haine1}
P_k(\lambda) \sim  e^{g(\lambda)},
\end{equation}
where
\begin{equation}\label{g01}
g(\lambda) = \int\log(\lambda - \mu)d\nu_{0}(\lambda),
\end{equation}
and 
$d\nu_{0}(\lambda)$ is the {\it equilibrium measure } corresponding to the potential $\varphi(\lambda)$.
This, means that $d\nu_{0}(\lambda)$ is an extremal point of the ''energy''  functional,
$$
{\cal{E}} = \int_{\Sigma_0}\int_{\Sigma_0}\log|\lambda -\mu|d\nu(\lambda)d\nu(\mu)
-\int_{\Sigma_0}\varphi(\lambda)d\nu(\lambda),
$$
considered on the space of Borel measures on $\Sigma_0$ satisfying the restriction,
\begin{equation}\label{lagr}
\int_{\Sigma_0}d\nu(\lambda) = k
\end{equation}

It is not very difficult, at least on the heuristic level, to see that the Euler-Lagrange equations for the functional  ${\cal{E}}$ have the form,
\begin{equation}\label{euler1}
g_+(\lambda) +g_-(\lambda)  - \varphi(\lambda) = \mbox{constant},\quad \lambda \in J,
\end{equation}
where $J$ means the support of the measure $d\nu_{0}(\lambda)$ and subscrips $\pm$ indicates
the relevant boundary values of the function $g(\lambda)$.  Also, condition (\ref{lagr}) yields
the asymptotic condition,
\begin{equation}\label{lagr2}
g(\lambda) \sim k\log \lambda, \quad \lambda  \to \infty.
\end{equation}
Remember that $\varphi(\lambda) = m\log(\lambda -1) +n\log(\lambda +1)$. This means, in particular,
a rather slow grows of the potential $\varphi(\lambda)$ at the infinity and hence a natural assumption
that the support $J$ of the equilibrium measure $d\nu_{0}(\lambda)$ should in fact coincide
with the whole semi-line $\Sigma_0$. Therefore, one can look at the problem (\ref{euler1}) - (\ref{lagr2})
as at  a scalar Riemann-Hilbert problem posed on the semi-line $\Sigma_0 =[0, -i\infty)$. The problem
can be solved by standard techniques which yields the following formula
for the $g$ - function.
\begin{equation}\label{g}
g(\lambda) = m\log(1 + \sqrt{\lambda}) + n\log(i + \sqrt{\lambda}).
\end{equation}
Here,  $\sqrt{\lambda}$ is defined on the plane cut along $\Sigma_0 = [0, -i\infty)$ and
the branch is fixed by the condition $-\pi/2 < \arg \lambda < 3\pi/2$.
For the  logarithmic function, $\log w$, its  principal branch, i.e. $-\pi < \arg w < \pi $
is taken.

We shall not attempt to transform the above heuristic considerations into a rigorous
proof of the asymptotic relation (\ref{haine1}). Instead, in accord with the method of \cite{DKMVZ}
we shall use them as a motivation for the  first  transformation of the $Y$ - RH problem:
%This  transformation is described by the formulae,
\begin{equation}\label{Tdef}
Y(\lambda) \Longrightarrow T(\lambda) : = Y(\lambda)e^{-g(\lambda)\sigma_3},
\end{equation}
with the function $g(\lambda)$ given by formula (\ref{g}). It is also significant, that  exactly the 
same function  $g(\lambda)$ appears in explicit solution of the $Y$ - Riemann-Hilbert
problem in the case ${a} =1$, see (\ref{g0}).

Let us see  how does the $Y$- RH problem
change under the transformation (\ref{Tdef}). 
As it will become clear soon,  the usefulness of this transformation is based
on the following properties of the  function (\ref{g}), first three of which have already appeared as the Euler-Lagrange
equations (\ref{euler1}), (\ref{lagr2}).
\begin{itemize}
\item $g(\lambda)$ is analytic in ${\Bbb C}\setminus [0, -i\infty]$,
\item as $\lambda \in [0, -i\infty]$,
\begin{equation}\label{gpm}
g_+(\lambda) + g_{-}(\lambda) = m\ln (1-\lambda) + n\ln (1+\lambda) +i\pi n
\equiv \varphi(\lambda) +i\pi n, \quad \mbox{mod}\, (2\pi),
\end{equation}
\item as $\lambda \to \infty$,
\begin{equation}\label{ginfty}
g(\lambda ) = \frac{m+n}{2}\ln \lambda + o(1),
\end{equation}
\item as $\lambda \to 0$,
\begin{equation}\label{gzero}
g(\lambda ) = \frac{i\pi}{2}n + (m-in)\sqrt{\lambda} + 0(\lambda).
\end{equation}
\end{itemize}
In view of the asymptotic formula (\ref{ginfty}), transformation (\ref{Tdef}) regularizes the behavior
at infinity in the setting of the Riemann-Hilbert problem. Indeed, for the new function $T(\lambda)$
we have that at $\lambda = \infty$,
$$
T(\lambda) = \hat{T}^{(\infty)}(\lambda)
 \begin{pmatrix}1&e^{\frac{i\pi{a}}{2}}\lambda^{-\frac{{a}}{2}}
\cr\cr
 0&1\end{pmatrix}(\lambda - 1)^{\frac{m}{2}\sigma_3}(\lambda+1)^{\frac{n}{2}\sigma_3}
 e^{-g(\lambda)\sigma_3}
 $$
\begin{equation}\label{Tinfty}
=  \hat{T}^{(\infty)}(\lambda)
 \begin{pmatrix}1&e^{\frac{i\pi{a}}{2}}\lambda^{-\frac{{a}}{2}}
\cr\cr
 0&1\end{pmatrix}\left(I + O\left(\frac{1}{\sqrt{\lambda}}\right)\right)_{diag}
% \end{equation}
%\begin{equation}\label{Tinfty2}
 =  \begin{pmatrix}\bullet&0\cr
 \bullet&\bullet\end{pmatrix}\Bigl(I + o(1)\Bigr),
\end{equation} 
where $ \hat{T}^{(\infty)}(\lambda) = \hat{Y}^{(\infty)}(\lambda)$ is holomorphic at $\lambda = \infty$.
The behavior at $\lambda = 0$ does not change much. Indeed, the asymptotic equations (\ref{Yzero})
and (\ref{Yzero2}) are  transformed into the equations,
\begin{equation}\label{Tzero}
T(\lambda)=  \hat{T}^{(0)}(\lambda)
 \begin{pmatrix}1&e^{\frac{i\pi{a}}{2}}\lambda^{-\frac{{a}}{2}}
 (\lambda - 1)^{-m}(\lambda+1)^{-n}
\cr\cr
 0&1\end{pmatrix}e^{-g(\lambda)\sigma_3}
\end{equation} 
and 
\begin{equation}\label{Tzero2}
T(\lambda) = \begin{pmatrix}1&\bullet\cr
 0&1\end{pmatrix}\begin{pmatrix}1&0\cr
 0&(-1)^{m}\end{pmatrix}e^{-\frac{i\pi}{2}n\sigma_3}\Bigl(I + o(1)\Bigr)
 \begin{pmatrix}1& e^{\frac{i\pi{a}}{2}}\lambda^{-\frac{{a}}{2}}
\cr\cr
 0&1\end{pmatrix}\Bigl(I + o(1)\Bigr), \quad \lambda \to 0,
\end{equation} 
respectively. Here,  $\hat{T}^{(0)}(\lambda) =  \hat{Y}^{(0)}(\lambda)$ is holomorphic at $\lambda = 0$,
and equation (\ref{ZY}) becomes
\begin{equation}\label{ZT}
Z^{{a}} = (-1)^{m+1}\hat{T}^{(0)}_{12}(0).
\end{equation} 
Simultaneously, the jump relations (\ref{Yjump})  transforms into the relations,
\begin{equation}\label{T1}
T_{+}(\lambda) = T_{-}(\lambda)\begin{pmatrix}e^{-h(\lambda)} &\omega(\lambda)\cr
                                     0 & e^{h(\lambda)}\end{pmatrix}, \quad \lambda \in \Sigma_0\setminus\{0\},
\end{equation}
where
\begin{equation}\label{h}   
h(\lambda) = g_{+}(\lambda) - g_{-}(\lambda)                                   
= m\log\frac{1 + \sqrt{\lambda}}{1-\sqrt{\lambda}} + n\log\frac{i + \sqrt{\lambda}}{i-\sqrt{\lambda}}, \quad \lambda \in \Sigma_0.
\end{equation} 

Put (cf. (\ref{H0gamma1})
\begin{equation}\label{Hdef}
H(\lambda) := \exp h(\lambda) = \left(\frac{1+\sqrt{\lambda}}{1-\sqrt{\lambda}}\right)^{m}
\left(\frac{i+\sqrt{\lambda}}{i-\sqrt{\lambda}}\right)^{n}.
\end{equation}
Observe, that this function admits an analytical continuation on ${\Bbb C}\setminus [0, i\infty]$.
Indeed, the continuation is given by formula (\ref{Hdef}) itself  with $\sqrt{\lambda}$ defined 
on the $\lambda$-plane with the cut along  $[0, i\infty)$
and the branch is fixed by the condition $-3\pi/2 < \arg \lambda < \pi/2$. We remind that in the case of the functions
$g(\lambda)$ and $\lambda^{{a}/4}$ the cut for  $\sqrt{\lambda}$ is $[0, -i\infty)$ and
$-\pi/2 < \arg \lambda < 3\pi/2$. The function $H(\lambda)$  has a pole at $\lambda = 1$ and a zero at $\lambda = -1$.
In what follows, a crucial role will be played by the following lemma.

\begin{lemma}\label{lemma1} For all $m, n > 0$, the positive  function $|H(\lambda)|$ is greater than $1$ in the first quadrant, 
and it is less than one in the second quadrant.
\end{lemma}
{\it Proof.} Follows immediately from the simple geometric fact that $|1+\sqrt{\lambda}| > |1-\sqrt{\lambda}|$ and
$|i+\sqrt{\lambda}| < |i-\sqrt{\lambda}|$ if $\lambda$ lies in the first quadrant. If $\lambda$ lies in the second
quadrant the inequalities are reversed. 

\subsection{Opening of lenses  and the second transformation $T \to S$ }
\begin{figure}[h]
  \begin{center}
     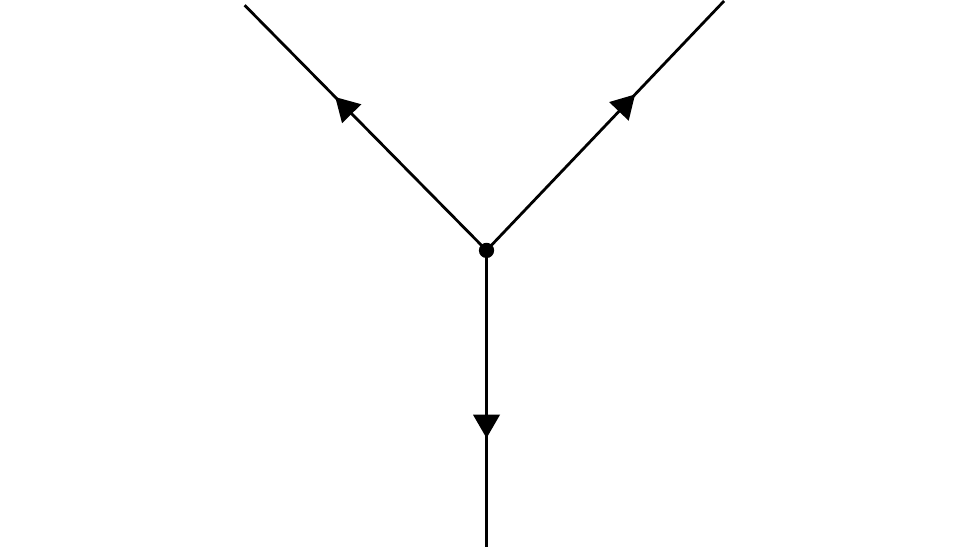
  \end{center}
  \caption{Contour for the $S$ - RH problem}\label{f.contour_S}
\end{figure}

As it is usual at this stage of implementation of the nonlinear steepest
descent method, we observe that
\begin{equation}\label{trianfactor}
\begin{pmatrix}e^{-h(\lambda)} &\omega(\lambda)\cr
                                     0 & e^{h(\lambda)}\end{pmatrix} =
\begin{pmatrix}1 &0\cr
                                     e^{h(\lambda)}\omega^{-1}(\lambda) & 1\end{pmatrix} 
 \begin{pmatrix}0 &\omega(\lambda)\cr
                                     -\omega^{-1}(\lambda) & 0\end{pmatrix} 
  \begin{pmatrix}1 &0\cr
                                     e^{-h(\lambda)}\omega^{-1}(\lambda) & 1\end{pmatrix} 
\end{equation}
and go from the function $T(\lambda)$ to the function $S(\lambda)$ defined by the
equations,
\begin{equation}\label{Sdef}
S(\lambda ) = T(\lambda)\begin{cases}
\begin{pmatrix}1 &0\cr
                                     -H^{-1}(\lambda)\omega^{-1}_{1}(\lambda) & 1\end{pmatrix}  &\lambda \in \Omega_r,\cr\cr
\begin{pmatrix}1 &0\cr
                                     H(\lambda)\omega^{-1}_{2}(\lambda) & 1\end{pmatrix}  &\lambda \in \Omega_l,\cr\cr
I & \lambda \notin \Omega_r\cup\Omega_l,
  \end{cases}
\end{equation}
where $\Omega_r$ ($\Omega_l$) is the region in the right (left) half-plane between the  rays 
$\Sigma_0$ and $\Sigma_1 = \{ \lambda: \Re \lambda = c\Im \lambda, c >0\}$($ \Sigma_2 = \{ \lambda: \Re \lambda = -c\Im \lambda, c >0\}$). 
The rays $\Sigma_1$ and $\Sigma_2$, similar to the ray $\Sigma_0$, are oriented toward infinity.
The functions $\omega_1$ and $\omega_2$ are given by the equations,
\begin{equation}\label{omega12}
\omega_1(\lambda) = 2i\sin{\frac{{a}\pi}{2}}\lambda^{-\frac{{a}}{2}},
\quad\mbox{and}\quad 
\omega_2(\lambda) = 2i\sin{\frac{{a}\pi}{2}}\lambda^{-\frac{{a}}{2}}e^{\pi i {a}}
\end{equation}
The regions $\Omega_l$ and $\Omega_r$ are depicted in Figure \ref{f.contour_S}.  The Riemann-Hilbert problem in
terms of the function $S(\lambda)$ reads as follows. 

\begin{itemize}
\item $S(\lambda)$  is analytic on $\Bbb{C}\setminus\Gamma$, $\Gamma = \Sigma_0\cup\Sigma_1\cup\Sigma_2$,
\item  The jump conditions are described by the equations,
\begin{enumerate}
\item as $\lambda \in \Sigma_0$,
\begin{equation}\label{SSigma0}
S_+(\lambda) =S_-(\lambda)
\begin{pmatrix}0 &\omega(\lambda)\cr
                                     -\omega^{-1}(\lambda) & 0\end{pmatrix}
\end{equation} 
\item as $\lambda \in \Sigma_1$,
\begin{equation}\label{SSigma1}
S_+(\lambda) =S_-(\lambda)
 \begin{pmatrix}1 &0\cr
                                     H^{-1}(\lambda)\omega^{-1}_{1}(\lambda) & 1\end{pmatrix} 
\end{equation}
\item as $\lambda \in \Sigma_2$,
\begin{equation}\label{SSigma2}
S_+(\lambda) =S_-(\lambda)
 \begin{pmatrix}1 &0\cr
                                     H(\lambda)\omega^{-1}_{2}(\lambda) & 1\end{pmatrix} 
\end{equation}
\end{enumerate}
\item as $\lambda \to \infty$,
\begin{equation}\label{Sinfty}
S(\lambda) =   \begin{pmatrix}\bullet&0\cr
 \bullet&\bullet\end{pmatrix}\Bigl(I + O\left(\frac{1}{\lambda}\right)\Bigr)
 \begin{pmatrix}1&e^{\frac{i\pi{a}}{2}}\lambda^{-\frac{{a}}{2}}
\cr\cr
 0&1\end{pmatrix}
 \end{equation}
$$
\times
(\lambda - 1)^{\frac{m}{2}\sigma_3}(\lambda+1)^{\frac{n}{2}\sigma_3}
 e^{-g(\lambda)\sigma_3}
 \begin{cases}
\begin{pmatrix}1 &0\cr
                                     -H^{-1}(\lambda)\omega^{-1}_{1}(\lambda) & 1\end{pmatrix}  &\lambda \in \Omega_r,\cr\cr
\begin{pmatrix}1 &0\cr
                                     H(\lambda)\omega^{-1}_{2}(\lambda) & 1\end{pmatrix}  &\lambda \in \Omega_l,\cr\cr
I & \lambda \notin \Omega_r\cup\Omega_l,
  \end{cases}
$$
\item as $\lambda \to 0$,
\begin{equation}\label{Szero}
S(\lambda) =  \begin{pmatrix}1&\bullet\cr
 0&1\end{pmatrix}\begin{pmatrix}1&0\cr
 0&(-1)^{m}\end{pmatrix}\Bigl(I + O(\lambda)\Bigr) 
\end{equation}
$$
 \begin{pmatrix}1&e^{\frac{i\pi{a}}{2}}\lambda^{-\frac{{a}}{2}}
 (\lambda - 1)^{-m}(\lambda+1)^{-n}
\cr\cr
 0&1\end{pmatrix}e^{-g(\lambda)\sigma_3}
 \begin{cases}
\begin{pmatrix}1 &0\cr
                                     -H^{-1}(\lambda)\omega^{-1}_{1}(\lambda) & 1\end{pmatrix}  &\lambda \in \Omega_r,\cr\cr
\begin{pmatrix}1 &0\cr
                                     H(\lambda)\omega^{-1}_{2}(\lambda) & 1\end{pmatrix}  &\lambda \in \Omega_l,\cr\cr
I & \lambda \notin \Omega_r\cup\Omega_l,
  \end{cases}
$$
\end{itemize} 

Let $U_{\delta} = \left\{|\lambda| < \delta < 1\right\}$ denote a small neighborhood of $\lambda = 0$. In this neighborhood, with the cut along the part of the ray $[0,-i\infty]$,
one can define the holomorphic function,
\begin{equation}\label{h0def}
h_0(\lambda) = m\log\frac{1 + \sqrt{\lambda}}{1-\sqrt{\lambda}} + n\log\frac{i + \sqrt{\lambda}}{i-\sqrt{\lambda}}
= 2(m-in)\sqrt{\lambda} + \sum_{k=0}^{\infty}a_k\lambda^{k +  \frac{1}{2}}.
\end{equation}
The function $h_0(\lambda)$ satisfies the following properties,
\begin{equation}\label{h0H}
\exp h_0(\lambda) = 
 \begin{cases}
H(\lambda) & \lambda \in U_{\delta},\,\,\, \Re \lambda >0,\cr\cr
H^{-1}(\lambda)&\lambda \in U_{\delta}, \,\,\,\Re \lambda <0.
  \end{cases}
\end{equation}
Moreover, one can also observe that, for all $\lambda \in U_{\delta}\cap[0,-i\infty)$,
\begin{equation}\label{gh1}
g(\lambda) - \frac{1}{2}h_0(\lambda) = \frac{m}{2}\log(1-\lambda) + \frac{n}{2}\log(1+\lambda) +\frac{i\pi}{2}n,
\end{equation}
where the branches, which are  holomorphic in $U_{\delta}$, are considered for the   logarithms in the right hand side. 

Similarly,
in the neighborhood $U_{1/\delta} = \left\{|\lambda| > 1/\delta > 1\right\}$ of the point $\lambda = \infty$, 
we can define the function $h_{\infty}(\lambda)$ ,
\begin{equation}\label{h0def2}
h_{\infty}(\lambda) = m\log\frac{1 + \sqrt{\lambda}}{1-\sqrt{\lambda}} + n\log\frac{i + \sqrt{\lambda}}{i-\sqrt{\lambda}}
= -i\pi(m+n) + 2(m+in)\frac{1}{\sqrt{\lambda}} + \sum_{k=0}^{\infty}b_k\lambda^{-k -  \frac{1}{2}}.
\end{equation}
In the neighborhood $U_{1/\delta}$, the function $h_{\infty}(\lambda)$ satisfies the properties similar to that
of $h_0(\lambda)$, i.e.
\begin{equation}\label{hinftyH}
\exp h_{\infty}(\lambda) = 
 \begin{cases}
H(\lambda) & \lambda \in U_{1/\delta},\,\,\, \Re \lambda >0,\cr\cr
H^{-1}(\lambda)& \lambda \in U_{1/\delta} \,\,\,\Re \lambda <0,
  \end{cases}
\end{equation}
and
\begin{equation}\label{gh2}
g(\lambda) - \frac{1}{2}h_{\infty}(\lambda) = \frac{m}{2}\log(1-\lambda) + \frac{n}{2}\log(1+\lambda) +\frac{i\pi}{2}n,
\end{equation}
for all $\lambda \in U_{1/\delta}\cap[0,-i\infty)$,

Equations (\ref{h0H} - \ref{gh1}) and (\ref{hinftyH} - \ref{gh2})  allow us to reformulate the 
$S$ - Riemann-Hilbert problem in the following, more compact way.

\vskip .2in
\noindent
{\bf Riemann-Hilbert problem for $S(\lambda)$}

\begin{itemize}
\item $S(\lambda)$  is analytic on $\Bbb{C}\setminus\Gamma$, $\Gamma = \Sigma_0\cup\Sigma_1\cup\Sigma_2$,
\item  The jump conditions are described by the equations,
\begin{enumerate}
\item as $\lambda \in \Sigma_0$,
\begin{equation}\label{SSigma00}
S_+(\lambda) =S_-(\lambda)
\begin{pmatrix}0 &\omega(\lambda)\cr
                                     -\omega^{-1}(\lambda) & 0\end{pmatrix}
\end{equation} 
\item as $\lambda \in \Sigma_1$,
\begin{equation}\label{SSigma11}
S_+(\lambda) =S_-(\lambda)
 \begin{pmatrix}1 &0\cr
                                     H^{-1}(\lambda)\omega^{-1}_{1}(\lambda) & 1\end{pmatrix} 
\end{equation}
\item as $\lambda \in \Sigma_2$,
\begin{equation}\label{SSigma22}
S_+(\lambda) =S_-(\lambda)
 \begin{pmatrix}1 &0\cr
                                     H(\lambda)\omega^{-1}_{2}(\lambda) & 1\end{pmatrix} 
\end{equation}
\end{enumerate}
\item as $\lambda \to \infty$,
\begin{equation}\label{Sinfty0}
S(\lambda) = \begin{pmatrix}\bullet&0\cr
 \bullet&\bullet\end{pmatrix}\left(I +O\left(\frac{1}{\lambda}\right)\right)
 \lambda^{-\frac{{a}}{4}\sigma_3}C\lambda^{\frac{{a}}{4}\sigma_3}e^{-\frac{1}{2}h_{\infty}(\lambda)\sigma_3},
\end{equation}
where the matrix $C$ is a piece-wise constant  matrix - valued function defined by the equations,
\begin{equation}\label{Cdef}
C = \begin{pmatrix}1 & e^{\frac{i\pi{a}}{2}}\cr
0&1\end{pmatrix}
 \begin{cases}
\begin{pmatrix}1 &0\cr
                                     -\omega^{-1}_{0} & 1\end{pmatrix}  &\lambda \in \Omega_r,\cr\cr
\begin{pmatrix}1 &0\cr
                                     \omega^{-1}_{0}e^{-i\pi{a}} & 1\end{pmatrix}  &\lambda \in \Omega_l,\cr\cr
I & \lambda \notin \Omega_r\cup\Omega_l,
  \end{cases}
\end{equation}
and
\begin{equation}\label{omega0}
\omega_0 = 2i\sin\frac{{a}\pi}{2},
\end{equation}
\item as $\lambda \to 0$,
\begin{equation}\label{Szero0}
S(\lambda) =  \begin{pmatrix}1&\bullet\cr
 0&1\end{pmatrix}\begin{pmatrix}1&0\cr
 0&(-1)^{m}\end{pmatrix}e^{-\frac{i\pi}{2}n\sigma_3}\Bigl(I + O(\lambda)\Bigr)
 \lambda^{-\frac{{a}}{4}\sigma_3}C\lambda^{\frac{{a}}{4}\sigma_3}e^{-\frac{1}{2}h_0(\lambda)\sigma_3}, 
\end{equation}
where $C$ is defined by the same equations (\ref{Cdef}).
\end{itemize} 
In terms of the $S$ - RH problem, the discrete $Z^{{a}}$ is given by the equation,
\begin{equation}\label{ZgammaS}
Z^{{a}} =(-1)^{m+1}e^{-\frac{i\pi}{2}n}\hat{S}^{(0)}_{12}(0),
\end{equation}
where 
\begin{equation}\label{Shatzero}
\hat{S}^{(0)}(\lambda)\equiv S(\lambda)e^{\frac{1}{2}h_0(\lambda)\sigma_3}
\lambda^{-\frac{{a}}{4}\sigma_3}C^{-1}\lambda^{\frac{{a}}{4}\sigma_3}
\end{equation}
--  is the holomorphic (and invertible) matrix factor in the right hand
side of (\ref{Szero0}).  Similar factor in the right hand side of (\ref{Sinfty0}) we shall denote
$\hat{S}^{(\infty)}(\lambda)$, i.e.
\begin{equation}\label{Shatinfty}
\hat{S}^{(\infty)}(\lambda)\equiv S(\lambda)e^{\frac{1}{2}h_{\infty}(\lambda)\sigma_3}
\lambda^{-\frac{{a}}{4}\sigma_3}C^{-1}\lambda^{\frac{{a}}{4}\sigma_3}
\end{equation}

\vskip .4in 
The $S$ - Riemann-Hilbert problem is depicted in Figure 3. We can completely switch to this Riemann-Hilbert
problem in our analysis of the discrete conformal map $Z^{{a}}$. That is, in addition to Theorem \ref{theorem3},
we can formulate the following theorem. 

\begin{theorem}\label{theorem4} Let $S(\lambda)$ be the matrix valued function defined by the  discrete 
conformal map $Z^{{a}}$ according to the equations (\ref{Sdef}), (\ref{Tdef}),
(\ref{Y}), (\ref{X}), (\ref{Phi}), and (\ref{Psinmdef}). Then, the function $S(\lambda)$ is the unique solution 
of the Riemann-Hilbert factorization problem (\ref{SSigma00})  -- (\ref{Szero0}). The map $Z^{{a}}$ itself can be
recovered from the known function $S$ by relation (\ref{ZgammaS}).
\end{theorem}

\begin{remark}\label{rem2}
It is worth noticing that, since $0< {a} < 2$, the condition at $\lambda =0$ can be {\it a priori} relaxed. Indeed,
it is enough to demand that
\begin{equation}\label{Szero1}
S(\lambda) =  O(1)\lambda^{-\frac{{a}}{4}\sigma_3}O(1)\lambda^{\frac{{a}}{4}\sigma_3}, \quad \lambda\to 0.
\end{equation}
More detail behavior at $\lambda = 0$ which is featured in (\ref{Szero0}) will be then a consequence of (\ref{Szero1}) and
the jump relations. (Of course, one still needs to formulate properly the normalization condition at $\lambda =0$).
\end{remark}

Now, we can highlight the role of Lemma \ref{lemma1}. Due to this lemma,  as $m, n \to \infty$,  the jump matrices across the rays $\Sigma_1$ and $\Sigma_2$ become
exponentially closed to the identity matrix away from the points $\lambda =0$ and $\lambda =\infty$, and hence the 
$S$ - problem is getting {\it localized.} This suggests that the approximate solution 
of the $S$ - RH problem can be assembled  from the two {\it local  parametrices} - the
solutions  of the local Riemann-Hilbert problems at $\lambda = 0$  and $\lambda = \infty$,
and  the  {\it global parametrix}  -   the solution  of the Riemann-Hilbert  problem 
associated with the jump across the ray $\Sigma_0$. In the next three subsections we will 
construct these three parametrices explicitly, and in subsection \ref{sec3.6} we will  assemble them  into the piece-wise
analytic matrix valued function,  which we will denote $S^{(as)}(\lambda)$.
In subsection \ref{sec3.6}, we will  also show  that $S^{(as)}(\lambda)$ is indeed a {\it parametrix} for the solution
of the full $S$ -problem, i.e., that the matrix quotient, $R(\lambda) := S(\lambda)[S^{(as)}(\lambda)]^{-1}$
solves the Riemann-Hilbert problem whose jump matrices are uniformly closed to the identity.
The last fact, by the general  arguments of the  Riemann-Hilbert theory, will allow us to prove that
$ S^{(as)}(\lambda)$ is the genuine asymptotic solution of the $S$-Riemann-Hilbert problem.
The just described strategy is standard for the nonlinear steepest descent method. The difference 
comparing with the usual situation is technical - the global parametrix is not, simultaneously,  the parametrix at the infinity,
as it happens in the usual  applications of the Riemann-Hilbert method. The reason lies in the Fuchsian origin
of the Riemann-Hilbert problem we are dealing with.

We shall start with the construction of the global parametrix.

\subsection{Global parametrix.} 

The  global parametrix for the solution of the $S$- RH problem,
which we will denote $P^{(gl)}(\lambda)$,  is
defined as the solution of the following Riemann-Hilbert problem posed on the
ray $\Sigma_0$.

\begin{itemize}
\item $P^{(gl)}(\lambda)$  is analytic on $\Bbb{C}\setminus\Sigma_0$,
\item  The jump conditions are described by the equation,
\begin{equation}\label{Pgljump}
P^{(gl)}_+(\lambda) =P^{(gl)}_-(\lambda)
\begin{pmatrix}0 &\omega(\lambda)\cr
                                     -\omega^{-1}(\lambda) & 0\end{pmatrix},\quad \lambda \in \Sigma_0\setminus\{0\}
\end{equation} 
%\item as $\lambda \to \infty$,
%\begin{equation}\label{Pglinfty}
%P^{(gl)}
%\end{equation}
\end{itemize} 
We note that in the setting of the $P^{(gl)}$ - RH problem we do not prescribe any
special behavior either at $\lambda = 0$ or at $\lambda = \infty$. Hence the parametrix
$P^{(gl)}(\lambda)$ is defined up to the left multiplication by the matrix valued function
analytic on $\Bbb{C}\setminus\{0\}$. This non-uniquiness, however, will not affect
the construction of the approximate solution to the $S$ - RH problem.

A solution of the $P^{(gl)}$ - RH problem can be easily found. Indeed, 
we notice that for all $\lambda \in \Sigma_0\setminus\{0\}$,
$$
\begin{pmatrix}0 &\omega(\lambda)\cr
                                     -\omega^{-1}(\lambda) & 0\end{pmatrix}
= \lambda_{-}^{-\frac{{a}}{4}\sigma_3}
\begin{pmatrix}0 &2i\sin\frac{{a}\pi}{2}e^{\frac{\pi i}{2}{a}}\cr
  -\frac{1}{2i\sin\frac{{a}\pi}{2}}e^{-\frac{\pi i}{2}{a}} & 0\end{pmatrix}
\lambda_{+}^{\frac{{a}}{4}\sigma_3}  
$$
\begin{equation}\label{glob1}
= \lambda_{-}^{-\frac{{a}}{4}\sigma_3}
\begin{pmatrix}0 &e^{i\pi{a}} -1\cr
  -\frac{1}{e^{i\pi{a}} -1} & 0\end{pmatrix}
\lambda_{+}^{\frac{{a}}{4}\sigma_3}
= \lambda_{-}^{-\frac{{a}}{4}\sigma_3}\eta^{\sigma_3}
\begin{pmatrix}0 &1\cr
  -1 & 0\end{pmatrix}\eta^{-\sigma_3}
\lambda_{+}^{\frac{{a}}{4}\sigma_3},
\end{equation}
where $\eta=\sqrt{e^{\pi i{a}} -1}$. Diagonalizing the matrix $\begin{pmatrix}0 &1\cr
  -1 & 0\end{pmatrix}$, we also have that
\begin{equation}\label{glob2}
\begin{pmatrix}0 &1\cr
  -1 & 0\end{pmatrix} = T^{-1}i\sigma_3 T = T^{-1}\lambda_{-}^{\frac{\sigma_3}{4}}\lambda_{+}^{-\frac{\sigma_3}{4}}T,
  \quad\mbox{where}\quad
  T= \begin{pmatrix}\frac{1}{2} &- \frac{i}{2}\cr\cr
 \frac{1}{2} & \frac{i}{2}\end{pmatrix}.
 \end{equation}
Combining (\ref{glob1}) and (\ref{glob2}) we arrive at the following representation for the jump
matrix of the $P^{(gl)}$ - RH problem.
$$
\begin{pmatrix}0 &\omega(\lambda)\cr
                                     -\omega^{-1}(\lambda) & 0\end{pmatrix}
= \lambda_{-}^{-\frac{{a}}{4}\sigma_3}\eta^{\sigma_3}
T^{-1}\lambda_{-}^{\frac{\sigma_3}{4}}\lambda_{+}^{-\frac{\sigma_3}{4}}T
\eta^{-\sigma_3}
\lambda_{+}^{\frac{{a}}{4}\sigma_3}.
$$
This equation suggests that the the global parametric can be taking in the form,
\begin{equation}\label{Sgl}
P^{(gl)}(\lambda) =   \lambda_{+}^{-\frac{\sigma_3}{4}}T
\eta^{-\sigma_3}
\lambda_{+}^{\frac{{a}}{4}\sigma_3} 
=\lambda^{-\frac{\sigma_3}{4}}
\begin{pmatrix}\frac{1}{2} &- \frac{i}{2}\cr\cr
 \frac{1}{2} & \frac{i}{2}\end{pmatrix}\eta^{-\sigma_3}\lambda^{\frac{{a}}{4}\sigma_3},\quad
 \eta=\sqrt{e^{\pi i{a}} -1}.
\end{equation}

We shall now concentrate on constructing the parametrices to the solution
of the $S$ - problem at points $\lambda =0$ and $\lambda =\infty$

\subsection{Parametrix at $\lambda =0$.}\label{sect3.4}

Expansion (\ref{h0def}) implies that in the neighborhood $U_{\delta}$,
\begin{equation}\label{h0series}
h^{2}_0(\lambda) = 4(m-in)^2\Bigl(\lambda + \sum_{j \geq 2}c_j \lambda^j\Bigr),
\end{equation} 
where the coefficients  $c_j$ satisfies the uniform estimate,
\begin{equation}\label{cjest}
|c_k| \leq \frac{c}{k}, \quad k >1, \quad n,m > 0.
\end{equation}
Therefore, the equation,
\begin{equation}\label{xidef}
\xi(\lambda) = h^{2}_0(\lambda) \equiv 4(m-in)^2\Bigl(\lambda + \sum_{j \geq 2}c_j \lambda^j\Bigr),
\end{equation} 
determines a conformal change of variables in the neighborhood  $U_{\delta}$:
\begin{equation}\label{ximap}
U_{\delta} \to D_{r}(0) \equiv \{\xi: |\xi| < r^2\delta\}, \quad r = \sqrt{m^2 + n^2}.
\end{equation}

\begin{figure}[h]
  \begin{center}
     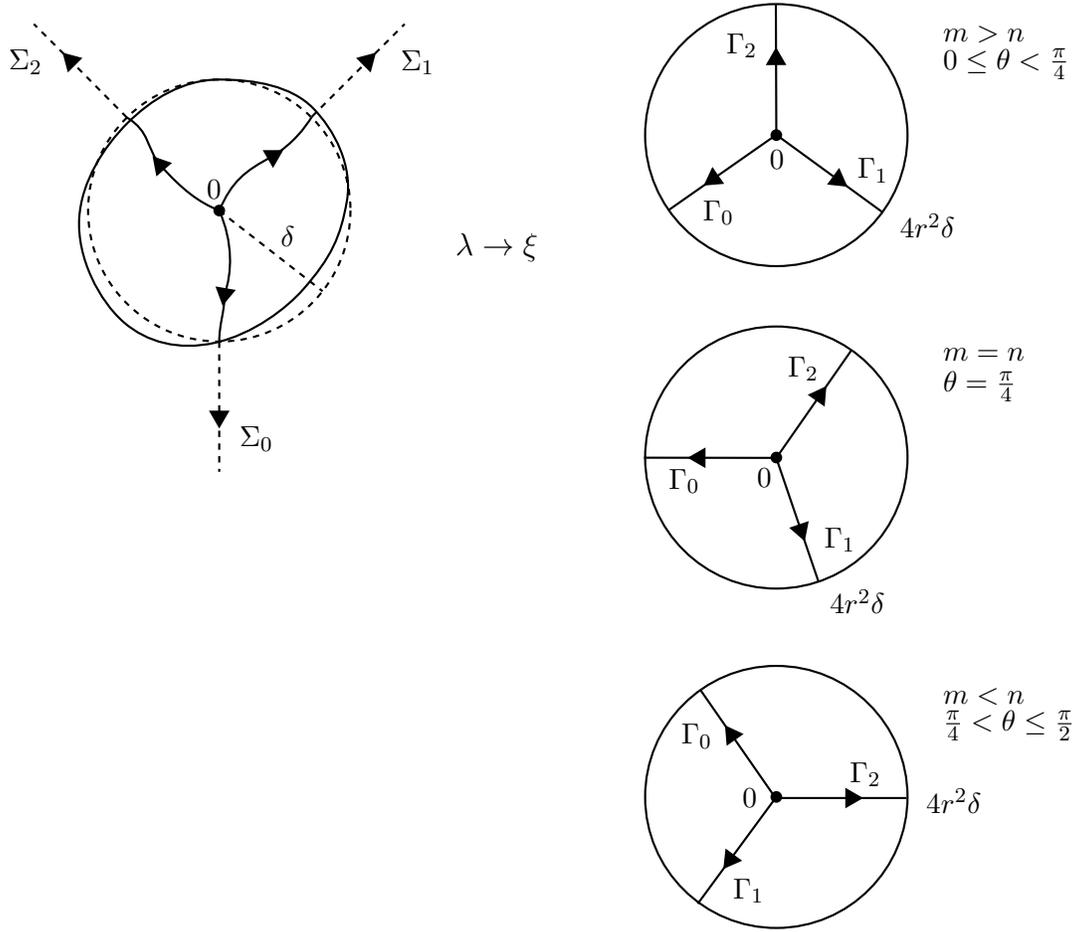
  \end{center}
  \caption{The local map $\lambda \rightarrow \xi$}\label{f.local_map}
\end{figure}

The action of the map $\lambda \to \xi$ on the part of the contour  $\Gamma$ 
of the $S$ - RH problem, which is inside of the neighborhood $U_{\delta}$ is indicated
in Figure \ref{f.local_map}. We shall assume that the rays $\Sigma_k$  are actually
slightly deformed so that   inside of the neighborhood $U_{\delta}$ they coincide
with the pre-images of the rays $\Gamma_k$ which satisfy the following conditions,
\begin{equation}\label{Gamma0}
\arg\xi|_{\Gamma_{0+}} = -\frac{\pi}{2} -2\theta,
\end{equation}
\begin{equation}\label{Gamma12}
\arg\xi|_{\Gamma_1} = \frac{\pi}{4} - 2\theta, \quad
\arg\xi|_{\Gamma_2} = \frac{3\pi}{4} - 2\theta,
\end{equation}
where
\begin{equation}\label{thetadef}
\theta = - \arg (m -i n), \quad 0\leq \theta \leq \frac{\pi}{2}.
\end{equation}
Define $\sqrt{\xi}$ on the $\xi$-plane cut along the ray $\Gamma_0$ and fixed
by the condition, 
$$
-\frac{\pi}{2} -2\theta < \arg\xi < \frac{3\pi}{2} -2\theta.
$$
Then, we will have that, for all $\theta$,
\begin{equation}\label{epsilon0}
-\frac{\pi}{2} + \frac{\pi}{8} \leq  \arg\sqrt{\xi}|_{\Gamma_{1}}<   \arg\sqrt{\xi}|_{\Gamma_{2}}
\leq \frac{\pi}{2} - \frac{\pi}{8},
\end{equation}
that is, for all $\theta$, the images, under
the map $\xi \mapsto \sqrt{\xi}$,  of the rays $\Gamma_1$ and $\Gamma_2$ and of the
sector between them lie in the right half plane $\Re \sqrt{\xi} > 0$.

\begin{figure}[h]
  \begin{center}
     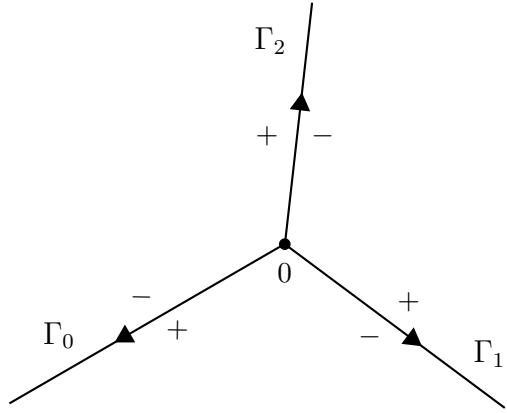
  \end{center}
  \caption{The contour for the $\Phi^{(0)}$ - RH problem}\label{f.contour_Phi}
\end{figure}

Observe that inside of the  neighborhood $U_{\delta}$, cut along the curve $\Sigma_0$, we have that
 $h_0(\lambda) = \sqrt{\xi(\lambda)}$. Therefore, the jump matrix $G_{S}$ of the $S$ - RH problem  inside 
of the  neighborhood $U_{\delta}$ can be written down in the form,
\begin{equation}\label{Sjump000}
G_{S}(\lambda) = e^{ \frac{1}{2} \sqrt{\xi(\lambda)}\sigma_3} \lambda^{-\frac{{a}}{4}\sigma_3}L
 \lambda^{\frac{{a}}{4}\sigma_3} e^{-\frac{1}{2}\sqrt{\xi(\lambda)}\sigma_3},
\end{equation}
where the piecewise constant matrix $L$ is given by the equations,
$$
L =  \begin{pmatrix}0 &\omega_0e^{\frac{i\pi{a}}{2}}\cr
                                     -\omega^{-1}_0e^{-\frac{i\pi{a}}{2}} & 0\end{pmatrix}\equiv L_0,\quad 
\lambda \in \Sigma_0\cap U_{\delta},
$$
$$
L =  \begin{pmatrix}1 &0\cr
                                     \omega^{-1}_0 & 1\end{pmatrix} \equiv L_1,\quad 
\lambda \in \Sigma_1\cap U_{\delta},
$$
$$
L =  \begin{pmatrix}1 &0\cr
                                     e^{-\pi i{a}}\omega^{-1}_0 & 1\end{pmatrix}\equiv L_2,\quad 
\lambda \in \Sigma_2\cap U_{\delta},
$$
and $\omega_0$ is defined in (\ref{omega0}).
Therefore, the map, $\lambda \to \xi$, transforms   the $U_{\delta}$ - part  of the 
$S$ - RH problem into the following model RH problem which is formulated for 
a matrix function $\Phi^{(0)}(\xi)$ defined on the $\xi$-plane.
\begin{itemize}
\item $\Phi^{(0)}(\xi)$  is analytic on $\Bbb{C}\setminus\Gamma_{\xi}$, $\Gamma_{\xi} = \Gamma_0
\cup\Gamma_1\cup\Gamma_2$,
\item  The jump conditions are described by the equations,
\begin{equation}\label{Phijamp}
\Phi^{(0)}_+ = \Phi_{-}^{(0)}e^{ \frac{1}{2} \sqrt{\xi}\sigma_3} \xi^{-\frac{{a}}{4}\sigma_3}L
 \xi^{\frac{{a}}{4}\sigma_3} e^{-\frac{1}{2}\sqrt{\xi}\sigma_3},
\end{equation}
where $L = L_k$ if $\xi \in \Gamma_k$, $k = 0, 1, 2$.
\item as $\xi \to \infty$,
\begin{equation}\label{Phiinfty}
\Phi^{(0)}(\xi) = \xi^{-\frac{1}{4}\sigma_3}\begin{pmatrix}\frac{1}{2} &- \frac{i}{2}\cr\cr
 \frac{1}{2} & \frac{i}{2}\end{pmatrix}
 \left(I + O\left(\frac{1}{\sqrt{\xi}}\right)\right)\eta^{-\sigma_3}\xi^{\frac{{a}}{4}\sigma_3}
 \end{equation}
\item as $\xi \to 0$,
\begin{equation}\label{Phizero}
\Phi^{(0)}(\xi) =  \hat{\Phi}^{(0)}(\xi)\xi^{-\frac{{a}}{4}\sigma_3}C\xi^{\frac{{a}}{4}\sigma_3}e^{-\frac{1}{2}\sqrt{\xi}\sigma_3}
\end{equation}
where the matrix-valued functions $\hat{\Phi}^{(0)}(\xi)$ is holomorphic at $\xi = 0$, 
\begin{equation}\label{B0}
\hat{\Phi}^{(0)}(\xi)= B\Bigl( I + O(\xi)\Bigr), \quad \det\hat{\Phi}^{(0)}(\xi) \equiv \frac{i}{2},
\end{equation}
and the piece-wise constant matrix - valued function $C$ is the same as in (\ref{Cdef}) with 
$\Omega_{r,l}$ replaced by their images via the map $\lambda \mapsto \xi(\lambda)$, i.e.
\begin{equation}\label{CPhidef}
C = \begin{pmatrix}1 & e^{\frac{i\pi{a}}{2}}\cr
0&1\end{pmatrix}
 \begin{cases}
\begin{pmatrix}1 &0\cr
                                     -\omega^{-1}_{0} & 1\end{pmatrix}  &-\frac{\pi}{2}-2\theta < \arg\xi < \frac{\pi}{4} -2\theta,\cr\cr
\begin{pmatrix}1 &0\cr
                                     \omega^{-1}_{0}e^{-i\pi{a}} & 1\end{pmatrix}  &\frac{3\pi}{4}-2\theta < \arg\xi < \frac{3\pi}{2} -2\theta,\cr\cr
I & \frac{\pi}{4}-2\theta < \arg\xi < \frac{3\pi}{4} -2\theta ,
  \end{cases}
\end{equation}

\end{itemize}
The branch of the function $\xi^{1/4}$ is define on the $\xi$-plane cut along the ray $\Gamma_0$
and fixed by the condition $\xi^{1/4} > 0$ as $\xi >0$, i.e.
\begin{equation}\label{xibranch}
-\frac{\pi}{2} - 2\theta < \arg{\xi} < \frac{3\pi}{2} - 2\theta.
\end{equation}
The problem is depicted in the Figure \ref{f.contour_Phi} (for the case $m > n$).

The same remark as in the case of the $S$-problem can be made, i.e., since $0 < {a} < 2$, in the setting
of the $\Phi^{(0)}$ - RH problem it is enough to demand, that
\begin{equation}\label{Phizero1}
\Phi^{(0)}(\xi) =  O(1)\xi^{-\frac{{a}}{4}\sigma_3}O(1)\xi^{\frac{{a}}{4}\sigma_3}, \quad \xi \to 0.
\end{equation}

The normalization condition (\ref{Phiinfty}) comes from the fact that we want the ``interior''  
function $\Phi^{(0)}(\xi(\lambda))$ to
match asymptotically, as $n, m \to \infty$,  the ``exterior''  function $P^{(gl)}(\lambda)$ 
at the boundary of $U_{\delta}$. In other words, to specify
the behavior of $\Phi^{(0)}(\xi)$ as $\xi \to \infty$, we must look at the behavior
of $P^{(gl)}(\lambda)$ at $\lambda = 0$. To this end, we notice that the function $P^{(gl)}(\lambda)$ 
can be written in the neighborhood $U_{\delta}$ in the form,
\begin{equation}\label{Sglat0}
P^{(gl)}(\lambda) = E(\lambda) \xi(\lambda)^{-\frac{1}{4}\sigma_3}\begin{pmatrix}\frac{1}{2} &- \frac{i}{2}\cr\cr
 \frac{1}{2} & \frac{i}{2}\end{pmatrix}
 \eta^{-\sigma_3}\lambda^{\frac{{a}}{4}\sigma_3},
\end{equation}
where
\begin{equation}\label{Edef}
E(\lambda) = \left(\frac{\xi(\lambda)}{\lambda}\right)^{\frac{1}{4}\sigma_3}
\end{equation}
is holomorphic at $\lambda = 0$. Indeed, in view of (\ref{xidef}),  we have that
\begin{equation}\label{Edef2}
E(\lambda) 
=\left(2(m-in)\right)^{\frac{1}{2}\sigma_3}\Bigl( 1 + \sum_{j\geq 1}c_j\lambda^j\Bigr)^{\frac{1}{4}\sigma_3}
= \left(2(m-in)\right)^{\frac{1}{2}\sigma_3}\Bigl( I + \sum_{j\geq 1}C_j\lambda^j\Bigr),
\end{equation}
where $C_j$ are (diagonal) matrix  coefficients of the Taylor series indicated.
Equation (\ref{Sglat0}) explains the choice of the normalization
condition at $\xi = \infty$ which we made in the model problem (\ref{Phijamp} - \ref{Phizero}).
The holomorphic factor $E(z)$ has no relevance to the setting of the 
Riemann-Hilbert problem in the $\xi$-plane; it will be restored latter on, when we 
start actually assembling the parametrix for $S(\lambda)$ in $U_{\delta}$.

\begin{remark}\label{Phi0uniq} It should be noticed that solution of the Riemann-Hilbert problem
(\ref{Phijamp}) - (\ref{Phizero}) is not unique and is defined up to the transformation,
\begin{equation}\label{gauge0}
\Phi^{(0)}(\xi) \to \begin{pmatrix} 1&0\cr\cr
\kappa&1\end{pmatrix}\Phi^{(0)}(\xi),
\end{equation}
where $\kappa$ is an arbitrary complex number. As with the setting of the Riemann-Hilbert
problem for the global parametrix, this non-uniqness will not
affect the construction of the approximate solution to the $S$ - RH problem. In fact, the uniqness can be formally
achieved if the error $O(\xi^{-1/2})$ in the normalization condition (\ref{Phiinfty})
is replaced by the error $O(\xi^{-1})$. However, as it follows from the explicit solution
of the problem, which is presented in Appendix A, this error can not be achieved for the generic
value of ${a}$. It also can be observed, that with the help of the gauge transformation (\ref{gauge0})
the normalization condition  (\ref{Phiinfty}) at infinity can be replaced by the condition,
\begin{equation}\label{Phiinfty0}
\Phi^{(0)}(\xi) =
 \left(I + O\left(\frac{1}{\xi}\right)\right)
 \xi^{-\frac{1}{4}\sigma_3}\begin{pmatrix}\frac{1}{2} &- \frac{i}{2}\cr\cr
 \frac{1}{2} & \frac{i}{2}\end{pmatrix}\eta^{-\sigma_3}\xi^{\frac{{a}}{4}\sigma_3},
\end{equation}
as $\xi \to \infty$. With this modification, the setting of the  Riemann-Hilbert problem for the function
 $\Phi^{(0)}(\xi)$ will provide  the uniqness property  of its  solution. We prefer, however, to 
stay with condition (\ref{Phiinfty}) and keep in mind the possibility of the gauge 
transformation (\ref{gauge0}).
\end{remark}

\begin{figure}[h]
  \begin{center}
     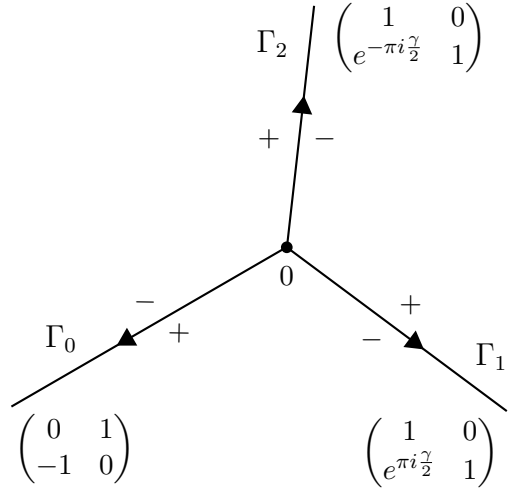
  \end{center}
  \caption{The contour and jump-matrices for the $\Psi^{(0)}$ - RH problem}\label{f.contour_Psi}
\end{figure}

Similar to the model problems appearing in \cite{DIZ} and \cite{Ku}, the model problem (\ref{Phijamp}) - (\ref{Phizero}) admits an explicit solution in
terms of the Bessel functions. In order to see this, let us make the following  simplifying 
substitution,
\begin{equation}\label{Psidef}
\Phi^{(0)}(\xi) = \Psi^{(0)}(\xi) \eta^{-\sigma_3}\xi^{\frac{{a}}{4}\sigma_3}e^{-\frac{1}{2}\sqrt{\xi}\sigma_3}.
\end{equation}
In terms of the function $\Psi^{(0)}(\xi)$, the Riemann-Hilbert problem (\ref{Phijamp}) - (\ref{Phizero}) 
reads:
\begin{itemize}
\item $\Psi^{(0)}(\xi)$  is analytic on $\Bbb{C}\setminus\Gamma_{\xi}$, $\Gamma_{\xi} = \Gamma_0
\cup\Gamma_1\cup\Gamma_2$\item  The jump conditions are described by the equations,
\begin{equation}\label{Psijamp}
\Psi^{(0)}_+ = \Psi_{-}^{(0)}L^{(0)}
\end{equation}
where $L^{(0)} = L^{0}_k$ if $\xi \in \Gamma_k$, $k = 0, 1, 2$, and
\begin{equation}\label{L0k}
L^{(0)}_0 = \eta^{-\sigma_3}L_0\eta^{\sigma_3}= \begin{pmatrix}0&1\cr\cr-1&0\end{pmatrix},
\quad 
L^{(0)}_{1,2} = \eta^{-\sigma_3}L_{1,2}\eta^{\sigma_3}= \begin{pmatrix}1&0\cr\cr e^{\pm\frac{\pi i{a}}{2}}&1\end{pmatrix}
\end{equation}
\item as $\xi \to \infty$,
\begin{equation}\label{Psiinfty}
\Psi^{(0)}(\xi) = \xi^{-\frac{1}{4}\sigma_3}\begin{pmatrix}\frac{1}{2} &- \frac{i}{2}\cr\cr
 \frac{1}{2} & \frac{i}{2}\end{pmatrix}
 \left(I + O\left(\frac{1}{\sqrt{\xi}}\right)\right)e^{\frac{1}{2}\sqrt{\xi}\sigma_3}
 \end{equation}
\item as $\xi \to 0$,
\begin{equation}\label{Psizero}
\Psi^{(0)}(\xi) =   B_0\Bigl( I + O(\xi)\Bigr)\xi^{-\frac{{a}}{4}\sigma_3}C_0
\end{equation}
where  $B_0$ is related with the matrix $B$ from  (\ref{B0}) by the relation,
\begin{equation}\label{B00}
B_0 =  B\eta^{\sigma_3},
\end{equation}
and the piece-wise constant matrix - valued function $C_0$ is defined by the equations, 
\begin{equation}\label{CPsidef}
C_0 = \begin{pmatrix}1 & \frac{1}{2i\sin{\frac{\pi{a}}{2}}}\cr
0&1\end{pmatrix}
 \begin{cases}
\begin{pmatrix}1 &0\cr
                                     -e^{\frac{i\pi{a}}{2}} & 1\end{pmatrix}  &-\frac{\pi}{2}-2\theta < \arg\xi < \frac{\pi}{4} -2\theta,\cr\cr
\begin{pmatrix}1 &0\cr
                                     e^{-\frac{i\pi{a}}{2}} & 1\end{pmatrix}  &\frac{3\pi}{4}-2\theta < \arg\xi < \frac{3\pi}{2} -2\theta,\cr\cr
I & \frac{\pi}{4}-2\theta < \arg\xi < \frac{3\pi}{4} -2\theta ,
  \end{cases}
\end{equation}
\end{itemize}
As before, the condition at $\xi =0$ can be replaced by 
\begin{equation}\label{Psizero0}
\Psi^{(0)}(\xi) =  O(1)\xi^{-\frac{{a}}{4}\sigma_3}O(1), \quad \xi \to 0.
\end{equation}
The problem is depicted in the Figure \ref{f.contour_Psi}.

A distinguished feature of this Riemann-Hilbert
problem is $\xi$ - independence of its jump matrices. Following the standard arguments (see e.g.
\cite{FIKN}) we derive from this fact that the ``logariphmic derivative''  of the solution
$\Psi^{(0)}(\xi))$ of the problem, 
$$
A(\xi) := \frac{d\Psi^{(0)}(\xi)}{d\xi}(\Psi^{(0)}(\xi))^{-1},
$$
is continious accross the contour $\Gamma_{\xi}$ and hence is analytic on ${\Bbb C} \setminus\{0\}$.
Moreover, the (differentiable in $\xi$ !) asymptotoc expansions (\ref{Psiinfty}) and
(\ref{Psizero}) tell us that,
\begin{equation}\label{Axiinfty}
A(\xi) = \frac{1}{4}\begin{pmatrix}0&0\cr\cr
1&0\end{pmatrix} + \frac{1}{4\xi}\begin{pmatrix}-1&1\cr\cr
0&1\end{pmatrix} +\xi^{-\frac{1}{4}\sigma_3}O\left(\frac{1}{\xi}\right)\xi^{\frac{1}{4}\sigma_3},
\end{equation}
as $\xi \to \infty$, and
\begin{equation}\label{Axizero}
A(\xi) = -\frac{{a}}{4\xi}B_0\sigma_3 B^{-1}_0.
\end{equation}
as $\xi \to 0$. Combaining these estimates with the analyticity of $A(\xi)$ on
${\Bbb C}\setminus\{0\}$, we arrive at the conclusion that 
 $A(\xi)$ is a rational function admiting the following
representation,
\begin{equation}\label{Axi0}
A(\xi) = \frac{1}{4}\begin{pmatrix}0&0\cr\cr
1&0\end{pmatrix} + \frac{1}{4\xi}\begin{pmatrix}\alpha&1\cr\cr
\beta&-\alpha\end{pmatrix},
\end{equation}
where $a$ and $b$ are some complex numbers satisfying (as it follows from (\ref{Axizero}))
the determinant constraint,
$$
\alpha^2+\beta = {a}^2.
$$
Using the gauge transformation (\ref{gauge0}) with $\kappa = -\alpha$ we can actually eliminate 
the diagonal entries of the matrix $A(\xi)$ and 
reduce $A(\xi)$ to the form,
\begin{equation}\label{Axi1}
A(\xi) = \frac{1}{4}\begin{pmatrix}0&0\cr\cr
1&0\end{pmatrix} + \frac{1}{4\xi}\begin{pmatrix}0&1\cr\cr
{a}^2&0\end{pmatrix}.
\end{equation}
Hence, the solution of the Riemann-Hilbert prioblem (\ref{Psijamp}) - (\ref{Psizero}),
if exists, can be choosen in such a way that it satisfies the matrix linear differential equation,
\begin{equation}\label{Axi2}
\frac{d\Psi^{(0)}(\xi)}{d\xi} = 
\frac{1}{4}\begin{pmatrix}0&\frac{1}{\xi}\cr\cr
1+\frac{{a}^2}{\xi}&0\end{pmatrix}\Psi^{(0)}(\xi).
\end{equation}
Put 
$$
\psi_1(\xi) := \Psi^{(0)}_{1j}(\xi),\quad \psi_2(\xi) := \Psi^{(0)}_{2j}(\xi),
$$
for $j =1$ or $j=2$. Then from (\ref{Axi2}) it follows that
\begin{equation}\label{psi2}
\psi_2(\xi) = 4\xi\frac{d\psi_1(\xi)}{d\xi},
\end{equation}
while the function $\psi_1(\xi)$ satisfies the second order linear ODE,
\begin{equation}\label{psi1}
\frac{d^2\psi_1 }{d\xi^2}+\frac{1}{\xi}\frac{d\psi_1}{d\xi}
-\frac{1}{16\xi}\left(1 + \frac{{a}^2}{\xi}\right)\psi_1=0.
\end{equation}
By the change of variables,
$$
z = \frac{i}{2}\sqrt{\xi}, \quad \psi_1(\xi) = y(z), 
$$
equation (\ref{psi1}) becomes the standart Bessel equation,
\begin{equation}\label{bessel}
\frac{d^2y }{dz^2}+\frac{1}{z}\frac{dy}{dz}
+\left(1 - \frac{{a}^2}{4z^2}\right)y=0.
\end{equation}
Therefore, for the solution of the model Riemann-Hilbert 
problem (\ref{Psijamp}) - (\ref{Psizero}) the following ansatz might be 
suggested,
$$
\Psi^{(0)}(\xi) = \begin{pmatrix}1&0\cr\cr
0&4\xi\end{pmatrix}
\begin{pmatrix}H^{(2)}_{-{a}/2}\left(\frac{i}{2}\sqrt{\xi}\right)&H^{(1)}_{-{a}/2}\left(\frac{i}{2}\sqrt{\xi}\right)\cr\cr
\frac{d}{d\xi}H^{(2)}_{-{a}/2}\left(\frac{i}{2}\sqrt{\xi}\right)&\frac{d}{d\xi}H^{(1)}_{-{a}/2}\left(\frac{i}{2}\sqrt{\xi}\right)\end{pmatrix}C^{(0)},
$$
where $H^{(1,2)}_{-{a}/2}(z)$ are the Hankel functions forming a basis for the solution space of (\ref{bessel}), 
and  $C^{(0)}$ is the constant matrix whose choice could depend on the sector on the $\xi$ - plane. 
Next proposition specifies exactly how the matrix $C^{(0)}$ should be chosen.
\begin{prop}\label{modanswer} The following
formulae define a solution of the problem (\ref{Psijamp}) - (\ref{Psizero}).
$$
\Psi^{(0)}(\xi) = \frac{\sqrt{\pi}}{2}\begin{pmatrix}\frac{1}{2}&0\cr\cr
0&2\xi\end{pmatrix}
\begin{pmatrix}H^{(2)}_{-{a}/2}\left(\frac{i}{2}\sqrt{\xi}\right)&H^{(1)}_{-{a}/2}\left(\frac{i}{2}\sqrt{\xi}\right)\cr\cr
\frac{d}{d\xi}H^{(2)}_{-{a}/2}\left(\frac{i}{2}\sqrt{\xi}\right)&\frac{d}{d\xi}H^{(1)}_{-{a}/2}\left(\frac{i}{2}\sqrt{\xi}\right)\end{pmatrix}
e^{\frac{\pi i{a}}{4}\sigma_3}
$$
\begin{equation}\label{Psi0form}
\times  \begin{cases}
I &-\frac{\pi}{2}-2\theta < \arg\xi < \frac{\pi}{4} -2\theta,\cr\cr
\begin{pmatrix}1 &0\cr
                                     2\cos\frac{\pi{a}}{2} & 1\end{pmatrix}  &\frac{3\pi}{4}-2\theta < \arg\xi < \frac{3\pi}{2} -2\theta,\cr\cr
\begin{pmatrix}1 &0\cr
                                     e^{\frac{i\pi{a}}{2}} & 1\end{pmatrix}  & \frac{\pi}{4}-2\theta < \arg\xi < \frac{3\pi}{4} -2\theta ,
\end{cases}
\end{equation}
which in addition satisfies the following specification of the asymptotic
condition (\ref{Psiinfty}),
\begin{equation}\label{Psispec}
\Psi^{(0)}(\xi) = \xi^{-\frac{1}{4}\sigma_3}
 \left(I + \frac{1}{\sqrt{\xi}}\Psi_1+O\left(\frac{1}{\xi}\right)\right)\begin{pmatrix}\frac{1}{2} &- \frac{i}{2}\cr\cr
 \frac{1}{2} & \frac{i}{2}\end{pmatrix}e^{\frac{1}{2}\sqrt{\xi}\sigma_3},
\end{equation}
where the constant matrix $\Psi_1$ is off-diagonal and is given by the equations,
\begin{equation}\label{Psi1}
\Psi_1 = \begin{pmatrix}0 & \psi_1\cr\cr
 \psi_1 -1 & 0\end{pmatrix}, \quad \psi_1 =\frac{1-{a}^2}{4}.
 \end{equation}
\end{prop}
The proof of the proposition is based on the known algebraic and asymptotics properties
of the Hankel functions and it is presented in detail in Appendix A.

Having constructed the function $\Psi^{(0)}(\xi)$ and hence the solution of the model problem
$\Phi^{(0)}(\xi)$ the local parametrix at the point $\lambda =0$ is defined by the equations,
\begin{equation}\label{P0def}
P^{(0)}(\lambda) =  E(\lambda)\Phi^{(0)}(\xi(\lambda))\left(\frac{\lambda}{\xi(\lambda)}\right)^{\frac{{a}}{4}\sigma_3} 
\end{equation}
Taking into account the holomorphicity of $E(\lambda)$ in $U_{\delta}$ (see (\ref{Edef})), we conclude that
inside of the neighborhood $U_{\delta}$, the function $P^{(0)}(\lambda)$
has exactly the same jumps as the  solution of the $S$-problem is supposed to have. Indeed, just as it is with the
function $E(\lambda)$, 
the right factor $\left(\frac{\lambda}{\xi(\lambda)}\right)^{\frac{{a}}{4}\sigma_3} $ is holomorphic in  $U_{\delta}$
and replaces the functions $\xi^{\pm \frac{{a}}{4}\sigma_3}$ in the $\Phi^{(0)}$ -  jump matrix (\ref{Phijamp})
by  the functions $\lambda^{\pm \frac{{a}}{4}\sigma_3}$. In this way,  the 
$\Phi^{(0)}$ -  jump matrix (\ref{Phijamp}) transforms into $S$-jump matrix (\ref{Sjump000}). By the same reason,
the singular factors of the right hand side of (\ref{Phizero}) transforms into the singular factors of the
right hand side of (\ref{Szero0}). In other words, if $S(\lambda)$ is the solution of the $S$ - problem, 
then
\begin{equation}\label{SP0}
S(\lambda)[P^{(0)}(\lambda)]^{-1} = \mbox{holomorphic function in  $U_{\delta}$}.
\end{equation}
At the same time, on the boundary of the neighborhood, $\xi(\lambda) \to \infty$ as $m, n \to \infty$. Therefore,
the function $\Psi^{(0)}(\xi(\lambda))$ can be replaced there by its asymptotics (\ref{Psispec}), and
we can see that  on the boundary of the neighborhood the following matching relation with the global
parametrix $P^{(gl)}(\lambda)$ takes place (cf. (\ref{Sglat0})).
\begin{equation}\label{match1}
P^{(0)}(\lambda) = \left(I + \frac{1}{\sqrt{\xi(\lambda)}}\lambda^{-\frac{1}{4}\sigma_3}
\Psi_1\lambda^{\frac{1}{4}\sigma_3} +  O\left(\frac{1}{r^2}\right)\right)P^{(gl)}(\lambda),
\quad \lambda \in \partial U_{\delta}, \quad n^2+m^2 \to \infty.
\end{equation}
We notice that again the factor $\left(\frac{\lambda}{\xi(\lambda)}\right)^{\frac{{a}}{4}\sigma_3}$
was important in bringing the leading asymptotic term of (\ref{Phiinfty}) to the form of (\ref{Sglat0}).

Near the point $\lambda =0$, the parametrix $P^{(0)}(\lambda)$ admits the representation (cf. (\ref{Szero})),
\begin{equation}\label{P001}
P^{(0)}(\lambda) = \hat{P}^{(0)}(\lambda)\lambda^{-\frac{{a}}{4}\sigma_3}
C\lambda^{\frac{{a}}{4}\sigma_3}e^{-\frac{1}{2}h_0(\lambda)\sigma_3}.
\end{equation}
In the last step of our evaluation of the asymptotics of the function $Z^{{a}}$ we
will need to know exactly the matrix $\hat{P}^{(0)}(0)$. The explicit  formula for
this object is presented in the following proposition.

\begin{prop}\label{P0pop}
The matrix factor $\hat{P}^{(0)}(0)$ in the right hand side of (\ref{P001}) is given
by the equations,
\begin{equation}\label{hatP0for}
\hat{P}^{(0)}(0) = \Delta^{\frac{1}{2}\sigma_3}B\Delta^{-\frac{{a}}{2}\sigma_3}
= \begin{pmatrix}-2^{{a}}\frac{\sqrt{\pi}}{\eta{a}\Gamma\left(-\frac{{a}}{2}\right)}\Delta^{\frac{1}{2} -\frac{{a}}{2}}&
-2^{-{a}-2}\frac{i\eta}{\sqrt{\pi}}\Gamma\left(-\frac{{a}}{2}\right)\Delta^{\frac{1}{2} +\frac{{a}}{2}}\cr\cr
2^{{a}}\frac{\sqrt{\pi}}{\eta\Gamma\left(-\frac{{a}}{2}\right)}\Delta^{-\frac{1}{2} -\frac{{a}}{2}}&
-2^{-{a}-2}\frac{i\eta{a}}{\sqrt{\pi}}\Gamma\left(-\frac{{a}}{2}\right)\Delta^{-\frac{1}{2} +\frac{{a}}{2}}
\end{pmatrix},
\end{equation}
where
\begin{equation}\label{Deltadef}
\Delta = 2(m-in)
\end{equation}
\end{prop}
The proof of the proposition  needs some extra
work with the Bessel functions, and it is moved to Appendix B.

\subsection{Parametrix at $\lambda =\infty$.}

The construction of the parametrix at $\lambda = \infty$ can be done in a complete analogy with the
construction of the parametrix at $\lambda = 0$. However, we can considerably reduce the calculations 
by using the symmetry of the problem with respect to the map $\lambda \mapsto 1/\bar{\lambda}$.

Let $U_{1/\delta}$ be the image of 
$U_{\delta}$ under the map $\lambda \mapsto 1/\bar{\lambda}$. We assume that the pieces of the contours $\Sigma_k$ inside of the neighborhood
$U_{1/\delta}$ are the images, under the map $\lambda \mapsto 1/\bar{\lambda}$,  of the respective pieces of 
$\Sigma_k$ inside of the neighborhood $U_{\delta}$. We notice, that the map preserves  the orientations
of the contours: the ``+'' - side of $\Sigma_{k}\cap U_{\delta}$ goes to the ``+'' side 
of $\Sigma_{k}\cap U_{1/\delta}$ and  the ``-'' - side of $\Sigma_{k}\cap U_{\delta}$ goes to the ``-'' side 
of $\Sigma_{k}\cap U_{1/\delta}$. Secondly, we observe that
\begin{equation}\label{h0hinfty}
\overline{h_0\left(\frac{1}{\bar{\lambda}}\right)}= h_{\infty}(\lambda) + i\pi(m+n).
\end{equation}
We also notice that $\arg{1/\bar{\lambda}} = \arg{\lambda}$ and hence  the branches of all the 
power functions are preserved, and, in particular, 
$$
\overline{\sqrt{\frac{1}{\bar{\lambda}}}} = \frac{1}{\sqrt{\lambda}}
$$

\begin{figure}[h]
  \begin{center}
     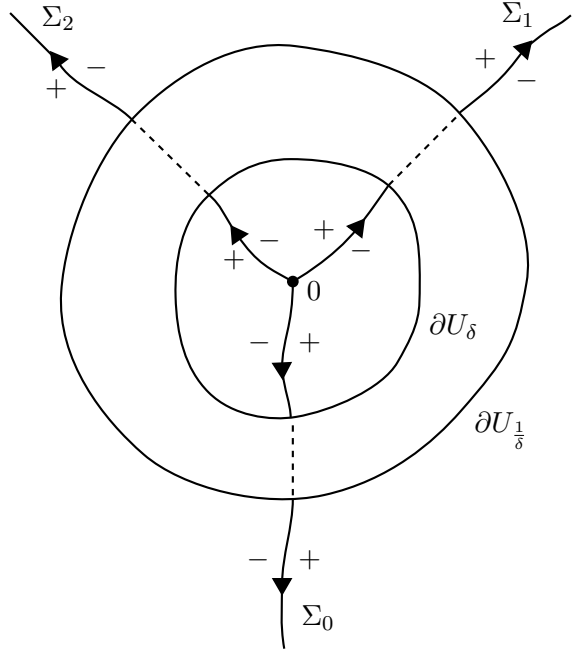
  \end{center}
  \caption{The local $P^{(0)}$ - and $P^{(\infty)}$ - RH problems}\label{f.local_P}
\end{figure}

Consider now again the zero parametrix $P^{(0)}(\lambda)$. By construction, it  solves the following local RH
problem in the neighborhood  $U_{\delta}$. 
\begin{itemize}
\item $P^{(0)}(\lambda)$ is analytic in $U_{\delta}\setminus \Bigl(\Sigma_{k}\cap U_{\delta}\Bigr)$
\item The jump condition is described by the equation,
\begin{equation}\label{P0jump}
P^{(0)}_{+}(\lambda) = P^{(0)}_{-}(\lambda)e^{\frac{1}{2}h_0(\lambda)\sigma_3}\lambda^{-\frac{{a}}{4}\sigma_3}
L_{k}\lambda^{\frac{{a}}{4}\sigma_3}e^{-\frac{1}{2}h_0(\lambda)\sigma_3},\quad
\lambda \in \Sigma_{k}\cap U_{\delta}
\end{equation}
\item as $\lambda \to 0$,
\begin{equation}\label{P00}
P^{(0)}(\lambda) = \hat{P}^{(0)}(\lambda)\lambda^{-\frac{{a}}{4}\sigma_3}
C\lambda^{\frac{{a}}{4}\sigma_3}e^{-\frac{1}{2}h_0(\lambda)\sigma_3},
\end{equation}
where the matrix-valued function $\hat{P}^{(0)}(\lambda)$ is holomorphic at $\lambda = 0$.
\item on the boundary of $U_{\delta}$, the following matching relation with the global parametrix
(\ref{Sgl}) takes place,
\begin{equation}\label{match2}
P^{(0)}(\lambda) = \left(I +  O\left(\frac{1}{r}\right)\right)P^{(gl)}(\lambda),
\quad \lambda \in \partial U_{\delta}, \quad n^2+m^2 \to \infty,
\end{equation}
which , in fact, can be specified as it is indicated in (\ref{match1}).
\end{itemize}
The problem is depicted in Figure \ref{f.local_P}.

Let us indicate explicitly the dependence of the parametrix $P^{(0)}(\lambda)$  and  the matrices $L$ and $C$ 
on the parameter ${a}$, i.e.,
we put,
$$
P^{(0)}(\lambda) \equiv P^{(0)}(\lambda; {a}), \quad L \equiv L({a}), \quad C\equiv C({a}).
$$
We argue, that the parametrix at $\lambda = \infty$ can be defined by the equation,
\begin{equation}\label{Pinftydef}
P^{(\infty)}(\lambda) = \sigma_1\overline{P^{(0)}\left(\frac{1}{\bar{\lambda}}; -{a}\right)}
\end{equation}
We have to check that so defined matrix-valued function solves the following
local RH problem in the neighborhood of infinity, $U_{1/\delta}$.
\begin{itemize}
\item $P^{(\infty)}(\lambda)$ is analytic in $U_{1/\delta}\setminus \Bigl(\Sigma_{k}\cap U_{1/\delta}\Bigr)$
\item The jump condition is described by the equation,
\begin{equation}\label{Pinftyjump}
P^{(\infty)}_{+}(\lambda) = P^{(\infty)}_{-}(\lambda)e^{\frac{1}{2}h_{\infty}(\lambda)\sigma_3}\lambda^{-\frac{{a}}{4}\sigma_3}
L_{k}\lambda^{\frac{{a}}{4}\sigma_3}e^{-\frac{1}{2}h_{\infty}(\lambda)\sigma_3},\quad
\lambda \in \Sigma_{k}\cap U_{1/\delta}
\end{equation}
\item as $\lambda \to \infty$,
\begin{equation}\label{Pinfty0}
P^{(\infty)}(\lambda) = \hat{P}^{(\infty)}(\lambda)\lambda^{-\frac{{a}}{4}\sigma_3}
C\lambda^{\frac{{a}}{4}\sigma_3}e^{-\frac{1}{2}h_{\infty}(\lambda)\sigma_3},
\end{equation}
where the matrix-valued function $\hat{P}^{(\infty)}(\lambda)$ is holomorphic at $\lambda = \infty$.
\item on the boundary of $U_{1/\delta}$, the following matching relation with the global parametrix
(\ref{Sgl}) takes place,
\begin{equation}\label{match3}
P^{(\infty)}(\lambda) = \left(I + O\left(\frac{1}{r}\right)\right)P^{(gl)}(\lambda),
\quad \lambda \in \partial U_{1/\delta}, \quad n^2+m^2 \to \infty.
\end{equation}
\end{itemize}
The problem is depicted in the same Figure 7.

The first condition is trivial; indeed, we have already indicated that under the map $ \lambda \mapsto 1/\bar{\lambda}$
the segments $\Sigma_{k}\cap U_{1/\delta}$ become the segments $\Sigma_{k}\cap U_{\delta}$ with the preservation
of the respective sides of the segments. In order to check the jump relations (\ref{Pinftyjump}), we should use (\ref{h0hinfty}) and the 
obvious equation, 
\begin{equation}\label{L-L}
\overline{L_k(-{a})} = L_k({a}) \equiv L_k.
\end{equation}
We would have that (taking into account that $m+n$ is even),
$$
P^{(\infty)}_{+}(\lambda) = \sigma_1\overline{P^{(0)}_{+}\left(\frac{1}{\bar{\lambda}}; -{a}\right)} 
$$
$$
=\sigma_1\overline{P^{(0)}_{-}\left(\frac{1}{\bar{\lambda}}; -{a}\right)} 
e^{\frac{1}{2}h_{\infty}(\lambda)\sigma_3 + \frac{1}{2}i\pi(m+n)\sigma_3}\lambda^{-\frac{{a}}{4}\sigma_3}
{\overline{L_{k}(-{a})}}\lambda^{\frac{{a}}{4}\sigma_3}e^{-\frac{1}{2}h_{\infty}(\lambda)\sigma_3 - \frac{1}{2}i\pi(m+n)\sigma_3}
$$
$$
= P^{(\infty)}_{-}(\lambda)e^{\frac{1}{2}h_{\infty}(\lambda)\sigma_3}\lambda^{-\frac{{a}}{4}\sigma_3}
L_{k}\lambda^{\frac{{a}}{4}\sigma_3}e^{-\frac{1}{2}h_{\infty}(\lambda)\sigma_3},\quad
\lambda \in \Sigma_{k}\cap U_{1/\delta}.
$$
Since the matrix $C$ satisfies the same relation (\ref{L-L}) as the matrices  $L_k$, we would have 
condition (\ref{Pinfty0}) at $\lambda = \infty$ with
\begin{equation}\label{Phatinfty}
\hat{P}^{(\infty)}(\lambda) = \sigma_1\overline{\hat{P}^{(0)}\left(\frac{1}{\bar{\lambda}}; -{a}\right)}.
\end{equation}
Finally, we observe that
$$
\sigma_1\overline{P^{(gl)}\left(\frac{1}{\bar{\lambda}}; -{a}\right)}=
\sigma_1\lambda^{\frac{\sigma_3}{4}}\begin{pmatrix}\frac{1}{2}&\frac{i}{2}\cr\cr
\frac{1}{2}&-\frac{i}{2}\end{pmatrix}\eta^{-\sigma_3}\lambda^{\frac{{a}}{4}\sigma_3}
$$
$$
= \sigma_1\lambda^{\frac{\sigma_3}{4}}\sigma_1\sigma_1\begin{pmatrix}\frac{1}{2}&\frac{i}{2}\cr\cr
\frac{1}{2}&-\frac{i}{2}\end{pmatrix}\eta^{-\sigma_3}\lambda^{\frac{{a}}{4}\sigma_3}
=\lambda^{-\frac{\sigma_3}{4}}\begin{pmatrix}\frac{1}{2}&-\frac{i}{2}\cr\cr
\frac{1}{2}&\frac{i}{2}\end{pmatrix}\eta^{-\sigma_3}\lambda^{\frac{{a}}{4}\sigma_3} = P^{(gl)}(\lambda)
$$
Therefore, on the boundary of $U_{1/\delta}$, we have,
$$
P^{(\infty)}(\lambda) = \left( I + O\left(\frac{1}{r}\right)\right)\sigma_1\overline{P^{(gl)}\left(\frac{1}{\bar{\lambda}}; -{a}\right)}
= \left( I + O\left(\frac{1}{r}\right)\right)P^{(gl)}(\lambda), \quad n^2+m^2 \to \infty.
$$
That is, the matching condition (\ref{match3}) is satisfied. This completes the proof that equation (\ref{Pinftydef}) indeed
defines a parametrix for the $S$ - RH problem in the neighborhood of $\lambda = \infty$. It should be also noticed
that from (\ref{match1}) the similar specification  of (\ref{match3}) follows,
\begin{equation}\label{match1inf}
P^{(\infty)}(\lambda) = \left(I + \frac{1}{\sqrt{\xi\left(\frac{1}{\lambda}\right)}}\lambda^{-\frac{1}{4}\sigma_3}
\sigma_1\overline{\Psi_1}(-{a})\sigma_1\lambda^{\frac{1}{4}\sigma_3} +  O\left(\frac{1}{r^2}\right)\right)P^{(gl)}(\lambda),
\quad \lambda \in \partial U_{1/\delta}, \quad n^2+m^2 \to \infty.
\end{equation}

In addition, formula (\ref{Phatinfty}), together with (\ref{hatP0for}) implies the following expression for the matrix
$\hat{P}^{(\infty)}(\infty)$.
\begin{equation}\label{hatPinffor}
\hat{P}^{(\infty)}(\infty) = \sigma_1\overline{\hat{P}^{(0)}(0; -{a})}
=
\begin{pmatrix}2^{-{a}}\frac{\sqrt{\pi}}{\eta\Gamma\left(\frac{{a}}{2}\right)}\bar{\Delta}^{-\frac{1}{2} +\frac{{a}}{2}}&
-2^{{a}-2}\frac{i\eta{a}}{\sqrt{\pi}}\Gamma\left(\frac{{a}}{2}\right)\bar{\Delta}^{-\frac{1}{2} -\frac{{a}}{2}}\cr\cr
2^{-{a}}\frac{\sqrt{\pi}}{\eta{a}\Gamma\left(\frac{{a}}{2}\right)}\bar{\Delta}^{\frac{1}{2} +\frac{{a}}{2}}&
2^{{a}-2}\frac{i\eta}{\sqrt{\pi}}\Gamma\left(\frac{{a}}{2}\right)\bar{\Delta}^{\frac{1}{2} -\frac{{a}}{2}}
\end{pmatrix}
\end{equation}

\subsection{Asymptotic solution of the $S$ - RH problem}\label{sec3.6}
It is convenient to pass from the matrix-valued function $S$ to
the function (cf. (\ref{Shatinfty})),
\begin{equation}\label{Stildedef}
\tilde{S}(\lambda) = \left[\hat{S}^{(\infty)}(\infty)\right]^{-1}S(\lambda).
\end{equation}
The function $\tilde{S}(\lambda)$ satisfies the same  $S$ - RH problem 
except that the condition at infinity (\ref{Sinfty0}) is replaced by the
more standard condition,
\begin{itemize}
\item as $\lambda \to \infty$,
\begin{equation}\label{Stildeinfty0}
\tilde{S}(\lambda) = \left(I +O\left(\frac{1}{\lambda}\right)\right)
 \lambda^{-\frac{{a}}{4}\sigma_3}C\lambda^{\frac{{a}}{4}\sigma_3}e^{-\frac{1}{2}h_{\infty}(\lambda)\sigma_3},
\end{equation}
\end{itemize}
The solution  $S(\lambda)$ of the $S$ - RH problem can be recovered from the solution $\tilde{S}(\lambda)$ of
the $\tilde{S}$ - RH  problem via the equation
\begin{equation}\label{StildeS}
S(\lambda) = M\tilde{S}(\lambda),
\end{equation}
where the matrix $M$ is uniquely determined by the properties,
\begin{equation}\label{M1}
M =  \begin{pmatrix}\bullet&0\cr
 \bullet&\bullet\end{pmatrix},\quad MD = \begin{pmatrix}1&\bullet\cr
 0&1\end{pmatrix}\begin{pmatrix}1&0\cr
 0&(-1)^{m}\end{pmatrix}e^{-\frac{i\pi}{2}n\sigma_3},
 \end{equation}
 where the matrix $D$ is the left constant matrix factor in the representation of
 the solution $\tilde{S}(\lambda)$ at $\lambda =0$ (cf. (\ref{Szero}),
 \begin{equation}\label{DStilde}
\tilde{S}(\lambda) =  D\Bigl(I + O(\lambda)\Bigr)
 \lambda^{-\frac{{a}}{4}\sigma_3}C\lambda^{\frac{{a}}{4}\sigma_3}e^{-\frac{1}{2}h_0(\lambda)\sigma_3}, 
\end{equation}

>From (\ref{M1}) it follows that 
 \begin{equation}\label{M2}
 M_{12} = 0, \quad M_{11} = \frac{1}{D_{11}}e^{-\frac{i\pi}{2}n},\quad M_{21}=(-1)^{m+1}e^{\frac{i\pi}{2}n}D_{21},
 \quad M_{22}=(-1)^{m}e^{\frac{i\pi}{2}n}D_{11}
 \end{equation}
 This in turn means that
$$
\hat{S}^{(0)}_{12} = M_{11}D_{12} = \frac{D_{12}}{D_{11}}e^{-\frac{i\pi}{2}n}.
$$
The last equation allows us to rewrite (\ref{ZgammaS}) in term of the $\tilde{S}$ - function,
 \begin{equation}\label{ZgammatildeS}
 Z^{{a}} = (-1)^{m+1}e^{-i\pi n}\frac{D_{12}}{D_{11}} = -\frac{D_{12}}{D_{11}},
 \end{equation}
 where we again took into account that $m+n$ is even.
We shall now present  the asymptotic solution of the  $\tilde{S}$ - RH problem.

\begin{figure}[h]
  \begin{center}
     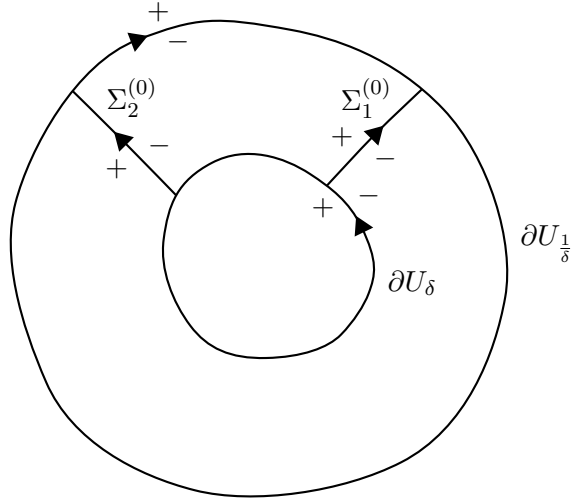
  \end{center}
  \caption{The contour for the $R$ - RH problem}\label{f.contour_R}
\end{figure}

Define the piecewise analytic  function,
\begin{equation}\label{Sasdef}
S^{(as)}(\lambda) = 
%\left[\hat{P}^{(\infty)}(\infty)\right]^{-1}
 \begin{cases}
P^{(0)}(\lambda) &\lambda \in U_{\delta},\cr\cr
P^{(\infty)}(\lambda) &\lambda \in U_{1/\delta},\cr\cr
P^{(gl)}(\lambda)  & \lambda \in {\bf C} \setminus (U_{\delta}\cup U_{1/\delta}),
\end{cases}
\end{equation}
and consider the matrix ratio,
\begin{equation}\label{Rdef}
R(\lambda) =\hat{P}^{(\infty)}(\infty) \tilde{S}(\lambda)[S^{(as)}(\lambda)]^{-1}.
\end{equation}
The function $R(\lambda)$ is the piece-wise analytic matrix-valued function whose jump-contour is
\begin{equation}\label{Rcontour}
\Sigma_R = \partial U_{\delta} \cup \partial U_{1/\delta}\cup \Sigma^{(0)}_1 \cup   \Sigma^{(0)}_2, 
\end{equation}
where $\Sigma^{(0)}_1$ and  $\Sigma^{(0)}_2$ denote the segments of the rays $\Sigma_1$ and  $\Sigma_2$,
respectively, included between the curves  $\partial U_{\delta}$ and  $\partial U_{1/\delta}$. It should be
noted that, since the functions $\tilde{S}(\lambda)$ and  the function $S^{(as)}(\lambda)$ share the same jump 
matrices on the ray $\Sigma_0$ and on the parts of the rays  $\Sigma_1$ and  $\Sigma_2$ which are
inside  the neighborhoods $U_{\delta}$ and $U_{1/\delta}$, the function $R(\lambda)$ is continuous across 
these  pieces of the contour $\Gamma$. On the contour $\Sigma_R$, the function $R(\lambda)$ solves the following Riemann-Hilbert problem.
\begin{itemize}
\item $R(\lambda)$  is analytic on $\Bbb{C}\setminus\Sigma_R$.
\item  The jump conditions are described by the equations,
\begin{enumerate}
\item as $\lambda \in \Sigma^{(0)}_{1}$,
\begin{equation}\label{SSigma110}
R_+(\lambda) =R_-(\lambda)P^{gl}(\lambda)
 \begin{pmatrix}1 &0\cr
                                     H^{-1}(\lambda)\omega^{-1}_{1}(\lambda) & 1\end{pmatrix}\left[P^{gl}(\lambda)\right]^{-1}, 
\end{equation}
\item as $\lambda \in \Sigma^{(0)}_{2}$,
\begin{equation}\label{SSigma220}
R_+(\lambda) =R_-(\lambda)P^{gl}(\lambda)
 \begin{pmatrix}1 &0\cr
                                     H(\lambda)\omega^{-1}_{2}(\lambda) & 1\end{pmatrix}\left[P^{gl}(\lambda)\right]^{-1}, 
\end{equation}
\item as $\lambda \in \partial U_{\delta}$,
\begin{equation}\label{Udjump}
R_+(\lambda) =R_-(\lambda)P^{gl}(\lambda)\left[P^{(0)}(\lambda)\right]^{-1}
\end{equation}
\item as $\lambda \in \partial U_{1/\delta}$,
\begin{equation}\label{Udjump2}
R_+(\lambda) =R_-(\lambda))P^{gl}(\lambda)\left[P^{(\infty)}(\lambda)\right]^{-1}
\end{equation}
\end{enumerate}
\item The function $R(\lambda)$ is normalized by the conduition,
\begin{equation}\label{inftyR}
R(\infty) = I
\end{equation}
\end{itemize} 
It is also worth noticing that at the node points of the graph $\Sigma_R$ the function
$R(\lambda)$ is bounded and its monodromy at each node point is trivial.  The $R$ - RH problem
is depicted in Figure \ref{f.contour_R}.

Let $G_{R}(\lambda)$ denote the $R$-jump matrix. Then, in view of Lemma \ref{lemma1},
we have that there exists a positive constant $c_0$ such that 
\begin{equation}\label{GRest1}
G_R(\lambda) = I + O\left(e^{-c_0r}\right), 
\end{equation}
for all $\lambda \in \Sigma^{(0)}_1\cup\Sigma^{(0)}_2$, as $n, m  \to \infty$. Simultaneously,
the estimates (\ref{match2}) and (\ref{match3}) imply that
\begin{equation}\label{GRest2}
G_R(\lambda) = I + O\left(\frac{1}{r}\right), 
\end{equation}
for all $\lambda \in \partial U_{\delta}\cup\partial U_{1/\delta}$, as $n, m  \to \infty$. Taking
into account (\ref{match1}) and (\ref{match1inf}), we can specify  estimate (\ref{GRest2})
as 
\begin{equation}\label{G1}
G_R(\lambda) = I  + G^{(0)}_1(\lambda) +  O\left(\frac{1}{r^2}\right),
\quad G^{(0)}_1(\lambda) \equiv \frac{1}{\sqrt{\xi(\lambda)}}\lambda^{-\frac{1}{4}\sigma_3}
\Psi_1\lambda^{\frac{1}{4}\sigma_3}
\end{equation}
if $\lambda \in \partial U_{\delta}$ and 
\begin{equation}\label{G2}
G_R(\lambda) = I  +  G^{(\infty)}_1(\lambda) +  O\left(\frac{1}{r^2}\right),
\quad G^{(\infty)}_1(\lambda)\equiv \frac{1}{\sqrt{\xi\left(\frac{1}{\lambda}\right)}}\lambda^{-\frac{1}{4}\sigma_3}
\sigma_1\overline{\Psi_1}(-{a})\sigma_1\lambda^{\frac{1}{4}\sigma_3}
\end{equation}
if $\lambda \in \partial U_{1/\delta}$. We note that, because of the off-diagonal
structure of the matrix $\Psi_1$ (cf. (\ref{Psi1})), the matrix functions 
$G^{(0)}_1(\lambda)$ and $G^{(\infty)}_1(\lambda)$ are holomorphic in
$U_{\delta}\setminus\{0\}$ and $U_{1/\delta}\setminus\{\infty\}$, respectively.

In their turn,
asymptotic relations (\ref{GRest1}) and (\ref{GRest2}) yield the estimate,
\begin{equation}\label{Grest2}
||I - G_R||_{L_1(\Sigma_R)\cap L_2(\Sigma_R)\cap L_{\infty}(\Sigma_R)} \leq \frac{L}{r}, \quad r>1,
\end{equation}
with some positive constant $L$. The standard arguments \cite{DZ} (see also Theorem 1.5 in \cite{I})
lead then to the asymptotic relation,
\begin{equation}\label{Ras}
R(\lambda) = I + O\left(\frac{1}{(1 + |\lambda|)r}\right), \quad n^2+m^2 \to \infty,
\end{equation}
 uniformly on every closed subset of ${\Bbb C}P^1$ outside of the contour 
$\Sigma_R$. Hence we arrive at the following asymptotic representation of the
solution  of the $\tilde{S}$ - problem.
\begin{theorem}\label{asStilde}. Let $\tilde{S}(\lambda)$ be the solution of the $\tilde{S}$ - problem.
Then,
\begin{equation}\label{tildeSas}
\tilde{S}(\lambda) = \left[\hat{P}^{(\infty)}(\infty)\right]^{-1}
\left(I + O\left(\frac{1}{(1 + |\lambda|)r}\right)\right)S^{(as)}(\lambda), 
\quad n^2+m^2 \to \infty,
\end{equation}
 uniformly on every closed subset of ${\Bbb C}P^1$ outside of the contour 
$\Sigma_R$. 
\end{theorem}

\subsection {Asymptotics of $Z^{{a}}$. The completion of proof of theorem \ref{theorem1}  for  the case of even $n+m$ }
The matrix factor $D$ from (\ref{DStilde}) is
given by the equation,
\begin{equation}\label{DStilde0}
D = \left[\hat{P}^{(\infty)}(\infty)\right]^{-1}R(0)\hat{P}^{(0)}(0).
\end{equation}
This, together with (\ref{Ras}) yields at once the   asymptotic equation,
\begin{equation}\label{DStildeas}
D = \left[\hat{P}^{(\infty)}(\infty)\right]^{-1}
\left(I + O\left(\frac{1}{r}\right)\right)\hat{P}^{(0)}(0), \quad n^2+m^2 \to \infty.
\end{equation}
We will need, however, a more detail information about the structure of the
estimate (\ref{DStildeas} ).
\begin{prop}\label{R0str}.The matrix entries of $R(0)$ satisfy the estimates,
$$
R_{11}(0) = 1 + O\left(\frac{1}{r^2}\right), \quad R_{22}(0)  = 1 + O\left(\frac{1}{r^2}\right),
\quad R_{12}(0) = O\left(\frac{1}{r}\right), \quad R_{21}(0)  = O\left(\frac{1}{r}\right).
$$
\end{prop}
{\it Proof.}
The matrix function $R(\lambda)$ admits the following integral representation (see again
\cite{DZ}, \cite{I}),
\begin{equation}\label{Requ1}
R(\lambda) = I + \frac{1}{2\pi i}\int_{\Sigma_{R}}\frac{\rho(\mu)(I - G_R(\mu))}{\mu - \lambda} d\mu,
\quad \lambda \in \Bbb{C}\setminus\Sigma_R,
\end{equation}
where  the matrix function $\rho(\lambda) \equiv R_-(\lambda)$ solves  the singular
integral equation,
\begin{equation}\label{Requ2}
\rho(\lambda) = I + \frac{1}{2\pi i}\int_{\Sigma_{R}}\frac{\rho(\mu)(I - G_R(\mu))}{\mu - \lambda_-} d\mu,
\quad \lambda \in \Sigma_R.
\end{equation}
In (\ref{Requ2}),  the singular Cauchy operator in the right hand side is defined by the formula,
$$
 \frac{1}{2\pi i}\int_{\Sigma_{R}}\frac{\rho(\mu)(I - G_R(\mu))}{\mu - \lambda_-} d\mu := 
\lim_{\lambda' \to \lambda,\, \, \lambda' \in - \,\mbox{side of }\,\, \Sigma_R} 
 \int_{\Sigma_{R}}\frac{\rho(\mu)(I - G_R(\mu))}{\mu - \lambda'} d\mu.
$$
Equation (\ref{Requ2}) is considered as an equation in $L_2(\Sigma_{R})$. 
>From the general theory (see again \cite{DZ}),  it follows that
estimate (\ref{Grest2}) implies the large $r$  solvability of equation (\ref{Requ2}) (which, in fact, we have a priori for all $r>0$ ) and
the estimate
\begin{equation}\label{Requ3}
||I - \rho||_{L_2(\Sigma_R)}\leq \frac{L}{r}, \quad r>1.
\end{equation}
Applying this estimate to (\ref{Requ1}), we come to the conclusion that
\begin{equation}\label{Requ4}
R(0) = I + \frac{1}{2\pi i}\int_{\Sigma_{R}}\frac{(I - G_R(\lambda))}{\lambda} d\lambda + O\left(\frac{1}{r^2}\right),
\end{equation}
as $n^2+ m^2 \to \infty$. Taking into account (\ref{GRest1}) and (\ref{G1}), (\ref{G2}) we see that one can
replace in (\ref{Requ4}) the contour of integration by the union $\partial U_{\delta} \cup \partial U_{1/\delta}$,
and the difference $I - G_R(\lambda)$ by $-G^{(0)}_1(\lambda)$ and $-G^{(\infty)}_1(\lambda)$. In other words, we have that
\begin{equation}\label{Requ5}
R(0) = I  -\mbox{res}_{\lambda = 0}\frac{1}{\lambda}G^{(0)}_1(\lambda)
- \mbox{res}_{\lambda = \infty}\frac{1}{\lambda}G^{(\infty)}_1(\lambda) + O\left(\frac{1}{r^2}\right),
\end{equation}
and remembering  the off-diagonal structure of the matrices $G^{(0)}_1(\lambda)$ and $G^{(\infty)}_1(\lambda)$,
the proposition follows. 

Denote
$$
p_{jk} = (\hat{P}^{(0)}(0))_{jk}, \quad q_{jk} = (\hat{P}^{(\infty)}(\infty))_{jk}.
$$
Then, taking into account that $\det\hat{P}^{(\infty)}(\infty) =  \det\hat{P}^{(0)}(0) = \det B = i/2$ (see
(\ref{hatP0for}) and (\ref{B0}), we would
have from (\ref{DStilde0}) that
\begin{equation}\label{D11}
D_{11} = -2iq_{22}p_{11}\left(R_{11}(0) -\frac{q_{12}}{q_{22}}R_{21}(0)
+\frac{p_{21}}{p_{11}}R_{12}(0) - \frac{q_{12}p_{21}}{q_{22}p_{11}}R_{22}(0)\right),
\end{equation}
and
\begin{equation}\label{D12}
D_{12} = -2iq_{22}p_{12}\left(R_{11}(0) -\frac{q_{12}}{q_{22}}R_{21}(0)
+\frac{p_{22}}{p_{12}}R_{12}(0) - \frac{q_{12}p_{22}}{q_{22}p_{12}}R_{22}(0)\right),
\end{equation}
Using Proposition \ref{R0str} and 
recalling explicit formulae for the matrices $\hat{P}^{(\infty)}(\infty)$ and $\hat{P}^{(0)}(0)$,  i.e., formulae (\ref{hatP0for})
and (\ref{hatPinffor}), respectively, we derive from the equations (\ref{D11}) and (\ref{D12})
the following estimates for the matrix entries $D_{11}$ and $D_{12}$ ,
\begin{equation}\label{D111}
D_{11} = -2iq_{22}p_{11}\left(1 +O\left(\frac{1}{r^2}\right)\right) ,
\end{equation}
and 
\begin{equation}\label{D121}
D_{12} = -2iq_{22}p_{12}\left(1 +O\left(\frac{1}{r^2}\right)\right).
\end{equation}
Substituting (\ref{D111}) and (\ref{D121}) into (\ref{ZgammatildeS}) we obtain that,
$$
Z^{{a}} = -\frac{p_{12}}{p_{11}}\left(1 +O\left(\frac{1}{r^2}\right)\right),
$$
or, looking one more time at (\ref{hatP0for}), 
\begin{equation}\label{Zgas1}
Z^{{a}} = -i2^{-2{a} -2}\frac{\eta^2 {a}}{\pi}\Gamma^{2}\left(-\frac{{a}}{2}\right)
\Delta^{{a}}\left(1 + O\left(\frac{1}{r^2}\right)\right), \quad n^2+m^2 \to \infty.
\end{equation}
Taking into account the definition (\ref{thetadef}) of the branch of the argument of 
$m-in$ and the assumption that $0< \arg (m+in) < \pi/2$, we see that
$$
\Delta^{{a}} = 2^{{a}}e^{-\frac{i\pi{a}}{2}}(n+im)^{{a}},
$$
and therefore,
$$
\eta^2\Delta^{{a}} = 2^{{a} +1}i\sin\frac{\pi{a}}{2}(n+im)^{{a}}.
$$
The last equation allows to rewrite (\ref{Zgas1}) as
\begin{equation}\label{Zgas2}
Z^{{a}} = \left(\frac{n+im}{2}\right)^{{a}}\frac{1}{\pi}\sin\frac{\pi{a}}{2}\Gamma^{2}\left(-\frac{{a}}{2}\right)\frac{{a}}{2}
\left(1 + O\left(\frac{1}{r^2}\right)\right), \quad n^2+m^2 \to \infty,
\end{equation}
Since,
$$
\frac{1}{\pi}\sin\frac{\pi{a}}{2}\Gamma^{2}\left(-\frac{{a}}{2}\right)\frac{{a}}{2}
=-\frac{\Gamma\left(-\frac{{a}}{2}\right)}{\Gamma\left(1+\frac{{a}}{2}\right)}\frac{{a}}{2}
= \frac{\Gamma\left(1-\frac{{a}}{2}\right)}{\Gamma\left(1+\frac{{a}}{2}\right)},
$$
equation (\ref{Zgas2}) is equivalent (\ref{ABconj}) and hence Theorem\ref{theorem1} is proven
for the case of the even sum $n+m$.

\subsection{Extension to the general case.}\label{oddcase}
We need to extend the validity of asymptotic formula (\ref{ABconj}) to the case of 
the odd value of the sum $n+m$. It is obvious that this will be achieved if we, still
assuming the evenness of the sum $n+m$, will be able to extract from the considerations of the
previous sections not only the asymptotics of $f_{n,m}$ but the asymptotics
of the quantities $f_{n+1,m}$ or $f_{n,m+1}$ as well. In order to have that, in virture of equations
(\ref{uv}), it is enough to find the asymptotic behavior
of the discrete derivatives $u_{n,m}$ and $v_{n,m}$.

We start with noticing that from (\ref{B}) and (\ref{C}) it follows that
\begin{equation}\label{unm}
u_{n,m} = \frac{(B_{n,m})_{11}}{(B_{n,m})_{21}},
\end{equation}
and
\begin{equation}\label{vnm}
v_{n,m} = \frac{(C_{n,m})_{11}}{(C_{n,m})_{21}}.
\end{equation}
Matrices $B_{n,m}$ and $C_{n,m}$, in their turn,  can be determined
through the left holomorphic factors in the representations
(\ref{-10}) and (\ref{10}) of the  function $\Psi_{n,m}(\lambda)$
near the points $-1$ and $1$, respectively. Indeed we have,
\begin{equation}\label{B1}
B_{n,m} = -n\hat{\Psi}^{(-1)}_{n,m}(-1)\begin{pmatrix}0&0\cr\cr
0&1\end{pmatrix}\left[\hat{\Psi}^{(-1)}_{n,m}(-1)\right]^{-1},
\end{equation}
and
\begin{equation}\label{C1}
C_{n,m} = -m\hat{\Psi}^{(1)}_{n,m}(1)\begin{pmatrix}0&0\cr\cr
0&1\end{pmatrix}\left[\hat{\Psi}^{(1)}_{n,m}(1)\right]^{-1},
\end{equation}
If we trace all the transformations which we made when moving from the original monodromy
problem (\ref{1}) - (\ref{0}) to the final $S$ - problem (\ref{SSigma00}) - 
(\ref{Szero0}), we will easily find out that
$$
\hat{\Psi}^{(-1)}_{n,m}(\lambda) =
S(\lambda)\begin{pmatrix}1&0\cr\cr
-H(\lambda)\omega^{-1}_{2}(\lambda)&1\end{pmatrix}
\begin{pmatrix}1&0\cr\cr
0&(\lambda -1)^{m}\end{pmatrix},
$$
and
$$
\hat{\Psi}^{(1)}_{n,m}(\lambda) =
S(\lambda)\begin{pmatrix}1&0\cr\cr
-H^{-1}(\lambda)\omega^{-1}_{1}(\lambda)&1\end{pmatrix}
\begin{pmatrix}1&0\cr\cr
0&(\lambda +1)^{n}\end{pmatrix}.
$$
Taking into account that $H(-1) =0$ and $H^{-1}(1) = 0$, we see that
$$
\hat{\Psi}^{(-1)}_{n,m}(-1) =
S(-1)
\begin{pmatrix}1&0\cr\cr
0&(-2)^{m}\end{pmatrix},
$$
and
$$
\hat{\Psi}^{(1)}_{n,m}(1) =
S(1)
\begin{pmatrix}1&0\cr\cr
0&2^{n}\end{pmatrix},
$$
and hence equations (\ref{B1}) and (\ref{C1}) can
be  rewritten directly in terms of the function $S(\lambda)$,
\begin{equation}\label{B2}
B_{n,m} = -nS(-1)\begin{pmatrix}0&0\cr\cr
0&1\end{pmatrix}S^{-1}(-1),
\end{equation}
and
\begin{equation}\label{C2}
C_{n,m} = -mS(1)\begin{pmatrix}0&0\cr\cr
0&1\end{pmatrix}S^{-1}(1).
\end{equation}
(We remind that we always suppress the indication of the dependence of
$S(\lambda)$ on $n$ and $m$.) 
As a consequence, the basic relations (\ref{unm}) and
(\ref{vnm}) for the discrete functions $u_{n,m}$ and $v_{n,m}$ 
can be replaced by the equations,
\begin{equation}\label{unm1}
u_{n,m} = \frac{(S(-1))_{12}}{(S(-1))_{22}},
\end{equation}
and
\begin{equation}\label{vnm2}
v_{n,m} = \frac{(S(1))_{12}}{(S(1))_{22}},
\end{equation}

>From the formulae (\ref{StildeS}) and (\ref{Rdef}) it follows  that the solution $S(\lambda)$
of the $S$ - RH problem can be written in the form of the product,
$$
S(\lambda) = M\left[\hat{P}^{(\infty)}(\infty)\right]^{-1}R(\lambda)S^{(as)}(\lambda)
$$
Taking into account (\ref{M1}) and (\ref{DStilde0}), the last equation can be transformed 
into the relation,
$$
S(\lambda) = \begin{pmatrix}1&-f_{n,m}\cr\cr
0&1\end{pmatrix}
\begin{pmatrix}1&0\cr\cr
0&(-1)^{m}\end{pmatrix}e^{-\frac{i/pi}{2}n\sigma_{3}}
\left[\hat{P}^{(0)}(0)\right]^{-1}R^{-1}(0)R(\lambda)S^{(as)}(\lambda),
$$
which in turn implies that,
$$
S(\pm 1) = S(\lambda) = \begin{pmatrix}1&-f_{n,m}\cr\cr
0&1\end{pmatrix}
\begin{pmatrix}1&0\cr\cr
0&(-1)^{m}\end{pmatrix}e^{-\frac{i/pi}{2}n\sigma_{3}}
$$
\begin{equation}\label{Spm1}
\times \left[\hat{P}^{(0)}(0)\right]^{-1}R^{(\pm)}P^{(gl)}(\pm 1),
\end{equation}
where
\begin{equation}\label{Rpm}
R^{(\pm)} = R^{-1}(0)R(\pm 1) = I + O\left(\frac{1}{r}\right), \quad r = \sqrt{n^2 +m^2} \to \infty.
\end{equation}
When deriving (\ref{Spm1}), we have used the fact that $\pm 1 \notin {\Bbb C} \setminus (U_{\delta}\cup U_{1/\delta})$
and hence $S^{(as)}(\pm 1) = P^{(gl)}(\pm 1)$ in accord with  definition (\ref{Sasdef}) of the parametrix $S^{(as)}$.
It also should be noticed that the matrix entries of $R^{(\pm)}$ admit the same type of
specification of the estimate (\ref{Rpm}) as in Proposition \ref{R0str}.

Equation (\ref{Spm1}) allows us to estimate the quantities $(S(\pm1))_{12}$ and $(S(\pm1))_{22}$
involved in the formulae (\ref{unm1}) and (\ref{vnm2}). To this end, we first  notice that,
as it follows from equation (\ref{Sgl}) and the convention about the branches of the multivalued functions used
in (\ref{Sgl}) (i.e., $-\pi/2 <\arg\lambda<3\pi/2$), we have that,
\begin{equation}\label{Pglpm1}
P^{(gl)}(\pm 1) =
\begin{pmatrix}\frac{1}{2\eta}e^{\frac{\sigma i\pi}{4}({a}-1)}&
-\frac{i\eta}{2}e^{-\frac{\sigma i\pi}{4}({a}+1)}\cr\cr
\frac{1}{2\eta}e^{\frac{\sigma i\pi}{4}({a}+1)}&
\frac{i\eta}{2}e^{-\frac{\sigma i\pi}{4}({a}-1)}
\end{pmatrix},
\end{equation}
where $\sigma = 0$ in the case $P^{(gl)}( 1)$ and $\sigma = 1$ in the case $P^{(gl)}(-1)$.
Secondly, taking into account that $\det \hat{P}^{(0)}(0) = i/2$ (cf. (\ref{B0})), we can write
\begin{equation}\label{Pglpm2}
\left[\hat{P}^{(0)}(0)\right]^{-1} = -2i
\begin{pmatrix} p_{22}&-p_{12}\cr\cr
-p_{21}&p_{11}\end{pmatrix},\quad p_{jk} = (\hat{P}^{(0)}(0))_{jk}.
\end{equation}
Substituting (\ref{Pglpm1}) and (\ref{Pglpm2}) into (\ref{Spm1}) and skipping
some strightforwrad though tedious calculations, we arrive at the
following representations for $(S(\pm1))_{12}$ and $(S(\pm1))_{22}$.
\begin{equation}\label{S12}
(S(\pm 1))_{12} = {\cal{A}} - f_{n,m}{\cal{B}},\quad (S(\pm 1))_{22} = {\cal{B}},
\end{equation}
where
\begin{equation}\label{calA}
{\cal{A}} = -\eta e^{-\frac{i\pi}{2}n -\frac{\sigma i\pi}{4}(1 +{a})}
\left(p_{22}\left(R^{(\pm)}_{11} -e^{\frac{\sigma i\pi}{2}}R^{(\pm)}_{12}\right)
-p_{12}\left(R^{(\pm)}_{21} -e^{\frac{\sigma i\pi}{2}}R^{(\pm)}_{22}\right)\right),
\end{equation}
and
\begin{equation}\label{calB}
{\cal{B}} = -\eta  (-1)^{m}e^{\frac{i\pi}{2}n -\frac{\sigma i\pi}{4}(1 +{a})}
\left(-p_{21}\left(R^{(\pm)}_{11} -e^{\frac{\sigma i\pi}{2}}R^{(\pm)}_{12}\right)
+p_{11}\left(R^{(\pm)}_{21} -e^{\frac{\sigma i\pi}{2}}R^{(\pm)}_{22}\right)\right).
\end{equation}
Substituting, in turn, these equations into the right hand sides of formulae
(\ref{unm1}) and (\ref{vnm2}), we obtain that (we remind that we are still
assuming that $m+n$ is even),
$$
u_{n,m} = -f_{n,m} 
+ \frac{p_{22}\left(R^{(-)}_{11} -iR^{(-)}_{12}\right)
-p_{12}\left(R^{(-)}_{21} -iR^{(-)}_{22}\right)}
{-p_{21}\left(R^{(-)}_{11} -iR^{(-)}_{12}\right)
+p_{11}\left(R^{(-)}_{21} -iR^{(-)}_{22}\right)},
$$
and
$$
v_{n,m} = -f_{n,m} 
+ \frac{p_{22}\left(R^{(+)}_{11} -R^{(+)}_{12}\right)
-p_{12}\left(R^{(+)}_{21} -R^{(+)}_{22}\right)}
{-p_{21}\left(R^{(+)}_{11} -R^{(+)}_{12}\right)
+p_{11}\left(R^{(+)}_{21} -R^{(+)}_{22}\right)},
$$
respectively. Remembering now equations (\ref{uv}), 
the last equations become in fact the equations for $f_{n+1, m}$
and $f_{n,m+1}$, respectively. That is we have,
\begin{equation}\label{fn1m}
f_{n+1,m}= 
\frac{p_{22}\left(R^{(-)}_{11} -iR^{(-)}_{12}\right)
-p_{12}\left(R^{(-)}_{21} -iR^{(-)}_{22}\right)}
{-p_{21}\left(R^{(-)}_{11} -iR^{(-)}_{12}\right)
+p_{11}\left(R^{(-)}_{21} -iR^{(-)}_{22}\right)},
\end{equation}
and
\begin{equation}\label{fnm1}
f_{n,m+1} = \frac{p_{22}\left(R^{(+)}_{11} -R^{(+)}_{12}\right)
-p_{12}\left(R^{(+)}_{21} -R^{(+)}_{22}\right)}
{-p_{21}\left(R^{(+)}_{11} -R^{(+)}_{12}\right)
+p_{11}\left(R^{(+)}_{21} -R^{(+)}_{22}\right)}.
\end{equation}

We are ready now to produce the asymptotic formulae   for $f_{n+1, m}$
and $f_{n,m+1}$. Indeed, taking from (\ref{hatP0for}) the exact expressions for $p_{jk} $
we derive from (\ref{fn1m}) and (\ref{fnm1}) the  relations,
\begin{equation}\label{fn1m2}
f_{n+1,m} = -i2^{-2{a} -2}\frac{\eta^2{a}}{\pi}\Gamma^2\left(-\frac{{a}}{2}\right)
\Delta^{{a}}\frac{1 -\frac{{a}}{\Delta}\kappa_-}
{1 +\frac{{a}}{\Delta}\kappa_-},
\end{equation}
and 
\begin{equation}\label{fnm12}
f_{n,m+1} = -i2^{-2{a} -2}\frac{\eta^2{a}}{\pi}\Gamma^2\left(-\frac{{a}}{2}\right)
\Delta^{{a}}\frac{1 -\frac{{a}}{\Delta}\kappa_+}
{1 +\frac{{a}}{\Delta}\kappa_+},
\end{equation}
where
$$
\kappa_- = 
\frac{R^{(-)}_{11} -iR^{(-)}_{12}}
{R^{(-)}_{21} -iR^{(-)}_{22}},\quad \kappa_+ = 
\frac{R^{(-)}_{11} -R^{(-)}_{12}}
{R^{(-)}_{21} -R^{(-)}_{22}}.
$$
Using  estimate (\ref{Rpm}) for the matrix entries of  $R^{(-)}$ we see that
\begin{equation}\label{kappaest}
\kappa_{-} = i +O\left(\frac{1}{r}\right), \quad \kappa_{+} = -1 +O\left(\frac{1}{r}\right),
\quad r \to \infty.
\end{equation}
Simultaneously, we observe that
\begin{equation}\label{Delta1}
\Delta^{{a}} = \Delta^{{a}}_1\left(1 + \frac{2i{a}}{\Delta_1} + O\left(\frac{1}{r^2_1}\right)\right),
\end{equation}
and 
\begin{equation}\label{Delta2}
\Delta^{{a}} = \Delta^{{a}}_2\left(1 -\frac{2{a}}{\Delta_2} + O\left(\frac{1}{r^2_2}\right)\right),
\end{equation}
where we have introduced the notations,
$$
\Delta_1 := 2(m - i(n+1)),\quad \Delta_2 := 2(m+1 - in),
\quad r_1 = \sqrt{(n+1)^2 + m^2}, \quad \mbox{and}\quad r_2 = \sqrt{n^2 + (m+1)^2}.
$$
Equations (\ref{kappaest}), (\ref{Delta1}), and (\ref{Delta2}) imply that
$$
\Delta^{{a}}\frac{1 -\frac{{a}}{\Delta}\kappa_-}
{1 +\frac{{a}}{\Delta}\kappa_-} = \Delta^{{a}}_1\left(1 +  O\left(\frac{1}{r^2_1}\right)\right),
$$
and
$$
\Delta^{{a}}\frac{1 -\frac{{a}}{\Delta}\kappa_+}
{1 +\frac{{a}}{\Delta}\kappa_+} = \Delta^{{a}}_2\left(1 +  O\left(\frac{1}{r^2_2}\right)\right),
$$
Therefore, formulae (\ref{fn1m2}) and (\ref{fnm12})  generate the asymptotic equations,
\begin{equation}\label{fn1m3}
f_{n+1,m} = -i2^{-2{a} -2}\frac{\eta^2{a}}{\pi}\Gamma^2\left(-\frac{{a}}{2}\right)
\Delta^{{a}}_1\left(1 +O\left(\frac{1}{r^2_1}\right)\right), r\to \infty
\end{equation}
and
\begin{equation}\label{fnm13}
f_{n,m+1} = -i2^{-2{a} -2}\frac{\eta^2{a}}{\pi}\Gamma^2\left(-\frac{{a}}{2}\right)
\Delta^{{a}}_2\left(1 +O\left(\frac{1}{r^2_2}\right)\right), r\to \infty
\end{equation}
Comparing these equations with (\ref{Zgas1}), we immediately  conclude that
\begin{equation}\label{fn1m4}
f_{n+1,m} = \frac{\Gamma\left(1-\frac{{a}}{2}\right)}{\Gamma\left(1+\frac{{a}}{2}\right)}
\left(\frac{n+1+im}{2}\right)^{{a}}
\left(1 +O\left(\frac{1}{(n+1)^2 +m^2}\right)\right), r\to \infty
\end{equation}
and
\begin{equation}\label{fnm14}
f_{n,m+1} = \frac{\Gamma\left(1-\frac{{a}}{2}\right)}{\Gamma\left(1+\frac{{a}}{2}\right)}
\left(\frac{n+i(m+1)}{2}\right)^{{a}}
\left(1 +O\left(\frac{1}{n^2 + (m+1)^2}\right)\right), r\to \infty.
\end{equation}
 This proves Theorem \ref{theorem1} for an arbitrary parity of the value of the sum $n+m$.

\section{Discrete logarithm and Green's functions}\label{s.log}

\begin{figure}[t]
\begin{center}
\includegraphics[width=0.6\linewidth, trim=0 200 0 200]{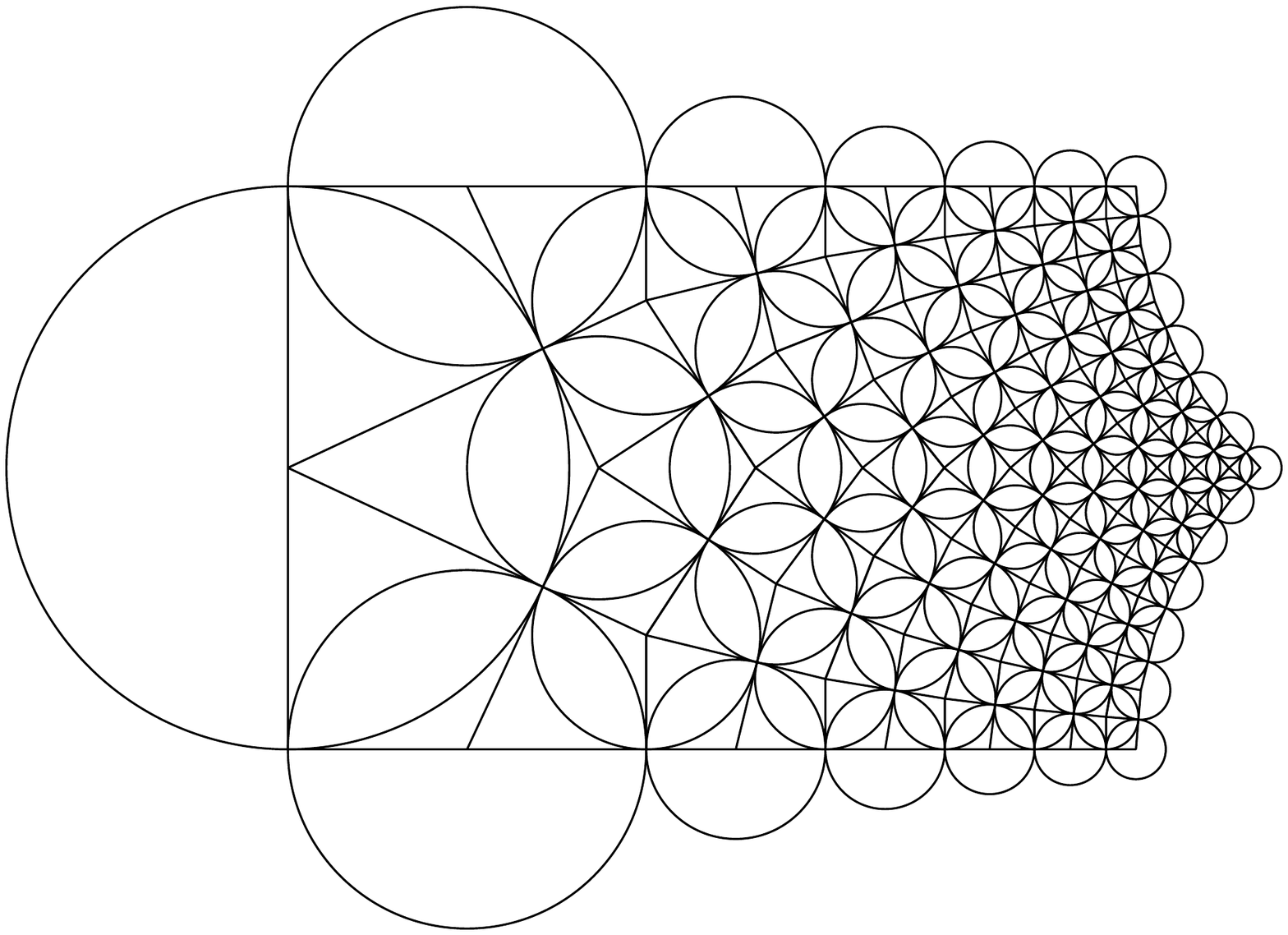} 
\end{center}
\caption{Discrete logarithm function $L(n,m)$ as an orthogonal circle pattern. {\small [Image by T.~Hoffmann]}}\label{f.log}
\end{figure}
Considered in this paper discrete function $Z^a$ with $0<a<2$ can be used to construct discrete analogs of logarithmic functions. The corresponding functions in the linear and nonlinear theories of discrete holomorphic functions were constructed in \cite{Ken} and \cite{AB} respectively. In this section we present the corresponding results and derive the asymptotics of these functions.

The circle pattern described by a discrete logarithm function $L(n,m)$ is presented in Figure \ref{f.log}.
As it was shown in \cite{AB} it can be obtained from the discrete $Z^a$ in the limit $a \to 0$ by the following formula (see \cite{AB, BMS, BobSurDDG} for more details):
\begin{equation}\label{L(n,m)}
L(n,m)= \lim_{a\to 0}\frac{Z^a (n,m)-1}{a}.
\end{equation}

There is another discrete version of the logarithmic function closely related to Green's function of the discrete Laplace operator on a isoradial graph, i.e. on a rhombic embedding of a quad-graph.
The system
\begin{eqnarray}
z_{n+1,m}-z_{n,m}=w_{n+1,m}w_{n,m},\quad z_{n,m+1}-z_{n,m}=iw_{n,m+1}w_{n,m}
\end{eqnarray}
describes a relation between solutions $z_{n,m}$ of the cross-ratio equation (\ref{def1}) and solutions $w_{n,m}$ of the Hirota equation
$$
w_{n,m}w_{n+1,m}+iw_{n+1,m}w_{n+1,m+1}-w_{n+1,m+1}w_{n,m+1}-iw_{n,m+1}w_{n,m}=0.
$$
The geometric meaning of the Hirota variables is the following: for even $n+m$ they are positive $w(n,m)\in{\mathbb R}_+$ and describe the radii of the corresponding circles, for odd $n+m$ they are unitary $w(n,m)\in S^1$ and describe the rotation angles at the intersection points of circles (see \cite{BMS, BobSurDDG} for details).
We denote by $W^a(n,m)$ the Hirota function corresponding to the the discrete $Z^a$, i.e. describing the radii and the rotation angles of the $Z^a$ circle pattern. Then as it was shown in \cite{BMS, BobSurDDG} the formula
\begin{equation}
\ell(n,m)=\frac{d}{da}W^{a-1}(n,m)_{|a=1}
\end{equation}
describes the discrete logarithm function in the linear theory. The latter satisfies the discrete Cauchy-Riemann equations
$$
\ell(n,m+1)-\ell(n+1,m)=i (\ell(n+1,m+1)-\ell(n,m)).
$$
At even $n+m$ this is Green's function of the discrete Laplace operator on an isoradial quad-graph introduced by
Kenyon \cite{Ken}.

\begin{theorem}\label{logas}
When $r^2 \equiv n^2+m^2 \to \infty$ the following asymptotic formulas hold for the nonlinear discrete logarithm (orthogonal circle pattern)
\begin{equation}\label{asympt_L}
L(n,m)=\log (n+im)+\gamma -\log 2 +O\left(\frac{\log r}{r^2}\right),
\end{equation}
and for the linear Green's function:
\begin{equation} \label{asympt_ell}
\ell(n,m)=\log \sqrt{n^2+m^2}+\gamma +\log 2 + O\left(\frac{\log r}{r}\right), \quad n+m \ {\rm even},
\end{equation}
where $\gamma$ is Euler's $\gamma$.
\end{theorem}
{\it Proof.}  The formal derivation of the asymptotic formulae (\ref{asympt_L}) and  (\ref{asympt_ell})
is easy.  Asymptotics (\ref{asympt_L}) is obtained by a direct differentiation of   estimate
(\ref{ABconj}) with respect to $a$ and putting then $a=0$. To obtain the second formula we observe the identity 
$W^a(n,m)=|Z^a_{n+1,m}-Z^a_{n,m}|$ at even $n+m$ due to the mentioned above geometric interpretation 
in terms of the radii of the circles.  After that the asymptotics (\ref{asympt_ell}) is a result of a simple computation 
including again the differentiation 
of estimate (\ref{ABconj}) with respect to $a$, this time at $a =1$. What  is needed  is the justification
of the legality of  differentiation of estimate (\ref{ABconj}). To this end it is enough to establish the following
two facts: (a) the validity of estimate (\ref{ABconj}) for the complex values of $a$ in the small neighborhoods
of the points $a=0$ and $a=1$ and (b) the analyticity of the map $Z^a$, at least for the large $n^2 + m^2$,
in these neighborhoods. In what follows we will show that these two facts indeed take place.

Applying  to the  Hankel asymptotic series the error term estimates (10.17.14) and (10.17.15) from \cite{DLMF}, one
can arrive at the following bound to  the error term in (\ref{Psispec})
%{\footnote{We remind that $H^{(1,2)}_{-\frac{a}{2}}(z)
%= \sqrt{\frac{2}{\pi z}}e^{\pm z}$.}}.
\begin{equation}\label{DLMF1}
\left|O\left(\frac{1}{\xi}\right)\right| \leq  \frac{\sqrt{2}\pi}{64|\xi|}|a^2 -1||a^2 -9|\exp\left\{ \frac{|a^2-1|}{4}\frac{\sqrt{2}\pi}{|\xi|}\right\}.
\end{equation}
This bound shows that the error term in  (\ref{Psispec}) and, as a consequence, the  error term in (\ref{match1}) 
are uniform in the small  complex neighborhoods of the points $a=0$ and $a=1$. This in turn implies the same uniformity
of the estimates (\ref{GRest2}) - (\ref{G2}) for the jump matrix $G_{R}(\lambda)$ of the $R$-RH problem.  In addition,
we notice that  estimate (\ref{GRest1}) is also uniform with respect to the complex $a$ in the indicated neighborhoods 
in view of the equation,
$$
\eta^{2}\omega^{-1}_{1,2} = e^{\pm \frac{\pi i a}{2}} \lambda^{\frac{a}{2}}.
$$
This means that the key estimate (\ref{Grest2}) is valid for the complex $a$ in the small neighborhoods of $a=0$ and $a=1$ with the universal
constant  $L$
and, as a consequence, that the final estimate (\ref{Ras})  for the solution $R(\lambda )$ of the $R$ - RH problem is
uniform in these neighborhoods.  This uniformity is obviously inherited by the estimates for $R_{jk}(0)$ given in 
Proposition \ref{R0str}. Let us notice that  
$$
\frac{q_{12}}{q_{22}} = -\frac{a}{\overline{\Delta}},\quad 
\frac{p_{21}}{p_{11}} = -\frac{p_{22}}{p_{12}} = -\frac{a}{\Delta}.
$$
Therefore, the estimates (\ref{D111}), (\ref{D121}) and, as a consequence, our final result - estimate (\ref{Zgas1}) for
the discrete map $Z^{a}$ are uniform in the small complex neighborhoods of $a=0$ and $a=1$. Let us now
show that the map $Z^{a}$ is analytic in these neighborhoods.

The analyticity of $Z^a$ with respect of  $a$, in fact its meromorphicity, is an immediate corollary of the formulae 
of  Section \ref{opol}. Indeed, equations (\ref{moments2}) shows that the moments $H_s$ are polynomials 
in $a$ and $e^{\frac{i\pi a}{2}}$; actually, they are linear functions in $e^{\frac{i\pi a}{2}}$ with  polynomial
in $a$ coefficients. In virtue of (\ref{Pldet}),  the polynomials $P_{l}(\lambda)$  are  meromorphic in $a$ and, in view of
(\ref{ZP}) so is the map $Z^a$. We only have to be sure that $a=0$ and $a=1$ are not, at least for sufficiently large
$n^2 +m^2$, its poles. This is true and follows from (\ref{Zgas1}). Together with the
uniformity of this estimate in $a$ in the small neighborhoods of $a=0$ and $a=1$ this allows us
to differentiate estimate (\ref{Zgas1}) which is equivalent to (\ref{ABconj})  with respect to $a$. The proof of Theorem \ref{logas} is completed.
\begin{remark} In order to be able to exploit   the analyticity - uniformity 
arguments for  justification of the differentiation with respect to $a$ when deriving  (\ref{asympt_ell}), one can
use the formula, 
$$
\frac{d}{da} W^{a}(n,m) = W^a(n,m)\Re \left(\frac{d}{da}\log \Bigl(Z^a_{n+1,m}-Z^a_{n,m}\Bigr)\right),
$$
which is valid for real $a$.
\end{remark}

Formula (\ref{asympt_ell}) is the asymptotics of the discrete Green function derived by Kenyon \cite{Ken}.

\section*{Acknowledgements}

This research was supported by the DFG Collaborative Research Center TRR 109, ``Discretization in Geometry and Dynamics''.
A. I.  also acknowledges support of the NSF grants DMS-1001777 and DMS-1361856, and of  SPbGU grant 11.38.215.2014.
  
\section{Appendix A. Proof of Proposition \ref{modanswer}}

The proof is formal: we will just check that the function $\Psi^{(0)}(\xi)$ determined by the right hand
side of the formula (\ref{Psi0form}) solves the Riemann-Hilbert problem (\ref{Psijamp}) - (\ref{Psizero}).

First we check the jump relations. The correct  jumps across  the rays $\Gamma_1$ and $\Gamma_2$
follows immediately from the  definition (\ref{Psi0form}). Consider then the jump across the ray $\Gamma_0$.
We have,
$$
\Psi^{(0)}_+(\xi) = \frac{\sqrt{\pi}}{2}\begin{pmatrix}\frac{1}{2}&0\cr\cr
0&2\xi\end{pmatrix}
\begin{pmatrix}H^{(2)}_{-{a}/2}\left(\frac{i}{2}\sqrt{\xi}_+\right)&H^{(1)}_{-{a}/2}\left(\frac{i}{2}\sqrt{\xi}_+\right)\cr\cr
\frac{d}{d\xi}H^{(2)}_{-{a}/2}\left(\frac{i}{2}\sqrt{\xi}_+\right)&\frac{d}{d\xi}H^{(1)}_{-{a}/2}\left(\frac{i}{2}\sqrt{\xi}_+\right)\end{pmatrix}
e^{\frac{\pi i{a}}{4}\sigma_3}
$$
and 
$$
\Psi^{(0)}_-(\xi) = \frac{\sqrt{\pi}}{2}\begin{pmatrix}\frac{1}{2}&0\cr\cr
0&2\xi\end{pmatrix}
\begin{pmatrix}H^{(2)}_{-{a}/2}\left(\frac{i}{2}\sqrt{\xi}_-\right)&H^{(1)}_{-{a}/2}\left(\frac{i}{2}\sqrt{\xi}_-\right)\cr\cr
\frac{d}{d\xi}H^{(2)}_{-{a}/2}\left(\frac{i}{2}\sqrt{\xi}_-\right)&\frac{d}{d\xi}H^{(1)}_{-{a}/2}\left(\frac{i}{2}\sqrt{\xi}_-\right)\end{pmatrix}
e^{\frac{\pi i{a}}{4}\sigma_3}
$$
$$
\times 
\begin{pmatrix}1 &0\cr
                                     2\cos\frac{\pi{a}}{2} & 1\end{pmatrix}
$$
The ray $\Gamma_0$ is the cut for all the multivalued functions involved. In particular,
$$
\sqrt{\xi}_- = \sqrt{\xi}_+ e^{i\pi}.
$$
The Hankel functions $H^{(1,2)}_{\nu}(z)$ are defined on the universal covering of ${\Bbb C}\setminus\{0\}$
and satisfy there the  relations (see e.g. \cite{BE}),
\begin{equation}\label{hankel1}
H^{(1)}_{\nu}(ze^{i\pi}) =-e^{-i\pi\nu}H^{(2)}_{\nu}(z),\quad 
H^{(1)}_{\nu}(ze^{-i\pi}) = 2\cos\pi\nu H^{(1)}_{\nu}(z) + e^{-i\pi\nu}H^{(2)}_{\nu}(z),
\end{equation}
\begin{equation}\label{hankel2}
H^{(2)}_{\nu}(ze^{-i\pi}) =-e^{i\pi\nu}H^{(1)}_{\nu}(z),\quad 
H^{(2)}_{\nu}(ze^{i\pi}) = 2\cos\pi\nu H^{(2)}_{\nu}(z) + e^{i\pi\nu}H^{(1)}_{\nu}(z).
\end{equation}
Therefore, 
$$
H^{(2)}_{-{a}/2}\left(\frac{i}{2}\sqrt{\xi}_-\right) = 2\cos\pi\frac{\pi{a}}{2} H^{(2)}_{-{a}/2}\left(\frac{i}{2}\sqrt{\xi}_+\right) + 
e^{-\frac{i\pi{a}}{2}}H^{(1)}_{-{a}/2}\left(\frac{i}{2}\sqrt{\xi}_+\right),
$$
$$
H^{(1)}_{-{a}/2}\left(\frac{i}{2}\sqrt{\xi}_-\right) =
-e^{\frac{i\pi{a}}{2}}H^{(2)}_{-{a}/2}\left(\frac{i}{2}\sqrt{\xi}_+\right),
$$
and the above formula for $\Psi^{(0)}_-(\xi) $ can be rewritten as,
$$
\Psi^{(0)}_-(\xi) = \frac{\sqrt{\pi}}{2}\begin{pmatrix}\frac{1}{2}&0\cr\cr
0&2\xi\end{pmatrix}
\begin{pmatrix}e^{-\frac{i\pi{a}}{2}}H^{(1)}_{-{a}/2}\left(\frac{i}{2}\sqrt{\xi}_+\right)&
-e^{\frac{i\pi{a}}{2}}H^{(2)}_{-{a}/2}\left(\frac{i}{2}\sqrt{\xi}_+\right)\cr\cr
e^{-\frac{i\pi{a}}{2}}\frac{d}{d\xi}H^{(1)}_{-{a}/2}\left(\frac{i}{2}\sqrt{\xi}_+\right)&
-e^{\frac{i\pi{a}}{2}}\frac{d}{d\xi}H^{(2)}_{-{a}/2}\left(\frac{i}{2}\sqrt{\xi}_+\right)\end{pmatrix}
\begin{pmatrix}1&0\cr\cr
-2e^{-\frac{i\pi{a}}{2}}\cos\frac{\pi{a}}{2}&1\end{pmatrix}
$$
$$
\times e^{\frac{\pi i{a}}{4}\sigma_3}
\begin{pmatrix}1 &0\cr
                                     2\cos\frac{\pi{a}}{2} & 1\end{pmatrix},
$$

\noindent
or

$$
\Psi^{(0)}_-(\xi) = \frac{\sqrt{\pi}}{2}\begin{pmatrix}\frac{1}{2}&0\cr\cr
0&2\xi\end{pmatrix}
\begin{pmatrix}e^{-\frac{i\pi{a}}{2}}H^{(1)}_{-{a}/2}\left(\frac{i}{2}\sqrt{\xi}_+\right)&
-e^{\frac{i\pi{a}}{2}}H^{(2)}_{-{a}/2}\left(\frac{i}{2}\sqrt{\xi}_+\right)\cr\cr
e^{-\frac{i\pi{a}}{2}}\frac{d}{d\xi}H^{(1)}_{-{a}/2}\left(\frac{i}{2}\sqrt{\xi}_+\right)&
-e^{\frac{i\pi{a}}{2}}\frac{d}{d\xi}H^{(2)}_{-{a}/2}\left(\frac{i}{2}\sqrt{\xi}_+\right)\end{pmatrix}
e^{\frac{\pi i{a}}{4}\sigma_3}
$$
$$
\times e^{-\frac{\pi i{a}}{4}\sigma_3}\begin{pmatrix}1&0\cr\cr
-2e^{-\frac{i\pi{a}}{2}}\cos\frac{\pi{a}}{2}&1\end{pmatrix} e^{\frac{\pi i{a}}{4}\sigma_3}
\begin{pmatrix}1 &0\cr
                                     2\cos\frac{\pi{a}}{2} & 1\end{pmatrix}
$$

$$
= \frac{\sqrt{\pi}}{2}\begin{pmatrix}\frac{1}{2}&0\cr\cr
0&2\xi\end{pmatrix}
\begin{pmatrix}e^{-\frac{i\pi{a}}{2}}H^{(1)}_{-{a}/2}\left(\frac{i}{2}\sqrt{\xi}_+\right)&
-e^{\frac{i\pi{a}}{2}}H^{(2)}_{-{a}/2}\left(\frac{i}{2}\sqrt{\xi}_+\right)\cr\cr
e^{-\frac{i\pi{a}}{2}}\frac{d}{d\xi}H^{(1)}_{-{a}/2}\left(\frac{i}{2}\sqrt{\xi}_+\right)&
-e^{\frac{i\pi{a}}{2}}\frac{d}{d\xi}H^{(2)}_{-{a}/2}\left(\frac{i}{2}\sqrt{\xi}_+\right)\end{pmatrix}
e^{\frac{\pi i{a}}{4}\sigma_3}
$$
$$
\times \begin{pmatrix}1&0\cr
-2\cos\frac{\pi{a}}{2}&1\end{pmatrix}
\begin{pmatrix}1 &0\cr
                                     2\cos\frac{\pi{a}}{2} & 1\end{pmatrix}
$$

\begin{equation}\label{jumpG0}
= \frac{\sqrt{\pi}}{2}\begin{pmatrix}\frac{1}{2}&0\cr\cr
0&2\xi\end{pmatrix}
\begin{pmatrix}e^{-\frac{i\pi{a}}{2}}H^{(1)}_{-{a}/2}\left(\frac{i}{2}\sqrt{\xi}_+\right)&
-e^{\frac{i\pi{a}}{2}}H^{(2)}_{-{a}/2}\left(\frac{i}{2}\sqrt{\xi}_+\right)\cr\cr
e^{-\frac{i\pi{a}}{2}}\frac{d}{d\xi}H^{(1)}_{-{a}/2}\left(\frac{i}{2}\sqrt{\xi}_+\right)&
-e^{\frac{i\pi{a}}{2}}\frac{d}{d\xi}H^{(2)}_{-{a}/2}\left(\frac{i}{2}\sqrt{\xi}_+\right)\end{pmatrix}
e^{\frac{\pi i{a}}{4}\sigma_3}
\end{equation}
>From (\ref{jumpG0}) it follows that,
$$
\Psi^{(0)}_-(\xi)\begin{pmatrix}0&1\cr\cr
-1&0\end{pmatrix}
$$
$$
= \frac{\sqrt{\pi}}{2}\begin{pmatrix}\frac{1}{2}&0\cr\cr
0&2\xi\end{pmatrix}
\begin{pmatrix}e^{-\frac{i\pi{a}}{2}}H^{(1)}_{-{a}/2}\left(\frac{i}{2}\sqrt{\xi}_+\right)&
-e^{\frac{i\pi{a}}{2}}H^{(2)}_{-{a}/2}\left(\frac{i}{2}\sqrt{\xi}_+\right)\cr\cr
e^{-\frac{i\pi{a}}{2}}\frac{d}{d\xi}H^{(1)}_{-{a}/2}\left(\frac{i}{2}\sqrt{\xi}_+\right)&
-e^{\frac{i\pi{a}}{2}}\frac{d}{d\xi}H^{(2)}_{-{a}/2}\left(\frac{i}{2}\sqrt{\xi}_+\right)\end{pmatrix}
\begin{pmatrix}0&1\cr\cr
-1&0\end{pmatrix}
e^{-\frac{\pi i{a}}{4}\sigma_3}
$$
$$
= \frac{\sqrt{\pi}}{2}\begin{pmatrix}\frac{1}{2}&0\cr\cr
0&2\xi\end{pmatrix}
\begin{pmatrix}e^{\frac{i\pi{a}}{2}}H^{(2)}_{-{a}/2}\left(\frac{i}{2}\sqrt{\xi}_+\right)&
e^{-\frac{i\pi{a}}{2}}H^{(1)}_{-{a}/2}\left(\frac{i}{2}\sqrt{\xi}_+\right)\cr\cr
e^{\frac{i\pi{a}}{2}}\frac{d}{d\xi}H^{(2)}_{-{a}/2}\left(\frac{i}{2}\sqrt{\xi}_+\right)&
e^{-\frac{i\pi{a}}{2}}\frac{d}{d\xi}H^{(1)}_{-{a}/2}\left(\frac{i}{2}\sqrt{\xi}_+\right)\end{pmatrix}
e^{-\frac{\pi i{a}}{4}\sigma_3}
$$
$$
= \frac{\sqrt{\pi}}{2}\begin{pmatrix}\frac{1}{2}&0\cr\cr
0&2\xi\end{pmatrix}
\begin{pmatrix}H^{(2)}_{-{a}/2}\left(\frac{i}{2}\sqrt{\xi}_+\right)&
H^{(1)}_{-{a}/2}\left(\frac{i}{2}\sqrt{\xi}_+\right)\cr\cr
\frac{d}{d\xi}H^{(2)}_{-{a}/2}\left(\frac{i}{2}\sqrt{\xi}_+\right)&
\frac{d}{d\xi}H^{(1)}_{-{a}/2}\left(\frac{i}{2}\sqrt{\xi}_+\right)\end{pmatrix}
e^{\frac{\pi i{a}}{4}\sigma_3}
= \Psi^{(0)}_+(\xi),\quad \xi \in \Gamma_0.
$$
Thus the function $\Psi^{(0)}(\xi)$ defined by (\ref{Psi0form}) satisfies all the prescribed 
jump condition. Next, we have to prove the asymptotics (\ref{Psiinfty}) and (\ref{Psizero}).
Consider (\ref{Psiinfty}) first.

The large $z$ behavior of the Hankel functions is given by the classical formulae (see \cite{BE} and \cite{DLMF}),
\begin{equation}\label{hankel3}
H_{\nu}^{(1)}(z) = \sqrt\frac{2}{\pi z} e^{i\left(z-\frac{\nu\pi}{2} -\frac{\pi}{4}\right)}
\left(1 + O\left(\frac{1}{z}\right)\right), 
\end{equation}
$$
z \to \infty, \quad -\pi < \arg z < 2\pi,
$$
and
\begin{equation}\label{hankel4}
H_{\nu}^{(2)}(z) = \sqrt\frac{2}{\pi z} e^{-i\left(z-\frac{\nu\pi}{2} -\frac{\pi}{4}\right)}
\left(1 + O\left(\frac{1}{z}\right)\right), 
\end{equation}
$$
z \to \infty, \quad -2\pi < \arg z < \pi.
$$
We remind that these  asymptotics are uniform in every sub-sector of the indicated sectors on the universal covering of
 ${\Bbb C}\setminus\{0\}$. We shall also assume that $i = e^{i\pi/2}$ in the all arguments of the Hankel  functions
 $H^{(1,2)}_{-{a}/2}\left(\frac{i}{2}\sqrt{\xi}\right)$. Consider the closed sector between the rays $\Gamma_0$ and $\Gamma_1$, i.e.,
 \begin{equation}\label{sect1}
 -\frac{\pi}{2} -2\theta \leq \arg \xi \leq \frac{\pi}{4} - 2\theta.
 \end{equation}
 For all $0\leq \theta \leq \pi/2$ we have that
 $$
  -\frac{3\pi}{2}  \leq \arg \xi \leq \frac{\pi}{4}, 
 $$
 and hence
 $$
  -\frac{3\pi}{4}  \leq \arg \sqrt {\xi} \leq \frac{\pi}{8}, 
 $$
while
$$
  -\frac{\pi}{4}  \leq \arg i\sqrt {\xi} \leq \frac{5\pi}{8}. 
$$
Therefore,  in sector (\ref{sect1}) and for all $\theta$ we can use for 
the functions $H^{(1,2)}_{-{a}/2}\left(\frac{i}{2}\sqrt{\xi}\right)$ 
formulae (\ref{hankel3}-\ref{hankel4}). This gives the following asymptotic
representations for these functions and their derivatives as 
$\xi \to \infty, \quad 
-\frac{\pi}{2} -2\theta \leq \arg \xi \leq \frac{\pi}{4} - 2\theta,$
\begin{equation}\label{as111}
H_{-{a}/2}^{(1)}\left(\frac{i}{2}\sqrt{\xi}\right) = -i\frac{2}{\sqrt{\pi}}e^{\frac{\pi i{a}}{4}}
\xi^{-\frac{1}{4}}e^{-\frac{1}{2}\sqrt{\xi}}\left(1+ O\left(\frac{1}{\sqrt{\xi}}\right)\right),
\end{equation}
\begin{equation}\label{as112}
\frac{d}{d\xi}H_{-{a}/2}^{(1)}\left(\frac{i}{2}\sqrt{\xi}\right) = \frac{i}{2\sqrt{\pi}}e^{\frac{\pi i{a}}{4}}
\xi^{-\frac{3}{4}}e^{-\frac{1}{2}\sqrt{\xi}}\left(1+ O\left(\frac{1}{\sqrt{\xi}}\right)\right),
\end{equation}
\begin{equation}\label{as113}
H_{-{a}/2}^{(2)}\left(\frac{i}{2}\sqrt{\xi}\right) = \frac{2}{\sqrt{\pi}}e^{-\frac{\pi i{a}}{4}}
\xi^{-\frac{1}{4}}e^{\frac{1}{2}\sqrt{\xi}}\left(1+ O\left(\frac{1}{\sqrt{\xi}}\right)\right),
\end{equation}
\begin{equation}\label{as114}
\frac{d}{d\xi}H_{-{a}/2}^{(2)}\left(\frac{i}{2}\sqrt{\xi}\right) = \frac{1}{2\sqrt{\pi}}e^{-\frac{\pi i{a}}{4}}
\xi^{-\frac{3}{4}}e^{\frac{1}{2}\sqrt{\xi}}\left(1+ O\left(\frac{1}{\sqrt{\xi}}\right)\right),
\end{equation}
In the same sector, the function
$\Psi^{(0)}(\xi)$ is given by the equation (cf. (\ref{Psi0form}),
\begin{equation}\label{as115}
\Psi^{(0)}(\xi) = \frac{\sqrt{\pi}}{2}\begin{pmatrix}\frac{1}{2}&0\cr\cr
0&2\xi\end{pmatrix}
\begin{pmatrix}H^{(2)}_{-{a}/2}\left(\frac{i}{2}\sqrt{\xi}\right)&H^{(1)}_{-{a}/2}\left(\frac{i}{2}\sqrt{\xi}\right)\cr\cr
\frac{d}{d\xi}H^{(2)}_{-{a}/2}\left(\frac{i}{2}\sqrt{\xi}\right)&\frac{d}{d\xi}H^{(1)}_{-{a}/2}\left(\frac{i}{2}\sqrt{\xi}\right)\end{pmatrix}
e^{\frac{\pi i{a}}{4}\sigma_3}.
\end{equation}
Combaining this formula with equations (\ref{as111}) - (\ref{as114}) we arrive at the desired large
$\xi$ behavior of $\Psi^{(0)}(\xi)$ in the sector (\ref{sect1}). Indeed, 
substituting (\ref{as111}) - (\ref{as114}) into the right hand side of (\ref{as115})and performing the trivial
matrix multiplications, we have,
$$
\Psi^{(0)}(\xi)=
\begin{pmatrix}\frac{1}{2}\xi^{-\frac{1}{4}}e^{\frac{1}{2}\sqrt{\xi}}\left(1+ O\left(\frac{1}{\sqrt{\xi}}\right)\right)
&-\frac{i}{2}\xi^{-\frac{1}{4}}e^{\frac{1}{2}\sqrt{\xi}}\left(1+ O\left(\frac{1}{\sqrt{\xi}}\right)\right)\cr\cr
\frac{1}{2}\xi^{\frac{1}{4}}e^{\frac{1}{2}\sqrt{\xi}}\left(1+ O\left(\frac{1}{\sqrt{\xi}}\right)\right)&
\frac{i}{2}\xi^{+\frac{1}{4}}e^{-\frac{1}{2}\sqrt{\xi}}\left(1+ O\left(\frac{1}{\sqrt{\xi}}\right)\right)
\end{pmatrix}
$$
\begin{equation}\label{assect1}
=\xi^{-\frac{1}{4}\sigma_3}\begin{pmatrix}\frac{1}{2}&-\frac{i}{2}\cr\cr
\frac{1}{2}&\frac{i}{2}\end{pmatrix}
\left(1+ O\left(\frac{1}{\sqrt{\xi}}\right)\right)e^{\frac{1}{2}\sqrt{\xi}\sigma_3},
\quad \xi \to \infty, \quad 
-\frac{\pi}{2} -2\theta \leq \arg \xi \leq \frac{\pi}{4} - 2\theta.
\end{equation}

Next we consider the sector between the rays $\Gamma_1$ and $\Gamma_2$, i.e.,
\begin{equation}\label{sect2}
 \frac{\pi}{4} -2\theta \leq \arg \xi \leq \frac{3\pi}{4} - 2\theta.
 \end{equation}
 For all $0\leq \theta \leq \pi/2$ we have that
 $$
  -\frac{3\pi}{4}  \leq \arg \xi \leq \frac{3\pi}{4}, 
 $$
 and hence
 \begin{equation}\label{ineq000}
  -\frac{3\pi}{8}  \leq \arg \sqrt {\xi} \leq \frac{3\pi}{8}, 
\end{equation}
while
$$
  \frac{\pi}{8}  \leq \arg i\sqrt {\xi} \leq \frac{7\pi}{8}. 
$$
Therefore,  in sector (\ref{sect2}) and for all $\theta$ we can again use for 
the functions $H^{(1,2)}_{-{a}/2}\left(\frac{i}{2}\sqrt{\xi}\right)$ 
formulae (\ref{hankel3}-\ref{hankel4}). The function $\Psi^{(0)}(\xi)$, however,
is now given by the equation (cf. (\ref{Psi0form}),
\begin{equation}\label{as120}
\Psi^{(0)}(\xi) = \frac{\sqrt{\pi}}{2}\begin{pmatrix}\frac{1}{2}&0\cr\cr
0&2\xi\end{pmatrix}
\begin{pmatrix}H^{(2)}_{-{a}/2}\left(\frac{i}{2}\sqrt{\xi}\right)&H^{(1)}_{-{a}/2}\left(\frac{i}{2}\sqrt{\xi}\right)\cr\cr
\frac{d}{d\xi}H^{(2)}_{-{a}/2}\left(\frac{i}{2}\sqrt{\xi}\right)&\frac{d}{d\xi}H^{(1)}_{-{a}/2}\left(\frac{i}{2}\sqrt{\xi}\right)\end{pmatrix}
e^{\frac{\pi i{a}}{4}\sigma_3}
\begin{pmatrix}1&0\cr\cr
e^{\frac{i\pi{a}}{2}}&1\end{pmatrix}.
\end{equation}
Therefore, instead of (\ref{assect1}, we shall get now,
$$
\Psi^{(0)}(\xi)=\xi^{-\frac{1}{4}\sigma_3}\begin{pmatrix}\frac{1}{2}&-\frac{i}{2}\cr\cr
\frac{1}{2}&\frac{i}{2}\end{pmatrix}
\left(1+ O\left(\frac{1}{\sqrt{\xi}}\right)\right)e^{\frac{1}{2}\sqrt{\xi}\sigma_3}
\begin{pmatrix}1&0\cr\cr
e^{\frac{i\pi{a}}{2}}&1\end{pmatrix}
$$
\begin{equation}\label{assect20}
= \xi^{-\frac{1}{4}\sigma_3}\begin{pmatrix}\frac{1}{2}&-\frac{i}{2}\cr\cr
\frac{1}{2}&\frac{i}{2}\end{pmatrix}
\left(1+ O\left(\frac{1}{\sqrt{\xi}}\right)\right)
\begin{pmatrix}1&0\cr\cr
e^{\frac{i\pi{a}}{2}}e^{-\sqrt{\xi}}&1\end{pmatrix}e^{\frac{1}{2}\sqrt{\xi}\sigma_3},
\end{equation}
$$
\xi \to \infty, \quad 
\frac{\pi}{4} -2\theta \leq \arg \xi \leq \frac{3\pi}{4} - 2\theta.
$$
At the same time, in the sector (\ref{sect2}) we have inequality (\ref{ineq000}).
Therefore, in the asymptotic formula (\ref{assect20}) the lower triangular
matrix in the right hand side can be droped, and we arrive at the desired large $\xi$ behavior
of the function $\Psi^{(0)}(\xi)$ in sector (\ref{sect2}).

Finally, we consider the sector between the rays $\Gamma_2$ and $\Gamma_0$, 
i.e.,
\begin{equation}\label{sect3}
 \frac{3\pi}{4} -2\theta \leq \arg \xi \leq \frac{3\pi}{2} - 2\theta.
 \end{equation}
This time,  for all $0\leq \theta \leq \pi/2$ we have that
 $$
  -\frac{\pi}{4}  \leq \arg \xi \leq \frac{3\pi}{2}, 
 $$
 
 \begin{equation}\label{ineq0001}
  -\frac{\pi}{8}  \leq \arg \sqrt {\xi} \leq \frac{\pi}{4}, 
\end{equation}
and
$$
  \frac{3\pi}{8}  \leq \arg i\sqrt {\xi} \leq \frac{5\pi}{4}. 
$$
This means that  we can continue to use asymptotic formula
(\ref{hankel3}) for the function $H^{(1)}_{-{a}/2}\left(\frac{i}{2}\sqrt{\xi}\right)$,
but can not use formula (\ref{hankel4}) for the function 
$H^{(2)}_{-{a}/2}\left(\frac{i}{2}\sqrt{\xi}\right)$. At the same time,
in sector (\ref{sect3}), the function $\Psi^{(0)}(\xi)$, 
is given by the equation (see again (\ref{Psi0form}),
$$
\Psi^{(0)}(\xi) = \frac{\sqrt{\pi}}{2}\begin{pmatrix}\frac{1}{2}&0\cr\cr
0&2\xi\end{pmatrix}
\begin{pmatrix}H^{(2)}_{-{a}/2}\left(\frac{i}{2}\sqrt{\xi}\right)&H^{(1)}_{-{a}/2}\left(\frac{i}{2}\sqrt{\xi}\right)\cr\cr
\frac{d}{d\xi}H^{(2)}_{-{a}/2}\left(\frac{i}{2}\sqrt{\xi}\right)&\frac{d}{d\xi}H^{(1)}_{-{a}/2}\left(\frac{i}{2}\sqrt{\xi}\right)\end{pmatrix}
e^{\frac{\pi i{a}}{4}\sigma_3}
\begin{pmatrix}1&0\cr\cr
2\cos{\frac{\pi{a}}{2}}&1\end{pmatrix}.
$$
or
$$
\Psi^{(0)}(\xi) = \frac{\sqrt{\pi}}{2}\begin{pmatrix}\frac{1}{2}&0\cr\cr
0&2\xi\end{pmatrix}
$$
\begin{equation}\label{assect300}
\times\begin{pmatrix}H^{(2)}_{-{a}/2}\left(\frac{i}{2}\sqrt{\xi}\right)
+2e^{-\frac{i\pi{a}}{2}}\cos{\frac{\pi{a}}{2}}H^{(1)}_{-{a}/2}\left(\frac{i}{2}\sqrt{\xi}\right)
&H^{(1)}_{-{a}/2}\left(\frac{i}{2}\sqrt{\xi}\right)\cr\cr
\frac{d}{d\xi}H^{(2)}_{-{a}/2}\left(\frac{i}{2}\sqrt{\xi}\right)
+2e^{-\frac{i\pi{a}}{2}}\cos{\frac{\pi{a}}{2}}\frac{d}{d\xi}H^{(1)}_{-{a}/2}
\left(\frac{i}{2}\sqrt{\xi}\right)
&\frac{d}{d\xi}H^{(1)}_{-{a}/2}\left(\frac{i}{2}\sqrt{\xi}\right)\end{pmatrix}
e^{\frac{\pi i{a}}{4}\sigma_3}.
\end{equation}
Observe now that from the second equation in (\ref{hankel1}) it follows that
$$
H^{(2)}_{-{a}/2}\left(\frac{i}{2}\sqrt{\xi}\right)
+2e^{-\frac{i\pi{a}}{2}}\cos{\frac{\pi{a}}{2}}H^{(1)}_{-{a}/2}\left(\frac{i}{2}\sqrt{\xi}\right)
= e^{-\frac{i\pi{a}}{2}}H^{(1)}_{-{a}/2}\left(\frac{i}{2}\sqrt{\xi}e^{-i\pi}\right). 
$$
Hence formula (\ref{assect300}) can be rewritten as,
$$
\Psi^{(0)}(\xi) = \frac{\sqrt{\pi}}{2}\begin{pmatrix}\frac{1}{2}&0\cr\cr
0&2\xi\end{pmatrix}
$$
\begin{equation}\label{assect301}
\times\begin{pmatrix}
e^{-\frac{i\pi{a}}{2}}H^{(1)}_{-{a}/2}\left(\frac{i}{2}\sqrt{\xi}
e^{-i\pi}\right)&H^{(1)}_{-{a}/2}\left(\frac{i}{2}\sqrt{\xi}\right)\cr\cr
e^{-\frac{i\pi{a}}{2}}\frac{d}{d\xi}H^{(1)}_{-{a}/2}
\left(\frac{i}{2}\sqrt{\xi}e^{-i\pi}\right)
&\frac{d}{d\xi}H^{(1)}_{-{a}/2}\left(\frac{i}{2}\sqrt{\xi}\right)\end{pmatrix}
e^{\frac{\pi i{a}}{4}\sigma_3}.
\end{equation}
We also observe that 
$$
 - \frac{5\pi}{8}  \leq \arg \left(i\sqrt {\xi}e^{-i\pi}\right) \leq \frac{\pi}{4}. 
$$
This means, that  we can use in the sector (\ref{sect3}) formula (\ref{hankel3})
for the  both $H^{(1)}_{-{a}/2}$ - functions in (\ref{assect301}), i.e., for 
the function $H^{(1)}_{-{a}/2}\left(\frac{i}{2}\sqrt{\xi}e^{-i\pi}\right)$ as well as for
the function  $H^{(1)}_{-{a}/2}\left(\frac{i}{2}\sqrt{\xi}\right)$. This gives us
in the sector (\ref{sect3}), in addition to (\ref{as111}) and (\ref{as112}), the asymptotic equations
(cf. (\ref{as113}) and (\ref{as114})), 
\begin{equation}\label{as213}
H_{-{a}/2}^{(1)}\left(\frac{i}{2}\sqrt{\xi}e^{-i\pi}\right) = \frac{2}{\sqrt{\pi}}e^{\frac{\pi i{a}}{4}}
\xi^{-\frac{1}{4}}e^{\frac{1}{2}\sqrt{\xi}}\left(1+ O\left(\frac{1}{\sqrt{\xi}}\right)\right),
\end{equation}
\begin{equation}\label{as214}
\frac{d}{d\xi}H_{-{a}/2}^{(1)}\left(\frac{i}{2}\sqrt{\xi}e^{-i\pi}\right) = \frac{1}{2\sqrt{\pi}}e^{\frac{\pi i{a}}{4}}
\xi^{-\frac{3}{4}}e^{\frac{1}{2}\sqrt{\xi}}\left(1+ O\left(\frac{1}{\sqrt{\xi}}\right)\right),
\end{equation}
as $\xi \to \infty$. Substituting these estimates, together with the estimates (\ref{as111})
and (\ref{as112}), into (\ref{assect301}) we obtain the desired large $\xi$ behavior of the
function $\Psi^{(0)}(\xi)$ in sector (\ref{sect3}). Indeed, we have that (cf. (\ref{assect1})),
$$
\Psi^{(0)}(\xi)=
\begin{pmatrix}\frac{1}{2}\xi^{-\frac{1}{4}}e^{\frac{1}{2}\sqrt{\xi}}\left(1+ O\left(\frac{1}{\sqrt{\xi}}\right)\right)
&-\frac{i}{2}\xi^{-\frac{1}{4}}e^{\frac{1}{2}\sqrt{\xi}}\left(1+ O\left(\frac{1}{\sqrt{\xi}}\right)\right)\cr\cr
\frac{1}{2}\xi^{\frac{1}{4}}e^{\frac{1}{2}\sqrt{\xi}}\left(1+ O\left(\frac{1}{\sqrt{\xi}}\right)\right)&
\frac{i}{2}\xi^{+\frac{1}{4}}e^{-\frac{1}{2}\sqrt{\xi}}\left(1+ O\left(\frac{1}{\sqrt{\xi}}\right)\right)
\end{pmatrix}
$$
\begin{equation}\label{assect3333}
=\xi^{-\frac{1}{4}\sigma_3}\begin{pmatrix}\frac{1}{2}&-\frac{i}{2}\cr\cr
\frac{1}{2}&\frac{i}{2}\end{pmatrix}
\left(1+ O\left(\frac{1}{\sqrt{\xi}}\right)\right)e^{\frac{1}{2}\sqrt{\xi}\sigma_3},
\quad \xi \to \infty, \quad 
\frac{3\pi}{4} -2\theta \leq \arg \xi \leq \frac{3\pi}{2} - 2\theta.
\end{equation}
This completes the proof of the fact that the function $\Psi^{(0)}(\xi)$ given
by the formula (\ref{Psi0form}) satisfies  the asymptotic condition (\ref{Psiinfty}).

Let us now show that the  asymptotic condition (\ref{Psiinfty}) can be actually written in the form
(\ref{Psispec}). To this end we need to calculate explicitly the term of order $1/\sqrt{\xi}$ in
(\ref{Psiinfty}). This term, as in fact the whole asymptitic series that can be written in the
right hand side of (\ref{Psiinfty}), does not depend on the sector in $\xi$-plane. Let us then choose
the sector $-\frac{\pi}{2} -2\theta <\arg\xi <  \frac{\pi}{4} -2\theta$ where the function $\Psi^{(0)}(\xi)$
is given by formula (\ref{as115}).  We will need the first  corrections to the asymptotic
equations (\ref{hankel3}) and (\ref{hankel4}). They are given by the formulae (see again \cite{BE} and \cite{DLMF} ), 
\begin{equation}\label{hankel3cor}
H_{\nu}^{(1)}(z) = \sqrt\frac{2}{\pi z} e^{i\left(z-\frac{\nu\pi}{2} -\frac{\pi}{4}\right)}
\left(1 + i\frac{4\nu^2 -1}{8z}+O\left(\frac{1}{z^2}\right)\right), 
\end{equation}
$$
z \to \infty, \quad -\pi < \arg z < 2\pi,
$$
and
\begin{equation}\label{hankel4cor}
H_{\nu}^{(2)}(z) = \sqrt\frac{2}{\pi z} e^{-i\left(z-\frac{\nu\pi}{2} -\frac{\pi}{4}\right)}
\left(1 - i\frac{4\nu^2 -1}{8z}+O\left(\frac{1}{z^2}\right)\right), 
\end{equation}
$$
z \to \infty, \quad -2\pi < \arg z < \pi,
$$
These formulae, in turn allow us to replace relations (\ref{as111}) - (\ref{as114}) by
the following more detail asymptotics,
\begin{equation}\label{as111cor}
H_{-{a}/2}^{(1)}\left(\frac{i}{2}\sqrt{\xi}\right) = -i\frac{2}{\sqrt{\pi}}e^{\frac{\pi i{a}}{4}}
\xi^{-\frac{1}{4}}e^{-\frac{1}{2}\sqrt{\xi}}\left(1-\frac{\psi_1}{\sqrt{\xi}} + O\left(\frac{1}{\xi}\right)\right),
\end{equation}
\begin{equation}\label{as112cor}
\frac{d}{d\xi}H_{-{a}/2}^{(1)}\left(\frac{i}{2}\sqrt{\xi}\right) = \frac{i}{2\sqrt{\pi}}e^{\frac{\pi i{a}}{4}}
\xi^{-\frac{3}{4}}e^{-\frac{1}{2}\sqrt{\xi}}\left(1- \frac{\psi_1-1}{\sqrt{\xi}}+ O\left(\frac{1}{\xi}\right)\right),
\end{equation}
\begin{equation}\label{as113cor}
H_{-{a}/2}^{(2)}\left(\frac{i}{2}\sqrt{\xi}\right) = \frac{2}{\sqrt{\pi}}e^{-\frac{\pi i{a}}{4}}
\xi^{-\frac{1}{4}}e^{\frac{1}{2}\sqrt{\xi}}\left(1+\frac{\psi_1}{\sqrt{\xi}}+O\left(\frac{1}{\xi}\right)\right),
\end{equation}
\begin{equation}\label{as114cor}
\frac{d}{d\xi}H_{-{a}/2}^{(2)}\left(\frac{i}{2}\sqrt{\xi}\right) = \frac{1}{2\sqrt{\pi}}e^{-\frac{\pi i{a}}{4}}
\xi^{-\frac{3}{4}}e^{\frac{1}{2}\sqrt{\xi}}\left(1+\frac{\psi_1-1}{\sqrt{\xi}}+O\left(\frac{1}{\xi}\right)\right),
\end{equation}
where {\footnote{ Actually the exact value of the coefficient $\psi_1$ is not that important.
What is important is  that it is the same value in all the eqautions  (\ref{as111cor}) - (\ref{as114cor}), and that
it appears in these equations where it appears.}}$\psi_1 =\frac{1-{a}^2}{4}$ (cf.\ref{Psi1}).
Substituting equations (\ref{as111cor}) - (\ref{as114cor}) into (\ref{as115}),  we have that  (cf. (\ref{assect1}))
$$
\Psi^{(0)}(\xi)=
\begin{pmatrix}\frac{1}{2}\xi^{-\frac{1}{4}}e^{\frac{1}{2}\sqrt{\xi}}\left(1 + \frac{\psi_1}{\sqrt{\xi}}+ O\left(\frac{1}{\xi}\right)\right)
&-\frac{i}{2}\xi^{-\frac{1}{4}}e^{\frac{1}{2}\sqrt{\xi}}\left(1-\frac{\psi_1}{\sqrt{\xi}}+ O\left(\frac{1}{\xi}\right)\right)\cr\cr
\frac{1}{2}\xi^{\frac{1}{4}}e^{\frac{1}{2}\sqrt{\xi}}\left(1+\frac{\psi_1-1}{\sqrt{\xi}} +O\left(\frac{1}{\xi}\right)\right)&
\frac{i}{2}\xi^{+\frac{1}{4}}e^{-\frac{1}{2}\sqrt{\xi}}\left(1-\frac{\psi_1-1}{\sqrt{\xi}} +O\left(\frac{1}{\xi}\right)\right)
\end{pmatrix}
$$

$$
=\xi^{-\frac{1}{4}\sigma_3}
\left[ \begin{pmatrix}\frac{1}{2}&-\frac{i}{2}\cr\cr
\frac{1}{2}&\frac{i}{2}\end{pmatrix} +\frac{1}{\sqrt{\xi}}
\begin{pmatrix}\frac{\psi_1}{2}&\frac{i}{2}\psi_1\cr\cr
\frac{\psi_1-1}{2}&-\frac{i}{2}(\psi_1 -1)\end{pmatrix} +O\left(\frac{1}{\xi}\right)\right]e^{\frac{1}{2}\sqrt{\xi}\sigma_3}
$$

$$
=\xi^{-\frac{1}{4}\sigma_3}
\left[ I +\frac{1}{\sqrt{\xi}}
\begin{pmatrix}\frac{\psi_1}{2}&\frac{i}{2}\psi_1\cr\cr
\frac{\psi_1-1}{2}&-\frac{i}{2}(\psi_1-1)\end{pmatrix}
\begin{pmatrix}1&1\cr\cr
i&-i\end{pmatrix} +O\left(\frac{1}{\xi}\right)\right]
\begin{pmatrix}\frac{1}{2}&-\frac{i}{2}\cr\cr
\frac{1}{2}&\frac{i}{2}\end{pmatrix}e^{\frac{1}{2}\sqrt{\xi}\sigma_3}
$$

$$
=\xi^{-\frac{1}{4}\sigma_3}
\left[ I +\frac{1}{\sqrt{\xi}}
\begin{pmatrix}0&\psi_1\cr\cr
\psi_1-1&0\end{pmatrix}
 +O\left(\frac{1}{\xi}\right)\right]
\begin{pmatrix}\frac{1}{2}&-\frac{i}{2}\cr\cr
\frac{1}{2}&\frac{i}{2}\end{pmatrix}e^{\frac{1}{2}\sqrt{\xi}\sigma_3},\quad \xi \to \infty,
$$
and this is the asymptotic equation (\ref{Psispec}).

To complete the proof of Proposition \ref{modanswer}, it is enough to notice that the
known expansions of the Hankel functions $H^{(1,2)}_{\nu}(z)$  at $z=0$ guarantee the
behavior indicated in (\ref{Psizero0}) and hence the asymptotic condition (\ref{Psizero}).
In fact, in Appendix B we derive representation (\ref{Psizero})
directly from  (\ref{Psi0form}) and  calculate explicitly  the relevant matrix $B_0$,
see equations (\ref{Psi0form4}) - (\ref{B0formula}) .

\section{Appendix B. Proof of Proposition \ref{P0pop}}
The proof of Proposition \ref{modanswer} is based on 
one the basic properties of the Bessel equation, which is
the possibility, rooted in the  relations (\ref{hankel1}) - (\ref{hankel2})
between the Hankel functions, to evaluated explicitely
the {\it Stokes multipiers } associated with the irregular singular
point $\lambda =\infty$. The proof of Proposition \ref{P0pop} exploits
another fundamental property of the Bessel equation, which is the possibility
to solve explicitely the {\it Connection Problem } associated with the two
singular points of the equation - the regular point at $\lambda =0$ and 
the  irregular point  at$\lambda =\infty$. This possibility is bassed 
on the classical relation between the Hankel and the Bessel functions (see again \cite{BE}),
\begin{equation}\label{HJ1}
H_{\nu}^{(1)}(z) = \frac{1}{i\sin\pi\nu}\left[J_{-\nu}(z) - J_{\nu}(z)e^{-i\nu\pi}\right],
\end{equation}
\begin{equation}\label{HJ2}
H_{\nu}^{(2)}(z) = \frac{1}{i\sin\pi\nu}\left[J_{\nu}(z)e^{i\nu\pi} - J_{-\nu}(z)\right],
\end{equation}

Using (\ref{HJ1}), (\ref{HJ2}), we can transform formula (\ref{Psi0form}) into the following representation of
the function  $\Psi^{(0)}(\xi)$ which is more suitable for the study of its behavior near $\xi =0$.
$$
\Psi^{(0)}(\xi) = \frac{\sqrt{\pi}}{2}\begin{pmatrix}\frac{1}{2}&0\cr\cr
0&2\xi\end{pmatrix}{\bf J}(\xi)
\begin{pmatrix}e^{-\frac{i\pi{a}}{2}}&-e^{\frac{i\pi{a}}{2}}\cr\cr
-1&1\end{pmatrix}
e^{\frac{\pi i{a}}{4}\sigma_3}\frac{i}{\sin{\frac{\pi{a}}{2}}}
$$
\begin{equation}\label{Psi0form1}
\times  \begin{cases}
I &-\frac{\pi}{2}-2\theta < \arg\xi < \frac{\pi}{4} -2\theta,\cr\cr
\begin{pmatrix}1 &0\cr
                                     2\cos\frac{\pi{a}}{2} & 1\end{pmatrix}  &\frac{3\pi}{4}-2\theta < \arg\xi < \frac{3\pi}{2} -2\theta,\cr\cr
\begin{pmatrix}1 &0\cr
                                     e^{\frac{i\pi{a}}{2}} & 1\end{pmatrix}  & \frac{\pi}{4}-2\theta < \arg\xi < \frac{3\pi}{4} -2\theta ,
\end{cases}
\end{equation}
where we denote,
\begin{equation}\label{bfJ}
{\bf J} (\xi) = \begin{pmatrix}J_{-{a}/2}\left(\frac{i}{2}\sqrt{\xi}\right)&J_{{a}/2}\left(\frac{i}{2}\sqrt{\xi}\right)\cr\cr
\frac{d}{d\xi}J_{-{a}/2}\left(\frac{i}{2}\sqrt{\xi}\right)&\frac{d}{d\xi}J_{{a}/2}\left(\frac{i}{2}\sqrt{\xi}\right)\end{pmatrix}.
\end{equation}
Observing that
$$
\begin{pmatrix}e^{-\frac{i\pi{a}}{2}}&-e^{\frac{i\pi{a}}{2}}\cr\cr
-1&1\end{pmatrix}
e^{\frac{\pi i{a}}{4}\sigma_3}\begin{pmatrix}1&0\cr\cr
e^{\frac{i\pi{a}}{2}}&1\end{pmatrix}
=\begin{pmatrix}-2i\sin\frac{\pi{a}}{2}e^{\frac{i\pi{a}}{4}}&-e^{\frac{i\pi{a}}{4}}\cr\cr
0&e^{-\frac{i\pi{a}}{4}}\end{pmatrix},
$$
equation (\ref{Psi0form1}) can be rewritten in the form,
$$
\Psi^{(0)}(\xi) = \frac{\sqrt{\pi}}{2}\begin{pmatrix}\frac{1}{2}&0\cr\cr
0&2\xi\end{pmatrix}{\bf J}(\xi)\frac{i}{\sin{\frac{\pi{a}}{2}}}
$$
\begin{equation}\label{Psi0form2}
\times
\begin{pmatrix}-2i\sin\frac{\pi{a}}{2}e^{\frac{i\pi{a}}{4}}&-e^{\frac{i\pi{a}}{4}}\cr\cr
0&e^{-\frac{i\pi{a}}{4}}\end{pmatrix}
 \begin{cases}
\begin{pmatrix}1 &0\cr
                                     -e^{\frac{i\pi{a}}{2}} & 1\end{pmatrix} &-\frac{\pi}{2}-2\theta < \arg\xi < \frac{\pi}{4} -2\theta,\cr\cr
\begin{pmatrix}1 &0\cr
                                 e^{-\frac{i\pi{a}}{2}}  & 1\end{pmatrix}  &\frac{3\pi}{4}-2\theta < \arg\xi < \frac{3\pi}{2} -2\theta,\cr\cr
I & \frac{\pi}{4}-2\theta < \arg\xi < \frac{3\pi}{4} -2\theta ,
\end{cases}
\end{equation}
Using the known convergent series  expansions of the Bessel function $J_{\nu}(z)$ at $z =0$ (see again \cite{BE}),
$$
J_{\nu}(z) = \sum_{j=0}^{\infty}\frac{(-1)^j}{j!\Gamma(\nu +j +1)}\left(\frac{z}{2}\right)^{2j+\nu},
$$
we derive from (\ref{Psi0form2})
the asymptotic representation  of $\Psi^{(0)}$ at $\xi = 0$. We have,
$$
\Psi^{(0)}(\xi) = \begin{pmatrix}2^{{a} -2}\frac{\sqrt{\pi}}{\Gamma\left(1 -\frac{{a}}{2}\right)}e^{-\frac{i\pi{a}}{4}}&
2^{-{a} -2}\frac{\sqrt{\pi}}{\Gamma\left(1 + \frac{{a}}{2}\right)}e^{\frac{i\pi{a}}{4}}\cr\cr
-2^{{a} -2}\frac{{a}\sqrt{\pi}}{\Gamma\left(1 -\frac{{a}}{2}\right)}e^{-\frac{i\pi{a}}{4}}&
2^{-{a} -2}\frac{{a}\sqrt{\pi}}{\Gamma\left(1 + \frac{{a}}{2}\right)}e^{\frac{i\pi{a}}{4}}\end{pmatrix}
\Bigl(I + O(\xi)\Bigr)\xi^{-\frac{{a}}{4}\sigma_3}
$$
\begin{equation}\label{Psi0form3}
\times
\frac{i}{\sin{\frac{\pi{a}}{2}}}\begin{pmatrix}-2i\sin\frac{\pi{a}}{2}e^{\frac{i\pi{a}}{4}}&-e^{\frac{i\pi{a}}{4}}\cr\cr
0&e^{-\frac{i\pi{a}}{4}}\end{pmatrix}
 \begin{cases}
\begin{pmatrix}1 &0\cr
                                     -e^{\frac{i\pi{a}}{2}} & 1\end{pmatrix} &-\frac{\pi}{2}-2\theta < \arg\xi < \frac{\pi}{4} -2\theta,\cr\cr
\begin{pmatrix}1 &0\cr
                                 e^{-\frac{i\pi{a}}{2}}  & 1\end{pmatrix}  &\frac{3\pi}{4}-2\theta < \arg\xi < \frac{3\pi}{2} -2\theta,\cr\cr
I & \frac{\pi}{4}-2\theta < \arg\xi < \frac{3\pi}{4} -2\theta ,
\end{cases}
\end{equation}
Noticing that 
$$
\frac{i}{\sin{\frac{\pi{a}}{2}}}\begin{pmatrix}-2i\sin\frac{\pi{a}}{2}e^{\frac{i\pi{a}}{4}}&-e^{\frac{i\pi{a}}{4}}\cr\cr
0&e^{-\frac{i\pi{a}}{4}}\end{pmatrix}
= \begin{pmatrix} 2e^{\frac{i\pi{a}}{4}}&0\cr\cr
0&\frac{i}{\sin\frac{\pi{a}}{2}}e^{-\frac{i\pi{a}}{4}}\end{pmatrix}
\begin{pmatrix}1&\frac{1}{2i\sin\frac{\pi{a}}{2}}\cr\cr
0&1\end{pmatrix},
$$
and taking into account some of the basic properties of the $\Gamma$ - function,
we conclude  from (\ref{Psi0form3}) that, as $\xi \to 0$,
$$
\Psi^{(0)}(\xi) = \begin{pmatrix}-2^{{a}}\frac{\sqrt{\pi}}{{a}\Gamma\left(-\frac{{a}}{2}\right)}&
-2^{-{a}-2}\frac{i}{\sqrt{\pi}}\Gamma\left(-\frac{{a}}{2}\right)\cr\cr
2^{{a}}\frac{\sqrt{\pi}}{\Gamma\left(-\frac{{a}}{2}\right)}&
-2^{-{a}-2}\frac{i{a}}{\sqrt{\pi}}\Gamma\left(-\frac{{a}}{2}\right)
\end{pmatrix}
\Bigl(I + O(\xi)\Bigr)\xi^{-\frac{{a}}{4}\sigma_3}
$$
\begin{equation}\label{Psi0form4}
\times
\begin{pmatrix}1&\frac{1}{2i\sin\frac{\pi{a}}{2}}\cr\cr
0&1\end{pmatrix}
 \begin{cases}
\begin{pmatrix}1 &0\cr
                                     -e^{\frac{i\pi{a}}{2}} & 1\end{pmatrix} &-\frac{\pi}{2}-2\theta < \arg\xi < \frac{\pi}{4} -2\theta,\cr\cr
\begin{pmatrix}1 &0\cr
                                 e^{-\frac{i\pi{a}}{2}}  & 1\end{pmatrix}  &\frac{3\pi}{4}-2\theta < \arg\xi < \frac{3\pi}{2} -2\theta,\cr\cr
I & \frac{\pi}{4}-2\theta < \arg\xi < \frac{3\pi}{4} -2\theta ,
\end{cases}
 \end{equation}
 
\begin{equation}\label{Psi0form5}
\equiv \begin{pmatrix}-2^{{a}}\frac{\sqrt{\pi}}{{a}\Gamma\left(-\frac{{a}}{2}\right)}&
-2^{-{a}-2}\frac{i}{\sqrt{\pi}}\Gamma\left(-\frac{{a}}{2}\right)\cr\cr
2^{{a}}\frac{\sqrt{\pi}}{\Gamma\left(-\frac{{a}}{2}\right)}&
-2^{-{a}-2}\frac{i{a}}{\sqrt{\pi}}\Gamma\left(-\frac{{a}}{2}\right)
\end{pmatrix}
\Bigl(I + O(\xi)\Bigr)\xi^{-\frac{{a}}{4}\sigma_3}C_0,
\end{equation}
where the matrix $C_0$ is the same as in (\ref{Psizero}).
The comparison of representation (\ref{Psi0form5}) with equation (\ref{Psizero}), yields the following 
explicit formula for the matrix $B_0$ in (\ref{Psizero}),
\begin{equation}\label{B0formula} 
B_0
= \begin{pmatrix}-2^{{a}}\frac{\sqrt{\pi}}{{a}\Gamma\left(-\frac{{a}}{2}\right)}&
-2^{-{a}-2}\frac{i}{\sqrt{\pi}}\Gamma\left(-\frac{{a}}{2}\right)\cr\cr
2^{{a}}\frac{\sqrt{\pi}}{\Gamma\left(-\frac{{a}}{2}\right)}&
-2^{-{a}-2}\frac{i{a}}{\sqrt{\pi}}\Gamma\left(-\frac{{a}}{2}\right)
\end{pmatrix},
\end{equation}
and the formula,
\begin{equation}\label{B0formula2} 
B
= \begin{pmatrix}-2^{{a}}\frac{\sqrt{\pi}}{\eta{a}\Gamma\left(-\frac{{a}}{2}\right)}&
-2^{-{a}-2}\frac{i\eta}{\sqrt{\pi}}\Gamma\left(-\frac{{a}}{2}\right)\cr\cr
2^{{a}}\frac{\sqrt{\pi}}{\eta\Gamma\left(-\frac{{a}}{2}\right)}&
-2^{-{a}-2}\frac{i\eta{a}}{\sqrt{\pi}}\Gamma\left(-\frac{{a}}{2}\right)
\end{pmatrix},
\end{equation}
for the $B$-matrix in (\ref{B0}).
We are now just one step from the formula for $\hat{P}^{(0)}(0)$.  Indeed, the definition
of  $P^{(0)}(\lambda)$ (see equation (\ref{P0def}) implies that
$$
\hat{P}^{(0)}(\lambda) = E(\lambda)\hat{\Phi}^{(0)}(\xi(\lambda))\left(\frac{\lambda}{\xi(\lambda)}\right)^{\frac{{a}}{4}\sigma_3}.
$$
Therefore,
\begin{equation}\label{hatP0for00}
\hat{P}^{(0)}(0) = \Delta^{\frac{1}{2}\sigma_3}B\Delta^{-\frac{{a}}{2}\sigma_3}, \quad 
\Delta = 2(m-in),
\end{equation}
and equations (\ref{hatP0for}) and (\ref{Deltadef}) follow from (\ref{B0formula2}).

\section{Appendix C. The ${a} =1$ case}\label{gamma1case}\label{apC}

As it has already been indicated in Remark \ref{gamma1remark}, in  
 the  case ${a} =1$
the unique solution of (\ref{def1}) - (\ref{incond}) is, as expected, $f_{n,m} = n+im \equiv Z^{{a}}|_{{a} =1}$.
Correspondingly, $u_{n,m} = 1$ and $v_{n,m} =i$ for all $n$ and $m$. This in turn
implies that, for all $n$ and $m$, the matrices $U_{n,m}$ and $V_{n,m}$  from the Lax pair (\ref{lax})
are given by the simple formulae,
\begin{equation}\label{UVgamma1}
U_{n,m}(\lambda)\equiv U(\lambda) = \begin{pmatrix}1&-1\cr\cr
\lambda&1\end{pmatrix},\quad 
V_{n,m}(\lambda)\equiv V(\lambda) = \begin{pmatrix}1&-i\cr\cr
i\lambda&1\end{pmatrix}.
\end{equation}
The corresponding function $\Psi_{n,m}(\lambda)$ is given by the equation (cf. (\ref{Psinmdef})),
\begin{equation}\label{Psigamma1}
\Psi_{n,m}(\lambda) = U^{n}V^{m}\lambda^{-\frac{1}{4}\sigma_3}.
\end{equation}
Matrices $U$ and $V$ are commute (as they should !) and their simultaneous diaganalization
can be written down as follows,
\begin{equation}\label{diaggamma1}
U(\lambda) = Q\begin{pmatrix}1-i\sqrt{\lambda}&0\cr\cr
0&1+i\sqrt{\lambda}\end{pmatrix}Q^{-1},\quad
V(\lambda) = Q\begin{pmatrix}1+\sqrt{\lambda}&0\cr\cr
0&1-\sqrt{\lambda}\end{pmatrix}Q^{-1},
\end{equation}
where
\begin{equation}\label{Qgamma1}
Q = \begin{pmatrix}1&1\cr\cr
i\sqrt{\lambda} &-i\sqrt{\lambda} \end{pmatrix}.
\end{equation}
Combining equations (\ref{diaggamma1}) with (\ref{Psigamma1}) we arrive at the
following explicit (i.e., no growing with $n$ and $m$ nontrivial matrix products) representation
of the function $\Psi_{n,m}$ in the case ${a} =1$.
$$
\Psi_{n,m}(\lambda) = Q\begin{pmatrix}(1-i\sqrt{\lambda})^n(1+\sqrt{\lambda})^m&0\cr\cr
0&(1+i\sqrt{\lambda})^n(1-\sqrt{\lambda})^m\end{pmatrix}Q^{-1}\lambda^{-\frac{1}{4}\sigma_3}
$$

\begin{equation}\label{Psigamma12}
 = Q
 \begin{pmatrix}1&0\cr\cr
 0&(1+\lambda)^n(1-\lambda)^m\end{pmatrix}
 e^{-\frac{i\pi}{2}n\sigma_3}e^{g(\lambda)\sigma_3}Q^{-1}
 \lambda^{-\frac{1}{4}\sigma_3},
 \end{equation}
 where 
\begin{equation}\label{g0}
g(\lambda) = m\log(1 + \sqrt{\lambda}) + n\log(i + \sqrt{\lambda}).
\end{equation}
 The corresponding 
 solution $Y(\lambda)$ of $Y$ - RH problem (\ref{Yjump}) - (\ref{Yzero}) is given
 by the equation,
 \begin{equation}\label{Ygamma1}
 Y(\lambda) = Q\begin{pmatrix}\frac{1}{2}&0\cr\cr
 \frac{1}{2}H^{-1}_0(\lambda)&\frac{i}{\sqrt{\lambda}}(-1)^m\end{pmatrix}
e^{-\frac{i\pi}{2}n\sigma_3} e^{g(\lambda)\sigma_3},
 \end{equation}
 where 
\begin{equation}\label{H0gamma1}
H_{0}(\lambda) =
 \left(\frac{1+\sqrt{\lambda}}{1-\sqrt{\lambda}}\right)^{m}
\left(\frac{i+\sqrt{\lambda}}{i-\sqrt{\lambda}}\right)^{n}.
\end{equation}
We note that $H^{-1}_{0}(\pm 1)  =0$ and hence the function (\ref{Ygamma1}), as it should,
has no singularities on ${\Bbb C}\setminus [0, -i\infty)$.

No  analog of the equations
(\ref{Psigamma12}) or  (\ref{Ygamma1}) is known for the generic non-commutative case ${a} \neq 1$.
However, as we have seen  in the  main text of the paper, the objects $g(\lambda)$ and $H_{0}(\lambda)$
which appear in the explicit formulae (\ref{Ygamma1}) for the solution of the 
Riemann-Hilbert problem (\ref{Yjump}) - (\ref{Yzero}) in the trivial ${a} =1$ case also play
central roles in the asymptotic analysis of this problem in the case of general ${a}$.
\begin{remark}\label{gammanot1}
For generic ${a}$, one can attempt to modify  ansatz (\ref{Psigamma1}) by replacing
the factor $\lambda^{-\frac{1}{4}\sigma_3}$ by the factor $\lambda^{-\frac{{a}}{4}\sigma_3}$.
This would lead to the replacement of equation  (\ref{Ygamma1}) by the equation,
\begin{equation}\label{Ygammagen}
Y(\lambda) = Q\begin{pmatrix}\frac{1}{2}&(-1)^{n}a_{-}H_0(\lambda)\cr\cr
\frac{1}{2}H^{-1}_0(\lambda)&(-1)^{m} a_{+}\end{pmatrix}
e^{-\frac{i\pi}{2}n\sigma_3} e^{g(\lambda)\sigma_3},
\end{equation}
wherek
$$
a_{\pm} = \frac{i}{2\sqrt{\lambda}}\left(\pm 1 -ie^{\frac{i\pi}{2}{a}}\lambda^{\frac{1-{a}}{2}}\right).
$$
The reader can easily check that the function $Y(\lambda)$ defined by this formula would
satisfy all the conditions of the Rimeann-Hilbert problem (\ref{Yjump}) - (\ref{Yzero}) except
it would not be analytic on ${\Bbb C}\setminus \Sigma_0$. Indeed, since $a_- $ is not zero
for all ${a}\neq 1$ the right hand side of (\ref{Ygammagen} )has poles at $ \lambda = \pm 1$.
\end{remark}

\end{document}